\documentclass[11pt]{amsart}
\usepackage{xcolor}

\usepackage{amsmath}
\usepackage{amssymb}
\usepackage{amsfonts}
\usepackage{amsthm}
\usepackage{enumerate}
\usepackage[mathscr]{eucal}
\usepackage{verbatim}
\usepackage{amsthm}
\usepackage{amscd}
\usepackage[mathscr]{eucal}
\usepackage{appendix}

\usepackage{multicol}

\theoremstyle{plain}
\newtheorem{theorem}{Theorem}

\newtheorem*{thmA}{Theorem A}
\newtheorem*{thmB}{Theorem B}
\newtheorem*{thmC}{Theorem C}

{}
\newtheorem{proposition}{Proposition}\numberwithin{proposition}{section}
\newtheorem{lemma}[proposition]{Lemma}
\numberwithin{corollary}{section}
\theoremstyle{definition}
\newtheorem{definition}{Definition}\numberwithin{definition}{section}
\theoremstyle{remark}
\newtheorem{remark}{Remark}

\numberwithin{equation}{section}

\newcommand{\N}{{\mathbb{N}}}
\newcommand{\n}{{\mathbb{R}^n}}
\newcommand{\bn}{{\mathbb{R}^{2n} }}

\newcommand{\one}{\uppercase\expandafter{\romannumeral1}}
\newcommand{\two}{\uppercase\expandafter{\romannumeral2}}
\newcommand{\three}{\uppercase\expandafter{\romannumeral3}}
\newcommand{\four}{\uppercase\expandafter{\romannumeral4}}
\newcommand{\spt}{\text{supp} \, }
\newcommand{\sch}{\mathcal{S}(\n)}
\newcommand{\z}{\mathbb{Z}}
\newcommand{\mA}{\mathcal{A}}

\newcommand{\SM}{\mathcal{SM}}
\newcommand{\MS}{\mathcal{MS}}
\newcommand{\DS}{\mathcal{SS}}
\newcommand{\RR}{\mathcal{R}}

\newcommand{\mmu}{\vec{\mu}}
\newcommand{\NN}{\mathcal{N}}

\newcommand{\tf}{\tfrac}
\newcommand{\vp}{\varphi}
\newcommand{\al}{\alpha}
\newcommand{\f}{\frac}
\newcommand{\sset}{\subset}

\newcommand{\q}{\quad}
\newcommand{\qq}{\qquad}

\newcommand{\bbz}{\mathbb{Z}}
\newcommand{\bbzn}{\mathbb Z^n}
\newcommand{\bbr}{\mathbb{R}}
\newcommand{\bbrn}{\mathbb R^n}

\newcommand{\bbn}{\mathbb{N}}

\def\000{\vec{\boldsymbol{0}}}

\newcommand{\supp}{\mathrm{supp}}

\newcommand{\wh}{\widehat}
\newcommand{\wt}{\widetilde}

\allowdisplaybreaks

\makeindex

\begin{document}

\author{Jiao Chen}
\address{School of Mathematical Sciences,
Chongqing Normal University, People's Republic of China}
\email{chenjiaobnu@163.com}

\author{Danqing He}
\address{School of Mathematical Sciences,
Fudan University, People's Republic of China}
\email{hedanqing@fudan.edu.cn}

\author{Guozhen Lu}
\address{Department of Mathematics, University of Connecticut, Storrs, CT 06269, USA}
\email{guozhen.lu@uconn.edu}

\author{Bae Jun Park}
\address{Department of Mathematics, Sungkyunkwan University, Suwon 16419, Republic of Korea}
\email{bpark43@skku.edu}

\author{Lu Zhang}
\address{School of Mathematics and Statistics,
Shaanxi Normal University, People's Republic of China}
\email{luzhang\_math@snnu.edu.cn}

\thanks{
D. He is supported by  National Key R$\&$D Program of China (No. 2021YFA1002500), the New Cornerstone Science Foundation, and Natural Science Foundation of Shanghai (No. 22ZR1404900 and 23QA1400300).
Jiao Chen's research is supported by the National Natural Science Foundation of China
(No. 11871126) and the Natural Science Foundation of Chongqing (No. CSTB2022NSCQ-JQX0004).
G. Lu's research is partly supported by a collaboration grant and Simons Fellowship from the Simons Foundation.
B. Park is supported in part by NRF grant 2022R1F1A1063637 and in part by POSCO Science Fellowship of POSCO TJ Park Foundation.
Lu Zhang's research is supported by the National Natural Science Foundation of China
(Grant No. 12271322), the Fundamental Research Funds for the Central Universities
(Grant No. GK202202007) and the Youth Innovation Team of Shaanxi Universities.
}

\title[Sharp H\"{o}rmander estimates for multi-parameter and  multi-linear multipliers]{ A sharp H\"{o}rmander estimate for multi-parameter and  multi-linear Fourier multiplier operators}

\maketitle

\begin{abstract}In this paper,
we investigate the H\"ormander type theorems for the multi-linear and multi-parameter Fourier multipliers. When the multipliers are characterized by $L^u$-based Sobolev norms for $1<u\le 2$ , our results on the smoothness assumptions  are sharp
 in the multi-parameter and bilinear case. In the multi-parameter and multi-linear case, our results are almost sharp. Moreover, even in the one-parameter and multi-linear case, our results improve earlier ones in the literature.
\end{abstract}

\section{Introduction}

Let $m\in L^\infty(\mathbb R^n)$ and we define the corresponding multiplier operator
$$
L_m(f)=\int_{\mathbb R^n}\wh f(\xi)m(\xi)e^{2\pi i \langle x,\xi\rangle}\;d\xi,$$
where $\wh f(\xi)=\int_{\mathbb R^n}f(x)e^{-2\pi i \langle x, \xi\rangle}\;dx$ is the Fourier transform of $f$. An important topic in harmonic analysis is to  characterize conditions on a Fourier multiplier $m$  for which $L_m$ is bounded on $L^p$, namely $\|L_m(f)\|_{L^p(\mathbb R^n)}\le C\|f\|_{L^p(\mathbb R^n)}$ for some absolute constant $C$.
The Mihlin multiplier theorem \cite{Mikhlin} says that if $m$ satisfies the differential estimates
\begin{equation}\label{Mcondition}
\big|\partial^{\alpha}m(\xi)\big|\lesssim_{\alpha}|\xi|^{-|\alpha|} \q \text{for all multi-indices $\alpha$ with $|\alpha|\le [n/2]+1$},
\end{equation}
then $L_m$ is bounded on $L^p(\bbrn)$ for all $1<p<\infty$.
An improvement of this result was provided by H\"ormander \cite{Hoe} who replaced the sufficient condition \eqref{Mcondition} by
\begin{equation}\label{Hcondition}
\sup_{j\in\mathbb{Z}}\big\Vert m(2^j \cdot) \,\wh{\vartheta}\big\Vert_{L_s^2(\bbrn)}<\infty\q \text{ for } ~s>n/2
\end{equation}
where $\wh{\vartheta}$ is an appropriate bump function on $\bbrn$ supported in an annulus of size $1$ and $L_s^2(\bbrn)$ is the Sobolev space on $\bbrn$. Weighted H\"{o}rmander type theorem was established by Kurtz \cite{Kurtz} and Kurtz and Wheeden \cite{KurtzWheeden}.
With such a H\"ormander condition, the $H^p$ boundedness  was later established by Calder\'on and Torchinsky \cite{CT}.  It was further improved by Baernstein and Sawyer \cite{BS} by considering the Besov regularity of the multipliers.

\medskip

In this paper, we consider a generalization of  H\"ormander's result in the multi-linear and multi-parameter setting. 
Coifman and Meyer \cite{CM} first studied a multi-linear version of Mihlin's result, and their result was later reinforced by Kenig and Stein \cite{kenig1999multilinear} and Grafakos and Torres \cite{grafakos2002multilinear}. If $m$ satisfies
$$|\partial^\al m(\xi)|\le C_\al |\xi|^{-|\al|}$$
for sufficiently large $L>0$ and all multi-indices $\alpha\in \mathbb{N}^{ln}$ with $|\alpha|\le L$, where $\xi=(\xi_1,\cdots,\xi_l)\in {(\bbrn)^l}$ then the corresponding multi-linear operator
\begin{equation}\label{multilinearmdef}
T_m\big(f_1,\cdots,f_l\big)(x):=\int_{(\bbrn)^l}\wh f_1(\xi_1)\cdots \wh f_l(\xi_l)m(\xi)e^{2\pi i\langle x,\xi_1+\cdots+\xi_l\rangle}\;d\xi
\end{equation}
is bounded from $L^{p_1}(\mathbb R^n)\times \cdots \times L^{p_l}(\mathbb R^n)$ to $L^p(\mathbb R^n)$ for $1<p_1,\cdots,p_l<\infty$ and $1/l<r<\infty$ with  $\tf1p=\tf1{p_1}+\cdots+\tf1{p_l}$. This so called  Coifman-Meyer theorem led to several useful tools such as the Bony decomposition \cite{Bony} and the Kato-Ponce inequality \cite{KP}. For flag trilinear Coifman-Meyer type theorem, see Muscalu \cite{Muscalu}.

A generalization of H\"ormander's multiplier theorem to the multi-linear setting, as far as we know, was initiated by Tomita \cite{Tom}.
Let $\vartheta^{(l)}$ be the $l$-linear analogue of $\vartheta$, that appeared in \eqref{Hcondition}, so that it is a Schwartz function on $\mathbb{R}^{ln}$ whose Fourier transform is supported in an annulus of constant size. Tomita showed that, in the multi-linear setting, if $s>ln/2$, then
$$
\big\|T_m(f_1,\cdots,f_l)\big \|_{L^p(\mathbb R^{ln})}
\lesssim \sup_{j\in\mathbb Z}\big\|m(2^j \cdot)\, \vartheta^{(l)}\big\|_{L_s^2(\mathbb R^{ln})} \|f_1\|_{L^{p_1}(\mathbb R^n)} \cdots\|f_l\|_{L^{p_l}(\mathbb R^n)}
$$
for $1<p_1,\cdots,p_l,p<\infty$ with $\tf1p=\tf1{p_1}+\cdots+\tf1{p_l}$.
Grafakos and Si \cite{GS} later improved and generalized the boundedness of $T_m$ to the range  $p<1$, using $L^u$-based Sobolev norms for $1<u\le 2$, instead of $L_s^2$ norm.
We also refer to \cite{FujitaTomita} for weighted multi-linear H\"{o}rmander multiplier theorems and \cite{Grafakos2013, Gr_Mi_Ng_Tom, Gr_Ng, MT}, which deal with several variant boundedness results using product type Sobolev norm conditions.

A comprehensive result including all of the previous results has recently been established in \cite{LHH}.

\begin{thmA}[{\cite{LHH}}]\label{t400}
Let $0<p_1,\cdots,p_l,p<\infty$, and $\tf1p=\tf1{p_1}+\cdots+\tf1{p_l}$. Suppose that
\begin{equation}\label{24condition}
s>{ln\over 2}, \quad  {1\over p}-{1\over 2}<{s\over n}+ \sum_{i\in I}({1\over p_i}-{1\over 2}), \quad j=1,\dots,d,
\end{equation}
where $I$ runs over all subsets of $J_l:=\{1,2,\dots,l\}$.  Then  there holds
\begin{equation}\label{24result}
\|T_m(f_1,\dots,f_l)\|_{L^p(\bbrn)}\lesssim  \sup_{j\in\mathbb Z}\big\|m(2^j \cdot)\, \vartheta^{(l)}\big\|_{L_s^2(\mathbb R^{ln})} \prod_{i=1}^l\|f_i\|_{H^{p_i}(\bbrn)}.
\end{equation}
\end{thmA}

Here, $H^p(\bbrn)$ denotes the Hardy space on $\bbrn$, which naturally extends $L^p(\bbrn)$ for $0<p\le 1$.
We note that the above theorem is almost sharp in the sense that  \eqref{24result} fails if at least one of inequalities in \eqref{24condition} has the opposite inequality sign.

\medskip

We now turn our attention to the multi-parameter theory, which is  more challenging than the  classical single-parameter theory of singular integrals.
The study of multi-parameter theory  goes back to \cite{JMZ} on the differentiability associated with the collection of rectangles, which was reproved by Cordoba and R. Fefferman \cite{CF}. Atomic decomposition theory of multi-parameter Hardy spaces was developed by Chang and R. Fefferman \cite{CHF} and $L^p$ boundedeness of multi-parameter singular operators were established  by R. Fefferman and Stein \cite{RFS} and Journ\'e \cite{J, J1, J2}. The boundedness on multi-parameter Hardy spaces was proved by R. Fefferman \cite{Fefferman-Annals} and Pipher \cite{P}. There is an extensive list of works on the development of multi-parameter theory in the past decades, we refer the interested reader to the many references   in the more recent  monographs by  Street  \cite{Street},  Nagel, Ricci,  Stein and  Wainger \cite{NRSW}, and  Shen and the first and fifth authors \cite{LSZ}.

 Regarding the multi-parameter mulitplier $L_m$ with the multi-parameter H\"ormander condition, their $L^p$ and $H^p$ estimates were studied in  Carbery and Seeger \cite{CS}, and the first author and third author \cite{CJL1}, which sharpened the classical Marcinkiewitz multiplier theorem (see \cite{Marc, Stein1}). Such results were then further improved by the first author, Huang and third author \cite{CHL1,CHL2} using the Littlewood-Paley-Stein square functions, atomic decomposition and Journ\'e's covering lemma \cite{J, P} by considering the multi-parameter Besov regularity. These results extend the works of Calder\'on and Torchinsky \cite{CT} and Baernstein and Sawyer \cite{BS} on H\"{o}rmander multiplier theorems on one-parameter Hardy spaces to the multi-parameter Hardy spaces.

For a bounded function $m$ on $\mathbb{R}^{4n}$, the associated bilinear bi-parameter operator $T_m$ is defined by
\begin{align}\label{op}
T_m\big(f,g \big)(x,y)&:=\int_{(\bbrn)^4}m(\xi_1,\xi_2,\eta_1,\eta_2)\wh{f}(\xi_1,\xi_2)\wh{g}(\eta_1,\eta_2)\\
&\qq\qq \times e^{2\pi i\langle x,\xi_1+\eta_1\rangle}e^{2\pi i \langle y,\xi_2+\eta_2\rangle}    \;d\xi_1d\xi_2 d\eta_1 d\eta_2, \qq x,y\in\bbrn\nonumber,
\end{align}
where $f$ and $g$ are Schwartz functions on $\bbr^{2n}$. A typical operator of this kind is the bi-parameter paraproduct that appeared naturally in  the Leibniz rule for fractional derivatives.
Regarding this operator,
Muscalu, Pipher, Tao and Thiele proved the following theorem  in the celebrated work \cite{MPTT1}.

\begin{thmB}[\cite{MPTT1}]
Let $m\in L^{\infty}(\bbr^{4})$ satisfies
$$\big|\partial_{\xi_1}^{\alpha_1} \partial_{\xi_2}^{\alpha_2} \partial_{\eta_1}^{\beta_1} \partial_{\eta_2}^{\beta_2} m(\xi_1,\xi_2,\eta_1,\eta_2)\big|\lesssim \frac{1}{(|\xi_1|+|\eta_1|)^{|\alpha_1|+|\beta_1|}}\frac{1}{(|\xi_2|+|\eta_2|)^{|\alpha_2|+|\beta_2|}}$$
for all indices $\alpha_1,\alpha_2,\beta_1,\beta_2>0$.
Then the bilinear operator $T_m$ is bounded from $L^{p}(\bbr^{2})\times L^{q}(\bbr^{2}) $ into $L^r(\bbr^{2})$ for $1<p,q\leq \infty$ and $0<r<\infty$ with ${1\over p}+{1\over q}={1\over r}.$
\end{thmB}
A main obstacle to handle the above bi-parameter setting in \cite{MPTT1} was the lack of bi-parameter Calder\'on-Zygmund theory, and the authors reduced the problem to bi-parameter paraproducts and applied techniques in time-frequency analysis. However, the argument in \cite{MPTT1} was based on Journ\'e's covering lemma for rectangles, which brought difficulty if one wants to extend it to general $d$-parameter cases for $d>2$.  Later on, the multi-parameter analogue was established in \cite{MPTT2}, where an elegant decomposition (see Lemma \ref{compact}) for the auxiliary functions were applied in order to avoid the dependence on Journ\'e's lemma. Moreover, to solve the difficulty when the target space $L^r$ is quasi-Banach ($0<r<1$), the authors in  \cite{MPTT1,MPTT2} took advantage of a weak-type duality argument and the multi-linear interpolation. An alternative proof of Theorem  B was given by the first and third authors using a different method in \cite{chen2014hormander} (see more detailed discussions below.)
More recently Li, Martikainen, and Vuorinen \cite{LMV} developed the bilinear bi-parameter Calder\'on-Zygmund theory, whose result implies Theorem~B. However, their theory cannot directly imply the bilinear bi-parameter version of H\"ormander's multiplier theorem we discuss below. Moreover, a bi-parameter trilinear flag-type multiplier was also studied by Muscalu and Zhai \cite{MZ}, the third author, Pipher and the fifth author \cite{LPZ}.

\medskip

In \cite{MPTT1,MPTT2}  the symbol $m$ is required to be differentiable sufficiently many times. It becomes a natural question to seek for the limited smoothness condition to obtain the boundedness of associated operators. To answer this problem completely (up to some endpoints) is the main purpose of this work.

To describe the problem accurately, we introduce some notations first.
Let $l,d\in\mathbb N$. Associated with an $l$-linear and $d$-parameter Fourier multiplier $m(\xi_1,\dots,\xi_l)$, where $\xi_i=(\xi_{i,1},\dots,\xi_{i,d})$,  
 $\xi_{i,j}\in \mathbb{R}^{n},1\leq i \leq l,1\leq j \leq d$,
 we define an operator
 \begin{gather*}
 \label{op1}
 T_m\big(f_1,\dots,f_l\big)(x):=\int_{\mathbb{R}^{dn}\times\dots\times \mathbb{R}^{dn}} m(\xi_1,\dots,\xi_l) \wh{ f_1} (\xi_1)\dots \wh{ f_l}(\xi_l) e^{2\pi i \langle x,\xi_1+\dots+\xi_l\rangle} \; d\xi_1 \cdots d\xi_l
 \end{gather*}
for $x\in \mathbb{R}^{dn}$, where $f_{i}\in \mathscr{S}(\mathbb{R}^{dn})$ for $1\leq i \leq l$.
For each $i=1,\dots,d$, let $\Theta_i$ be a Schwartz function on $\mathbb{R}^{ln}$, whose Fourier transform is supported in the  annulus $\big\{\vec {\xi^i}=(\xi_{1,i},\dots,\xi_{l,i})\in\mathbb{R}^{ln}: 1/2\le |\vec {\xi^i}|\le 2\big\}$ and denote by
$$\mA_{s_1,\dots,s_d}^{u}[m]:=\sup_{k_1,\dots,k_d\in\z}\big\Vert m(2^{k_1}\vec {\xi^1},\dots,2^{k_d} \vec {\xi^d})\wh{\Theta_1}(\vec {\xi^1})\cdots\wh{\Theta_d}(\vec {\xi^d})\big\Vert_{L^u_{s_1,\dots,s_d}(\mathbb{R}^{ln}\times\dots\times \mathbb{R}^{ln})},$$
where the $d$-parameter  Sobolev norm $\|\cdotp \|_{L^u_{s_1,\dots,s_d}(\mathbb{R}^{ln}\times\dots\times \mathbb{R}^{ln})}$ for $F$ on $\mathbb{R}^{ln}\times\dots\times \mathbb{R}^{ln}$ is defined as
\begin{equation*}
\| F \|_{L^u_{s_1,\dots,s_d}({\mathbb{R}^{ln}\times\dots\times \mathbb{R}^{ln}})}:=\Big\|\big[(1+4\pi^2 |\vec {\xi^1}|^2)^{s_1/ 2}\cdots  (1+4\pi^2 |\vec {\xi^d}|^2)^{s_d/ 2} \widehat{F} \big]^\vee\Big\|_{L^u(\mathbb{R}^{ln}\times\dots\times \mathbb{R}^{ln})}.
\end{equation*}

Concerning the boundedenss of $T_m$, the first and third author in \cite{chen2014hormander} proved the following result in the spirit of H\"{o}rmander multiplier theorem.
\begin{thmC}[\cite{chen2014hormander}]
Let $1<p_1,\dots, p_l<\infty$, and $0<p<\infty$ with ${1\over p_1}+\cdots{1\over {p_l}}={1\over p}$.
Suppose that
$$\min\{s_1,\dots, s_d\}>\tf{dn}2 \q \text{and}\q \min\{s_1,\dots, s_d\}>\frac{dn}{\min\{ p_1,\dots, p_l\}}.$$ Then the multi-linear multi-parameter operator $T_m$ satisfies
\begin{equation}
\big\|T_m(f_1,\dots, f_l)\big\|_{L^q(\bbr^{2n})}\lesssim  \mA_{s_1,\dots,s_d}^2[m] 
\|f_1\|_{L^{p_1}(\bbr^{dn})}\cdots \|f_l\|_{L^{p_l}(\bbr^{dn})}.
\end{equation}

\end{thmC}

Their proof takes advantage of the Littlewood-Paley-Stein theory and offers a more elementary and direct approach than using the deep time-frequency analysis given in \cite{MPTT1,MPTT2}. Moreover, it provides more relaxed smoothness assumption on the multipliers than in \cite{MPTT1, MPTT2}. This multi-parameter result was further extended to $1\leq u \leq 2$  by the second author, Grafakos, Nguyen, and Yan in \cite{HGY}. However, the restrictions for the smoothness indices $s_1,\dots,s_d$ in  \cite{chen2014hormander,HGY} are not sharp. In fact, this can be seen from the optimal smoothness condition established in \cite{LHH} for the single-parameter setting, which indicates that the very similar conclusion might be made in the multi-parameter setting.
For instance, \cite{LHH} shows that  in order for (single-parameter) bilinear  operator to have $L^p \times L^q \to L^r$ estimate for all $1<p,q<\infty$, we simply need $s_1,s_2>{3n\over 2}$, while both $s_1$ and $s_2$ are assumed to be greater than $2n$  in  \cite{chen2014hormander}.

\hfill

Our first result concerning the boundedness of multi-linear multi-parameter H\"ormander multipliers is
  given as follows.
\begin{theorem}   \label{ml}
	Let $1<u\leq 2$, $\frac{1}{u'}+\frac{1}{u}=1$, $\frac{1}{l}<p<\infty$, and $1<p_1,\dots, p_l< \infty$ with $\frac{1}{p}=\sum_{1\leq i \leq l}\frac{1}{p_i}$. Suppose that
	\begin{equation}
	\label{mlidx}
	s_j>{ln\over u}, \quad  {1\over p}-{1\over u'}<{s_j\over n}+ \sum_{i\in I}({1\over p_i}-{1\over u}), \quad j=1,\dots,d,
	\end{equation}
	where $I$ is an arbitrary subset of $J_l:=\{1,2,\dots,l\}$.  Then  there holds
	\begin{equation}
\label{mlmp}
  \|T_m(f_1,\dots,f_l)\|_{L^p(\bbrn)}\lesssim  \mathcal{A}_{s_1,\dots,s_d}^{u}[m] \prod_{i=1}^l\|f_i\|_{L^{p_i}(\bbrn)}.
\end{equation}
\end{theorem}

It seems that, even in the single-parameter setting,
Theorem   \ref{ml}  is new when $1<u<2$.
As we claimed above, this result
improves the single-parameter result for the multipliers characterized by $L^2$-based Sobolev norms in \cite{LHH}  to the multi-paramete and $L^u$-based ($1< u\leq 2$) setting, and relaxes the smoothness conditions for multi-linear multi-parameter multipliers in preceding works \cite{chen2014hormander,HGY}, which turns out to be almost optimal.

To avoid the notational complexity, we will discuss first the bilinear bi-parameter version, which contains all necessary ingredients to obtain the corresponding multi-linear multi-parameter result. For this reason, let us state the following result, while  the  general case will be discussed briefly in Section~\ref{mm}.

 \begin{theorem} \label{bl}
   Let $1<u\le2$,$\frac{1}{u'}+\frac{1}{u}=1$, $\frac{1}{2}<r<\infty$ and $1<p,q< \infty$ with $\frac{1}{r}=\frac{1}{p}+\frac{1}{q}$.
    Suppose that
    \begin{equation}\label{e001}
    s_1,s_2>\tf{2n}u, \q \text{ and }\q s_1,s_2>\frac{n}{r}-\frac{n}{u'}.
    \end{equation}
Then the bilinear operator $T_m$ satisfies
\begin{equation}\label{e002}
\big\|T_m(f,g)\big\|_{L^r(\bbr^{2n})}\lesssim  \mA_{s_1,s_2}^u[m] 
\|f\|_{L^p(\bbr^{2n})} \|g\|_{L^q(\bbr^{2n})}.
\end{equation}
\end{theorem}

The main novelty in our proof is that, in Lemma~\ref{lmnoncompactbip}, we provide a new $L^p$ estimate for hybrid operators  avoiding the compact Fourier support condition compared with the corresponding Lemma~3.2 of \cite{LHH}, which allows us to generalize their argument to the multi-parameter setting as the tail terms in \eqref{e021} are much easier to handle.
We recall that   the proof of \cite[Lemma~3.2]{LHH} relies on a Bernstein type inequality (see Lemma \ref{t301} below)  requiring a compact Fourier support condition, which brings a complicated argument in the single-parameter setting and
 an essential difficulty in the bi-parameter setting. In fact, in the bi-parameter case, to avoid using Journ\'e's lemma, we need the decomposition in Lemma~\ref{compact},  introducing smooth bump functions
compactly supported in the physical space, which are used in the definition of hybrid operators. These smooth bump functions play a crucial role in establishing Theorem~\ref{bl}.
To overcome the difficulty that  Bernstein type inequalities can no longer be directly applied, we decompose these bump functions further from the frequency side with the help of rapid decay of different smooth bump functions to obtain functions with compact Fourier supports. See \eqref{e0111} for the decomposition.

\hfill

The sharpness of the smoothness conditions on H\"ormander multipliers from linear to multi-linear settings were discussed extensively in the literature. For the linear version, one may consult \cite{GHHN, hirschman2, SL, W}. The papers \cite{Grafakos2016a, LHH} discussed counterexamples for multi-linear H\"ormander multipliers characterized by the classical Sobolev norm.
When the  classical Sobolev norm is replaced by the product type, we refer to \cite{GP, MT}.
In particular,  \cite{GP} obtained a sharp characterization by ruling out all endpoints. The sharpness of the smoothness on conditions on multi-parameter H\"ormander multipliers were discussed in \cite{CS, CHL1, CHL2, CJL1}.

Regarding the sharpness of the smoothness conditions  \eqref{mlidx} and \eqref{e001}, we obtain the  following theorem, from which one can conclude that Theorem \ref{bl} provides the optimal smoothness condition, and Theorem \ref{ml} gives an ``almost sharp" result.

\begin{theorem}\label{t329}
Let $1<p_1,\dots, p_l<\infty$, $\tf1p=\tf1{p_1}+\cdots+\tf1{p_l}$, $1<u<\infty$, and $0<s_1,\dots,s_d<\infty$.
Suppose that
\begin{equation}
\|T_m(f_1,\dots,f_l)\|_{L^p}\lesssim  \mathcal{A}_{s_1,\dots,s_d}^{u}[m] \prod_{i=1}^l\|f_i\|_{L^{p_i}},
\end{equation}
for all $m$ satisfying $\mathcal{A}_{s_1,\dots,s_d}^{u}[m]<\infty$.
Then it is necessary to have
\begin{equation}\label{e339}
\min\{\tf {s_1}n,\dots,\tf{s_d}n\}\ge \max\{\max_I\tf1p-\tf1{2}-\sum_{i\in I}(\tf1{p_i}-\tf12),\max_I\tf1p-\tf1{u'}-\sum_{i\in I}(\tf1{p_i}-\tf1u)\},
\end{equation}
and
\begin{equation}\label{e350}
\min\{\tf {s_1}n,\dots,\tf{s_d}n\}>\max\{\tf{l}u,\tf1p-\tf1{u'}\},
\end{equation}
where the maximum runs over all $I\sset\{1,2,\dots, l\}$.

\end{theorem}

\medskip

The necessary condition
$$
\min\{\tf {s_1}n,\dots,\tf{s_d}n\}\ge \max_I\tf1p-\tf1{u'}-\sum_{i\in I}(\tf1{p_i}-\tf1u)
$$
seems new in literature. One reason for why it was not noticed before, in our opinion, is that it is nontrivial only when $l\ge 3$, while the theory of $l$-linear H\"ormander multipliers with $l\ge 3$ is still far from being completed.

A key observation in establishing the necessary condition \eqref{e339} is that the Sobolev norm of the multiplier discussed in Theorem~\ref{t326} is closely related to the `codimension' of the support of the multiplier, which leads to the desired counterexample with the help of the argument from \cite{Grafakos2016a}. To establish the   condition  \eqref{e350}, we utilized the arguments in \cite{GP} that originated from an estimate for a variant of the Bessel potentials involving a logarithmic term in \cite{Gr_Park} and its discretized estimates in \cite{Park4}.

\medskip

This paper is organized as follows. In the next section, we provide  preliminaries for this paper, which mainly include certain maximal inequalities as well as the boundedness properties of  hybrid operators including  Max-Square, Square-Max, and double Square operators. These lemmas will be essential in handling  our decomposed pieces of the multiplier.
Section \ref{reductionsection} is devoted to a main reduction of \eqref{e002}. Specifically, we decompose $T_m(f,g)$ by adopting arguments of \cite{LHH} in the single-parameter setting, from which our problem can be reduced to Proposition \ref{pp42}. 
The proof of Proposition \ref{pp42} will be given in Section \ref{pop}.
 In Section \ref{mm}, we discuss the main ingredients so that we can extend  Theorem \ref{bl} to the multi-linear and multi-parameter setting, namely Theorem \ref{ml}. In Section \ref{ctexample}, we provide counterexamples to show Theorem \ref{t329}, and also discuss  certain endpoints {of the} condition \eqref{mlidx}.

\hfill

{\bf Notation}
We denote by $\bbn$, $\bbz$, and $\bbr$ the sets of natural numbers, integers, and real numbers, respectively, and let $\bbn_0:=\bbn\cup \{0\}$ and $\bbr_+:=\{r\in \bbr:r>0\}$.
We write $A\lesssim B$ if there exists an absolute constant $C$ such that $A\le CB$, and $A\approx B$ means that both $A\lesssim B$ and $B\lesssim A$ hold.
Moreover, $\mathscr{S}(\bbr^N)$ denotes the space of all Schwartz functions on $\bbr^N$ for $N\in\bbn$.

 \section{Preliminaries : Maximal inequalities}
 \label{prl}

Let $\mathcal{M}^{(n)}$ denote the Hardy-Littlewood maximal operator on $\bbr^n$, defined by
 $$\mathcal{M}^{(n)}F(x):=\sup_{x\in R} \frac{1}{|R|} \int_R |F(y)| ~ dy,\qq x\in\bbr^n$$
 for a locally integrable function $F$ on $\bbr^n$,
 where the supremum runs over all cubes $R$ in $\bbr^n$ containing $x$.
 For  $0<r<\infty$, we also define
\begin{equation*}
\mathcal{M}^{(n)}_{r}F(x):=\big(\mathcal{M}^{(n)}(|F|^r)(x)\big)^{1/ r}.
\end{equation*}
Then the maximal inequality
\begin{equation}\label{simplemaximalin}
\big\| \mathcal{M}_r^{(n)}F\big\|_{L^p(\bbr^n)} \lesssim   \|F\|_{L^p(\bbr^n)}
\end{equation}
holds for all $0<r<p\le \infty$.
In general, it turns out in \cite{FS} that for any $0<p<\infty$, $0<q\leq \infty$, and $0<r<p,q$,
\begin{equation}\label{multimaximalin}
   \big\| \big\{ \mathcal{M}_r^{(n)}F_j\big\}_{j\in \z}\big\|_{L^p(\ell^q)} \lesssim   \big\|\{F_j\}_{j\in \z}\big\|_{L^p(\ell^q)}.
\end{equation}

We also define the (product-type) strong maximal operators $\mathcal{M}^{(n),\otimes}$ acting on  $F\in \mathscr{S}(\bbr^{2n})$ by
$$\mathcal{M}^{(n),\otimes}F(x_1,x_2):=\sup_{x_1\in R_1}\sup_{x_2\in R_2}\frac{1}{|R_1||R_2|}\int_{R_1\times R_2}|F(y_1,y_2)|~dy_1dy_2 $$
and
$$\mathcal{M}_r^{(n),\otimes}F(x_1,x_2):=\big(\mathcal{M}^{(n),\otimes}\big(|F|^r\big)(x_1,x_2) \big)^{1/r}.$$
Then \eqref{simplemaximalin} implies that for $0<r<p\le \infty$
\begin{equation}\label{fsinq0}
\big\Vert \mathcal{M}_r^{(n),\otimes}F\big\Vert_{L^p(\bbr^{2n})}\lesssim \Vert F\Vert_{L^p(\bbr^{2n})}
\end{equation}
The following inequality is an immediate consequence of \eqref{multimaximalin}:\\
For  $0<p<\infty$, $0<q\leq \infty$, and $0<r<p,q$,
we have
\begin{equation}\label{fsinq}
   \big\| \big\{ \mathcal{M}_r^{(n),\otimes}F_j\big\}_{j\in \z}\big\|_{L^p(\ell^q)} \lesssim   \big\|\{F_j\}_{j\in \z}\big\|_{L^p(\ell^q)}.
\end{equation}
We remark that $$\mathcal{M}^{(2n)}F(x_1,x_2)\le \mathcal{M}^{(n),\otimes}F(x_1,x_2),$$
and the product-type maximal version $\mathcal{M}^{(n),\otimes}$ does not satisfy the weak type $(1,1)$ boundedness, while the original one $\mathcal{M}^{(2n)}$ does. See \cite[Exercise 2.1.6]{G249} for more details.

\hfill

The proof of Theorem \ref{bl} is based on decompositions via a (bi-parameter version of) Littlewood-Paley theory.
Let $\gamma$ be an integrable $\mathscr{C}^1$ function on $\bbr^n$ with $\int_{\bbr^n} \gamma(x)\, dx=0$ that satisfies
\begin{equation}\label{sizeest}
|\gamma(x)|+|\nabla \gamma(x)|\le B\frac{1}{(1+|x|)^{n+1}}
\end{equation} for some $B>0$. Define $\gamma_k:=2^{kn}\gamma(2^n\cdot )$.
Then for all  $1<p<\infty$,     there exists a constant $C_{n,p,\gamma}>0$ such that
\begin{equation}\label{LPtheory}
\big\Vert \big\{ \gamma_k\ast f\big\}_{k\in\bbz}\big\Vert_{L^p(\ell^2)}\le C_{n,p,\gamma}B\Vert f\Vert_{L^p(\bbr^n)}
\end{equation} for $f\in L^p(\bbr^n)$.
We refer to \cite[Theorem 5.1.2]{G249} for more details.
Our special interest will be given to a bi-parameter version of the above inequality as follows: \\
Let $\gamma^{(1)}$ and $\gamma^{(2)}$ be integrable $\mathscr{C}^1$ functions on $\bbr^n$ with the mean value zero satisfying \eqref{sizeest} with constants $B$ replaced by $B_1$ and $B_2$, respectively.
Then for $1<p<\infty$,  there exists a constant $C_{n,p,\gamma}'>0$ such that
\begin{equation}\label{biLPtheory}
\big\Vert \big\{ (\gamma^{(1)}_j\otimes\gamma^{(2)}_k)\ast F\big\}_{j,k\in\bbz}\big\Vert_{L^p(\ell^2)}\le C_{n,p,\gamma}' B_1 B_2\Vert F\Vert_{L^p(\bbr^{2n})}
\end{equation}
for $F\in L^p(\bbr^n\times \bbr^n)$.
This result follows from \eqref{LPtheory} by making use of Rademacher functions. We refer to \cite[Lemma~3.2]{muscalu2013classical} for more details.

\hfill

We now introduce some notations for balls and cubes in $\n$ that we will use frequently in the rest of this paper:\\
Suppose that $k\in \z,  l\in \mathbb{R}_+, m\in \z^n$.
  \begin{itemize}
     \item[(a)] Denote by $B_l(m)$ the ball in $\bbrn$ centered at $m$ with radius $l$. When the center is $0$, we simply write $B_l$.
     \item[(b)] Denote by $C_{l}(m):=m+[0,l]^n$ the cube with the lower left vertex $ m$ and the side length $ l$.
 For any cubes $I$ in $\bbrn$, let $\ell(I)$ denote the side-length of $I$ so that $\ell\big(C_{l}(m)\big)=l$.    All dyadic cubes in $\bbrn$ can be expressed as $C_{2^{-k}}(2^{-k}m)$ for some $k\in\bbz$ and $m\in\bbzn$ where $2^{-k}m$ and $2^{-k}$ would be its center and  side length. We simply denote $C_{2^{-k}}(2^{-k}m)$ by $I_k^m$ or $J_k^m$, and call
     $R_{j,k}^{\vec m}:=I^{m_1}_j\times J^{m_2}_k$ a dyadic rectangle for $j,k\in\bbz$ and ${\vec m}:= (m_1,m_2)\in (\bbzn)^2$.
      \item[(c)] For each $\vec M=(M_1,M_2)\in\mathbb \bbn_0\times \bbn_0$, 
      we denote by $R^{\vec m}_{j,k,\vec M}$ the dyadic rectangle $I^{m_1}_{j-M_1} \times J^{m_2}_{k-M_2}$.
     \item[(d)] Let $c>0$ be a positive constant. If $E\subset \n$ represents either  a ball  or a cube, then we use $c E$ to represent the set having the same center as $E$, but with radius or side length $c$ times. If $R=I\times J$ for cubes $I$ and $J$, then $cR:=cI \times cJ$.

         \item[(e)]
  We call $R_{-j,-k}^{0}=I^{0}_{-j}\times J^{0}_{-k}$ the dual rectangle of $R=R_{j,k}^{\vec m}=I^{m_1}_j\times J^{m_2}_k$, which is denoted by $R'$.

  \item[(f)] A weight $w_R(x_1,x_2)$ adapted to $R=R_{j,k}^{\vec m}=I^{m_1}_j\times J^{m_2}_k$ is defined by
  $$(1+2^{j}|x_1-2^{-j}m_1|)^{-N}(1+2^{k}|x_2-2^{-k}m_1|)^{-N}$$
  for a large enough $N$.

   \end{itemize}

\hfill

Then let us consider several auxiliary functions, so-called double maximal function, maximal-square function, square-maximal function, and double square function. For this one,
let $\varphi$ and $\psi$ represent radial Schwartz functions on $\n$ that satisfy
\begin{equation}\label{suppvarphipsi}
\spt \wh{\varphi}\subseteq \{\xi\in \n: |\xi|\lesssim 1\}, \qquad \spt \wh{\psi}\subseteq \{\xi\in \n: |\xi|\approx 1\}
\end{equation}
and let $\varphi_k:=2^{kn}\varphi(2^k\cdot)$ and $\psi_k:=2^{kn}\psi(2^k\cdot)$.
Functions satisfying these properties will be called $\Phi$-type and $\Psi$-type, respectively.
We use $\phi$ to represent a function that is either $\Phi$-type or $\Psi$-type function
so that $\phi_k$ indicates $\varphi_k$ or $\psi_k$.
We note that $\psi$ also satisfies all properties of $\gamma$ above and thus, Littlewood-Paley theories \eqref{LPtheory} and \eqref{biLPtheory} hold by replacing $\gamma$ by $\psi$.
That is, for $1<p<\infty$,
\begin{equation}\label{biLPtheory1}
\big\Vert \big\{ (\psi_j\otimes\psi_k)\ast F\big\}_{j,k\in\bbz}\big\Vert_{L^p(\ell^2)}\lesssim_{n,p,\psi} \Vert F\Vert_{L^p(\bbr^{2n})}.
\end{equation}
Moreover, for $1<p<\infty$, we have
\begin{align}
\big\Vert \big\{ \sup_{j\in\bbz}\big| (\varphi_j\otimes\psi_k)\ast F\big|\big\}_{k\in\bbz}\big\Vert_{L^p(\ell^2)}&\lesssim_{n,p,\varphi,\psi} \Vert F\Vert_{L^p(\bbr^{2n})},\label{redin1}\\
\big\Vert \big\{ \sup_{k\in\bbz}\big| (\psi_j\otimes\varphi_k)\ast F\big|\big\}_{j\in\bbz}\big\Vert_{L^p(\ell^2)}&\lesssim_{n,p,\varphi,\psi} \Vert F\Vert_{L^p(\bbr^{2n})},\nonumber\\
\Big\Vert \sup_{j\in\bbz}\big\Vert \big\{  (\varphi_j\otimes\psi_k)\ast F\big\}_{k\in\bbz}\big\Vert_{\ell^2}\Big\Vert_{L^p(\bbr^{2n})}&\lesssim_{n,p,\varphi,\psi} \Vert F\Vert_{L^p(\bbr^{2n})},\label{redin2}\\
\Big\Vert \sup_{k\in\bbz}\big\Vert \big\{  (\psi_j\otimes\varphi_k)\ast F\big\}_{j\in\bbz}\big\Vert_{\ell^2}\Big\Vert_{L^p(\bbr^{2n})}&\lesssim_{n,p,\varphi,\psi} \Vert F\Vert_{L^p(\bbr^{2n})},\nonumber\\
\Big\Vert \sup_{j,k\in\bbz}  \big| (\varphi_j\otimes\varphi_k)\ast F\big|\Big\Vert_{L^p(\bbr^{2n})}&\lesssim_{n,p,\varphi,\psi} \Vert F\Vert_{L^p(\bbr^{2n})}.\label{redin3}
\end{align}
The proofs of those estimates are standard. For example, by using \eqref{fsinq0}, one may easily obtain the last inequality \eqref{redin3}. By fixing two variables in order, with the help of the Fefferman-Stein inequality, \eqref{multimaximalin} yields \eqref{redin1}, and the remaining one \eqref{redin2} follows immediately from the pointwise estimate
$$\sup_{j\in\bbz}\big\Vert \big\{  (\varphi_j\otimes\psi_k)\ast F\big\}_{k\in\bbz}\big\Vert_{\ell^2}\le \big\Vert \big\{ \sup_{j\in\bbz}\big| (\varphi_j\otimes\psi_k)\ast F\big|\big\}_{k\in\bbz} \big\Vert_{\ell^2}. $$

\hfill

 We define hybrid square and maximal functions as follows:\\
For any $F\in \mathscr{S}(\n\times \n)$, fixed $\vec{M}=(M_1,M_2)\in (\bbn_0)^2$, $1\le u<\infty$, and $x\in \bbr^{2n}$,
\begin{equation}\label{hybridoperator}
\begin{aligned}
  \mathcal{SS}^{(u)}_{\vec{M}}(F)(x)&:=  \Big( \sum_{\substack{j,k\in \z \\ m_1,m_2\in \z^n}} \Big\|(\psi_j\otimes \psi_k) \ast F \Big({m_1+y_1 \over 2^{j-M_1}}, {m_2+y_2 \over 2^{k-M_2}}\Big)\Big\|_{L^u(|y|\lesssim 1)}^2 \chi_{R^{\vec m}_{j,k,\vec M}}(x)\Big)^{1/ 2} , \\
   \mathcal{MS}^{(u)}_{\vec{M}}(F)(x)&:=  \sup_{\substack{j \in \z\\ m_1 \in \z^n}} \Big( \sum_{\substack{k\in \z \\ m_2 \in \z^n}}   \Big\|\big(\varphi_j\otimes \psi_k\big) \ast F \Big({m_1+y_1 \over 2^{j-M_1}}, {m_2+y_2 \over 2^{k-M_2}}\Big)\Big\|_{L^u(|y|\lesssim 1)}^2 \chi_{R^{\vec m}_{j,k,\vec M}}(x) \Big)^{1/ 2},
  \\
  \mathcal{SM}^{(u)}_{\vec{M}}(F)(x)&:=  \Big( \sum_{\substack{j\in \z \\ m_1 \in \z^n}}  \sup_{\substack{k \in \z\\ m_2 \in \z^n}} \Big\|\big(\psi_j\otimes \varphi_k\big) \ast F \Big({m_1+y_1 \over 2^{j-M_1}}, {m_2+y_2 \over 2^{k-M_2}}\Big)\Big\|_{L^u(|y|\lesssim 1)}^2 \chi_{R^{\vec m}_{j,k,\vec M}}(x) \Big)^{1/ 2},
   \\
    \mathcal{MM}^{(u)}_{\vec{M}}(F)(x)&:=  \sup_{j,k\in \z} \sup_{ m_1,m_2 \in \n} \Big\|\big(\varphi_j\otimes \varphi_k\big) \ast F \Big({m_1+y_1 \over 2^{j-M_1}}, {m_2+y_2 \over 2^{k-M_2}}\Big)\Big\|_{L^u(|y|\lesssim 1)} \chi_{R^{\vec m}_{j,k,\vec M}}(x)
\end{aligned}
\end{equation}
where $y$ denotes $(y_1,y_2)\in(\bbrn)^2$.

\hfill

Then the next lemma is bi-parameter extensions of \cite[Lemma~3.1]{LHH}, which contain the boundedness of the above hybrid operators when $r=1$.
\begin{lemma} \label{lmmsl1}
Let  $1<p<\infty$ and  $M_1,M_2\in \bbn_0$.
The operators $ \mathcal{SS}^{(1)}_{\vec{M}},\  \mathcal{MS}^{(1)}_{\vec{M}},\  \mathcal{SM}^{(1)}_{\vec{M}},\  \mathcal{MM}^{(1)}_{\vec{M}}$ are bounded on $L^p(\n\times \n)$ uniformly in $M_1,M_2$. That is,
$$\big\Vert  \mathcal{SS}^{(1)}_{\vec{M}} F\big\Vert_{L^p(\bbr^{2n})},\, \big\Vert  \mathcal{MS}^{(1)}_{\vec{M}} F\big\Vert_{L^p(\bbr^{2n})},\, \big\Vert  \mathcal{SM}^{(1)}_{\vec{M}} F\big\Vert_{L^p(\bbr^{2n})},\, \big\Vert  \mathcal{MM}^{(1)}_{\vec{M}} F\big\Vert_{L^p(\bbr^{2n})} \lesssim \Vert F\Vert_{L^p(\bbr^{2n})}$$
uniformly in $M_1$ and $M_2$.
\end{lemma}

\begin{proof}
We first see that
\begin{align}\label{shiftlm1}
 &\Big\|\big(\phi_j\otimes \phi_k\big) \ast F \Big({m_1+y_1 \over 2^{j-M_1}}, {m_2+y_2 \over 2^{k-M_2}}\Big)\Big\|_{L^1(|y|\lesssim 1)}\nonumber\\
 &\lesssim \frac1{|R^{\vec m}_{j,k,\vec M}|}\int_{c R^{\vec m}_{j,k,\vec M}} \big| (\phi_j\otimes \phi_k) \ast F (z_1, z_2) \big| ~dz
\end{align}
and thus it follows that
 \begin{align}\label{maximalbound}
   &   \Big\|\big(\phi_j\otimes \phi_k\big) \ast F \Big({m_1+y_1 \over 2^{j-M_1}}, {m_2+y_2 \over 2^{k-M_2}}\Big)\Big\|_{L^1(|y|\lesssim 1)} \chi_{R^{\vec m}_{j,k,\vec M}}(x)  \nonumber\\
   &\lesssim  \mathcal{M}^{(n),\otimes} \big((\phi_j\otimes \phi_k) \ast F\big)(x)\chi_{R^{\vec m}_{j,k,\vec M}}(x)
 \end{align}
for some $c>0$, where the implicit constant is independent of $M_1,M_2$.
Then the maximal inequality \eqref{fsinq} and the double square estimate \eqref{biLPtheory1} clearly verify the estimate for $\mathcal{SS}^{(1)}_{\vec{M}}$.

Similarly, using \eqref{maximalbound}, \eqref{fsinq},  and \eqref{redin1}, we have
\begin{align*}
  \big\Vert  \mathcal{SM}^{(1)}_{\vec{M}} F\big\Vert_{L^p(\bbr^{2n})} &\lesssim  \bigg\|  \Big( \sum_{\substack{k\in \z}}      \big| \mathcal{M}^{(n),\otimes}\big( \sup_{\substack{j\in \z }}  |(\varphi_j\otimes \psi_k) \ast F|\big)\big|^2  \Big)^{1/2}\bigg\|_{L^p(\n\times \n)}\\
  &\lesssim  \big\|   \big\{ \sup_{\substack{j\in \z }}  |(\varphi_j\otimes \psi_k) \ast F|\big\}_{k\in\bbz} \big\|_{L^p(\ell^2)}\lesssim  \|F\|_{L^p(\n\times  \n)}.
\end{align*}

The proof of other two estimates is almost the same and so we omit it here.
\end{proof}

Considering the left-hand side of each estimate in the above lemma, if the inside $L^1$ norm  is replaced with an $L^q$ norm for $2\le q\le \infty$,
we need to lower the exponent of the maximal function  to apply the Fefferman-Stein inequality.
Due to the compact Fourier supports of $\varphi$ and $\psi$, we are able to fulfill this goal by making use of the Bernstein inequalities, at the cost of  the independence of the implicit constant   on $M_1$ and $M_2$.
\begin{lemma}\label{t301}
Let $0<p\le q\le\infty$.
Let $R=R_{j,k}^{\vec m}$, $R_{\vec t}=R_{j-t_1,k-t_2}^{\vec m}$, and $R_{\vec M}=R_{j,k,\vec M}^{\vec m}$.
Suppose that $\wh f$ is supported in the dual rectangle $R'_{\vec t}$ of $R_{\vec t}$. Then, for any $x\in R_{\vec M}$,
\begin{align}
\Big(\frac1{|R_{\vec M}|}\int_{R_{\vec M}}|f(y)|^qdy\Big)^{1/q} &\lesssim 2^{(M_1+M_2)n(1/p-1/q)} \max\{2^{-t_1n(1/p-1/q)},1\}\label{e301}\\
& \qq\times\max\{2^{-t_2n(1/p-1/q)},1\} \mathcal{M}_p^{(n),\otimes}  f(x)\chi_{R_{\vec M}}(x).\notag
\end{align}
\end{lemma}
A special case is contained in \cite[Lemma~3.2]{LHH}, compared to which our proof is slightly intuitive. Both proofs follow from the Bernstein inequalities.
\begin{proof}
We first consider the case when $t_1=t_2=0$.

By the classical Bernstein inequality (see, for instance, \cite[Proposition~5.5]{Wo}, \cite[p17]{Tri}, and references therein),
$\|f\|_{L^q(R)}\lesssim |R|^{1/q-1/p}\|f\|_{L^p(w_R)}$.
As $M_1,M_2\ge 0$, we can decompose $R_{\vec{M}}$ into rectangles $R'$ which are translations of $R$.
Moreover,
$$
\|f\|_{L^q(R')}\lesssim |R|^{1/q-1/p}\|f\|_{L^p(w_{R'})}
$$
holds for all $R'$.
This, combined with $p\le q$, implies further that
$$
\|f\|_{L^q(R_{\vec M})}
\lesssim |R|^{1/q-1/p}\|f\|_{L^p(w_{R_{\vec M}})}
$$
by summing over $R'$ as one easily observes $\sum_{R'\subset R_{\vec M}} w_{R'}\lesssim w_{R_{\vec M}}$. In particular, the implicit constant is independent of $\vec M$.

By definition we have
$w_{R_{\vec M}}(y)\lesssim w_{R^0_{\vec M}}(x-y)$
for any $x\in R_{\vec M}$, where $R^0_{\vec M}=R_{j,k,\vec M}^{\vec 0}$.
As a result, we can bound $\|f\|_{L^q(R_{\vec M})}$ by
\begin{align}
&|R|^{1/q-1/p}|R^0_{\vec M}|^{1/p}\Big(\frac1{|R^0_{\vec M}|}\int |f(y)|^pw_{R^0_{\vec M}}(x-y)dy\Big)^{1/p}\notag\\
\lesssim& |R_{\vec M}|^{1/q}2^{(M_1+M_2)n(1/p-1/q)} \mathcal{M}_p^{(n),\otimes}  f(x)\chi_{R_{\vec M}}(x).\label{e302}
\end{align}
This proves \eqref{e301} when $t_1=t_2=0$.

When $t_1,t_2\le 0$, the desired result follows from last case as
$$
R_{\vec M}=R_{j,k,\vec M}^{\vec m}=R_{j-t_1,k-t_2, (M_1-t_1,M_2-t_2)}^{(2^{-t_1}m_1,2^{-t_2}m_2)}.
$$

Next we consider the case when $t_1,t_2\ge 0$.
We control the left hand side of \eqref{e301} by
\begin{equation}\label{e310}
\|f\|_{L^\infty(R_{\vec M})}^{(q-p)/q}\left[\mathcal{M}_p^{(n),\otimes}  f(x)\right]^{p/q}\chi_{R_{\vec M}}(x).
\end{equation}
As $t_1,t_2\ge 0$, we can cover $R_{\vec M}$ by at most $2^{2n}$ number of rectangles of the form $R_{j-t_1,k-t_2,\vec M}^{\vec m'}$. Therefore $\|f\|_{L^\infty(R_{\vec M})}^{(q-p)/q}\chi_{R_{\vec M}}$  is bounded by
$$
\sum_{\vec m'}\|f\|_{L^\infty(R_{j-t_1,k-t_2,\vec M}^{\vec m'})}^{(q-p)/q}\chi_{R_{\vec M}}(x)
\lesssim 2^{(M_1+M_2)n(1/p-1/q)}\left[\mathcal{M}_p^{(n),\otimes}  f(x)\right]^{1-p/q}\chi_{R_{\vec M}}(x),
$$
where we use
\eqref{e302} with $q=\infty$.

The remaining cases can be proved similarly.
\end{proof}

\begin{remark}
To control $\big(\frac1{|R_{\vec M}|}\int_{R_{\vec M}}|f(y)|^q\; dy\big)^{1/q}$, a standard dilation argument would yield a growth $2^{(t_1+t_2)n/q}$ as $t_1,t_2\to+\infty$. The better bound in \eqref{e301} for this case is due to   \eqref{e310}, where $\|f\|_{L^\infty(R_{\vec M})}$  behaves better for our purpose, since only for $q=\infty$ we can find a constant $C$, independent of $R'$ and $R$, such that
$$
\Big(\frac{1}{|R|}\int_R|f|^q(y)dy\Big)^{1/q}\le C\Big(\frac{1}{|R'|}\int_{R'}|f|^q(y)dy\Big)^{1/q}
$$
for all $R\subset R'$ and $f$.
\end{remark}

We are now ready to state generalizations of Lemma~\ref{lmmsl1}.
\begin{lemma}\label{lemmassmm}
 Let $1<p<\infty$, $1\leq p_0 <\min \{u,p\}$, and  $M_1,M_2\in\bbn_0$.
 Then we have
 \begin{align}
 \big\Vert \mathcal{SS}^{(u)}_{\vec{M}} F\big\Vert_{L^p(\bbr^{2n})}& \lesssim_{p,p_0 } 2^{(M_1+M_2) n(1/p_0-1/u)} \|F\|_{L^p(\bbr^{2n})}, \label{lmss}\\
  \big\Vert \mathcal{MS}^{(u)}_{\vec{M}} F\big\Vert_{L^p(\bbr^{2n})}& \lesssim_{p,p_0 } 2^{(M_1+M_2) n(1/p_0-1/u)} \|F\|_{L^p(\bbr^{2n})}, \label{lmms}\\
   \big\Vert \mathcal{SM}^{(u)}_{\vec{M}} F\big\Vert_{L^p(\bbr^{2n})}& \lesssim_{p,p_0 } 2^{(M_1+M_2) n(1/p_0-1/u)} \|F\|_{L^p(\bbr^{2n})}, \label{lmsm}\\
    \big\Vert \mathcal{MM}^{(u)}_{\vec{M}} F\big\Vert_{L^p(\bbr^{2n})}& \lesssim_{p,p_0 } 2^{(M_1+M_2) n(1/p_0-1/u)} \|F\|_{L^p(\bbr^{2n})}. \label{lmmm}
 \end{align}

\end{lemma}
The proof of Lemma \ref{lemmassmm} is exactly the same as that of Lemma \ref{lmmsl1} simply by replacing \eqref{maximalbound} with the following inequality:\\
For any $1\le r<\min\{u,p\}$, by Lemma~\ref{t301}, we have
\begin{align}\label{lm23keylemma}
&\Big\|\big(\phi_j\otimes \phi_k\big) \ast F \Big({m_1+y_1 \over 2^{j-M_1}}, {m_2+y_2 \over 2^{k-M_2}}\Big)\Big\|_{L^u(|y|\lesssim 1)}\chi_{R^{\vec m}_{j,k,\vec M}}(x)\nonumber\\
\lesssim&\Big( \frac1{|R^{\vec m}_{j,k,\vec M}|}\int_{c R^{\vec m}_{j,k,\vec M}} \big| (\phi_j\otimes \phi_k) \ast F (z_1, z_2) \big|^u ~dz \Big)^{1/u}\chi_{R^{\vec m}_{j,k,\vec M}}(x)\notag\\
\lesssim& 2^{(M_1+M_2)n(1/p_0-1/u)}\mathcal{M}_{p_0}^{(n),\otimes} \big((\phi_j\otimes \phi_k) \ast F\big)(x)\chi_{R^{\vec m}_{j,k,\vec M}}(x).
\end{align}

\hfill

To deal with the case when $0<r\le 1$ in Theorem \ref{bl}, we follow the strategy introduced in \cite{MPTT2} for which further decomposition would be needed.
For this argument, we first recall the definition of bump functions and then introduce associated maximal inequalities.
 \begin{definition}\label{003}
 Let $I\subseteq \n$ be an arbitrary cube and $N\in\bbn$. A smooth function $\vartheta$ is said to be a bump function adapted to $I$ if there exits a sufficiently large number $M\le N$ such that
 \begin{equation}\label{e010}
 \big|\partial^\alpha\vartheta(x)\big|\lesssim_{\alpha,N}  {1\over{\ell(I)^{|\alpha|}}}{1\over \big(1+|x-x_I|/\ell(I)\big)^N}
 \end{equation}
 for  multi-indices $\alpha\in \mathbb{N}_0^n$ with $|\al|\le M$, where $x_I$ denotes the center of $I$. If $\vartheta_I$ is a bump adapted to $I$, we say that $|I|^{-1/p}\vartheta_I$ is an $L^p$-normalized bump adapted to $I$, for $1\leq p\leq \infty.$
 \end{definition}

 Then we have the following decomposition property, which is an $n$-dimensional analogue of \cite[Lemma 3.9]{MPTT2}.
 \begin{lemma}[\cite{MPTT2}]\label{compact}
  Let $I\subseteq \n$ be a fixed cube. Then every smooth bump function $\vartheta_I$ adapted to $I$ can be decomposed by
  $$\vartheta_I=\sum_{\mu\in \mathbb{N}}2^{-1000n\mu}\vartheta_I^\mu,$$
  where for each $\mu\in \mathbb{N}$, $\vartheta_I^\mu$ is a bump function adapted to $I$ with an additional property that $\supp (\vartheta_I^\mu)\subseteq 2^\mu I$.
  If  $\int_{\n} \vartheta_I \, dx=0$, then correspondingly $\vartheta_I^\mu$ can be chosen with $\int_{\n} \vartheta_I^\mu \,dx=0$.
\end{lemma}
The proof of the lemma is almost the same as that of \cite[Lemma3.9]{MPTT2} and thus it is omitted here.\\

We now apply the above decomposition to Schwartz functions $\varphi$ and $\psi$ satisfying \eqref{suppvarphipsi} as both are adapted to a certain cube with constant side-length, centered at the origin.
Let $Q$ be such a cube to which $\varphi$ and $\psi$ are adapted. Then for each $\mu\in\bbn$, there exist bump functions $\varphi^{\mu}$ and $\psi^{\mu}$ adapted to $Q$ such that
\begin{equation*}
\supp(\varphi^{\mu}), ~\supp(\psi^{\mu}) \subseteq 2^{\mu}Q,
\end{equation*}
\begin{equation*}
\varphi=\sum_{\mu\in\bbn}2^{-1000 n\mu}\varphi^{\mu},\q \text{ and }\q \psi=\sum_{\mu\in\bbn}2^{-1000 n\mu}\psi^{\mu}.
\end{equation*}
Consequently, the dilated functions $\varphi_k$ and $\psi_k$ can be written as
 \begin{eqnarray*}
   \varphi_k=\sum_{\mu \in \mathbb{N}}2^{-1000n\mu} \varphi_k^{\mu} \q \text{ and }\q \psi_k=\sum_{\mu \in \mathbb{N}} 2^{-1000n\mu} \psi_k^{\mu},
 \end{eqnarray*}
 where $\varphi_k^{\mu}:=2^{kn}\varphi^{\mu}(2^k\cdot)$ and $\psi_k^{\mu}:=2^{kn}\psi^{\mu}(2^k\cdot)$ so that $2^{-kn}\varphi_k^{\mu}$ and $2^{-kn}\psi_k^{\mu}$ are adapted to $2^{-k}Q$ and
 \begin{equation}\label{suppphi_k}
 \supp(\varphi_k^{\mu}),\, \supp(\psi_k^{\mu})\subseteq 2^{\mu-k}Q.
 \end{equation}
 Actually, it holds that for any $N\in\bbn$
 $$\big| \partial^{\alpha} \psi_k^{\mu}(x)\big|\lesssim_{\alpha,N}2^{kn}\frac{1}{(2^{-k}\ell(Q))^{|\alpha|}}\frac{1}{\big(1+\frac{|x|}{2^{-k}\ell(Q)}\big)^N}\lesssim 2^{kn}2^{k|\alpha|}\frac{1}{(1+2^k|x|)^N} \q \text{ if }~|\alpha|\le N$$
 and similarly,
 $$\big| \partial^{\alpha} \varphi_k^{\mu}(x)\big|\lesssim 2^{kn}2^{k|\alpha|}\frac{1}{(1+2^k|x|)^N} \q \text{ if }~|\alpha|\le N.$$
 We note that $\psi^{\mu}$ satisfies the mean value zero property
 \begin{equation}\label{meanzeropm}
 \int_{\bbrn}\psi^{\mu}(x)\, dx=0
 \end{equation} as $\psi$ does so.
  Moreover, from the construction of each $\psi^\mu$ (see \cite[Lemma 3.9]{MPTT2} for details), one can see
   \begin{equation}
 \label{lpcd}
 |\psi^\mu(x)|+|\nabla \psi^\mu(x)|\lesssim_M \frac{1}{(1+|x|)^M}, \quad \text{uniformly in } \mu \in \mathbb{N}
 \end{equation}
for any $M>0$. Therefore, as in \eqref{LPtheory} and \eqref{biLPtheory}, the Littlewood-Paley inequalities would be still valid after the above decomposition. That is, for any $1<p<\infty$ and all $f\in L^p(\bbrn)$ and $F\in L^p(\bbr^{2n})$,
 \begin{align}
\big\Vert \big\{ \psi^{\mu_1}_j\ast f\big\}_{j\in\bbz}\big\|_{L^p(\ell^2)}&\lesssim_p \|f\|_{L^p(\bbrn)}\label{lpmu1}\\
   \big\Vert \big\{ (\psi^{\mu_1}_j\otimes\psi^{\mu_2}_k)\ast F\big\}_{j,k\in\bbz}\big\|_{L^p(\ell^2)}&\lesssim_p \|F\|_{L^p(\bbr^{2n})}  \label{lpmu2}
 \end{align}
 uniformly in $\mu_1,\mu_2\in\bbn$.

 \hfill

Now we define hybrid operators related to $\varphi_k^{\mu}$ and $\psi_k^{\mu}$:\\
For $F\in\mathscr{S}(\bbrn\times\bbrn)$ and any $\mmu:=(\mu_1,\mu_2)\in \bbn^2$,
\begin{align*}
  \mathcal{SS}^{(u),\mmu}_{\vec{M}}(F)(x)&:=  \Big( \sum_{\substack{j,k\in \z \\ m_1,m_2\in \z^n}} \Big\|(\psi_j^{\mu_1}\otimes \psi_k^{\mu_2}) \ast F \Big({m_1+y_1 \over 2^{j-M_1}}, {m_2+y_2 \over 2^{k-M_2}}\Big)\Big\|_{L^u(|y|\lesssim 1)}^2 \chi_{R^{\vec m}_{j,k,\vec M}}(x)\Big)^{1/ 2} , \\
    \mathcal{SM}^{(u),\mmu}_{\vec{M}}(F)(x)&:=  \Big( \sum_{\substack{j\in \z \\ m_1 \in \z^n}}  \sup_{\substack{k \in \z\\ m_2 \in \z^n}} \Big\|\big(\psi_j^{\mu_1}\otimes \varphi_k^{\mu_2}\big) \ast F \Big({m_1+y_1 \over 2^{j-M_1}}, {m_2+y_2 \over 2^{k-M_2}}\Big)\Big\|_{L^u(|y|\lesssim 1)}^2 \chi_{R^{\vec m}_{j,k,\vec M}}(x) \Big)^{1/ 2},
   \\
   \mathcal{MS}^{(u),\mmu}_{\vec{M}}(F)(x)&:=    \sup_{\substack{j \in \z\\ m_1 \in \z^n}}\Big( \sum_{\substack{k\in \z \\ m_2 \in \z^n}}  \Big\|\big(\varphi_j^{\mu_1}\otimes \psi_k^{\mu_2}\big) \ast F \Big({m_1+y_1 \over 2^{j-M_1}}, {m_2+y_2 \over 2^{k-M_2}}\Big)\Big\|_{L^u(|y|\lesssim 1)}^2 \chi_{R^{\vec m}_{j,k,\vec M}}(x) \Big)^{1/ 2},
  \\
    \mathcal{MM}^{(u),\mmu}_{\vec{M}}(F)(x)&:=  \sup_{j,k\in \z} \sup_{ m_1,m_2 \in \n} \Big\|\big(\varphi_j^{\mu_1}\otimes \varphi_k^{\mu_2}\big) \ast F \Big({m_1+y_1 \over 2^{j-M_1}}, {m_2+y_2 \over 2^{k-M_2}}\Big)\Big\|_{L^u(|y|\lesssim 1)} \chi_{R^{\vec m}_{j,k,\vec M}}(x).
\end{align*}

Then as a counterpart of Lemma \ref{lemmassmm}, we have the following lemma.
\begin{lemma} \label{lmnoncompactbip}
Let  $1<p<\infty$ , $1\leq  p_0 <\min(u,p)$, and $M_1,M_2\in \bbn_0$.
Then there exists a large number $M$ such that for $F\in \mathscr{S}(\mathbb{R}^{2n})$ we have
 \begin{align*}
 \big\Vert    \mathcal{SS}^{(u),\mmu}_{\vec{M}}(F)   \big\Vert_{L^p(\bbr^{2n})}&\lesssim 2^{(\mu_1+\mu_2)M}\,  2^{(M_1+M_2)n(1/p_0-1/u)}  \|F\|_{L^p(\mathbb R^{2n})},\\
   \big\Vert    \mathcal{SM}^{(u),\mmu}_{\vec{M}}(F)   \big\Vert_{L^p(\bbr^{2n})}&\lesssim 2^{(\mu_1+\mu_2)M}\,  2^{(M_1+M_2)n(1/p_0-1/u)}  \|F\|_{L^p(\mathbb R^{2n})},\\
  \big\Vert    \mathcal{MS}^{(u),\mmu}_{\vec{M}}(F)   \big\Vert_{L^p(\bbr^{2n})}&\lesssim 2^{(\mu_1+\mu_2)M}\,  2^{(M_1+M_2)n(1/p_0-1/u)}  \|F\|_{L^p(\mathbb R^{2n})},\\
    \big\Vert    \mathcal{MM}^{(u),\mmu}_{\vec{M}}(F)   \big\Vert_{L^p(\bbr^{2n})}&\lesssim 2^{(\mu_1+\mu_2)M}\,  2^{(M_1+M_2)n(1/p_0-1/u)}  \|F\|_{L^p(\mathbb R^{2n})}.
 \end{align*}
 \end{lemma}
In order to prove Lemma \ref{lmnoncompactbip}, similar to \eqref{lm23keylemma}, we shall estimate the term
$$\Big\|\big(\phi_j^{\mu}\otimes \phi_k^{\mu}\big) \ast F \Big({m_1+y_1 \over 2^{j-M_1}}, {m_2+y_2 \over 2^{k-M_2}}\Big)\Big\|_{L^u(|y|\lesssim 1)}\chi_{R^{\vec m}_{j,k,\vec M}}(x),$$
but the arguments used to prove \eqref{lm23keylemma} cannot be directly applied as the Fourier transform of $\phi_k^{\mu}$ should not have compact support due to \eqref{suppphi_k}. Therefore, we perform further decomposition of $\phi_k^{\mu}$.

 Now  we pick a family of Schwartz functions $\{\Gamma^{\mu}_{k,l}\}_{l\in\N_0}$ such that
$\sum_{l\in \N_0} \widehat{\Gamma^{\mu}_{k,l}}=1$ with the compact Fourier support conditions
 $$  \supp (\widehat{\Gamma^{\mu}_{k,0}})\subseteq \{\xi\in\bbrn: |\xi|\lesssim 2^{-\mu+k}\}, $$
and
$$\supp (\widehat{\Gamma^{\mu}_{k,l}})\subseteq \{\xi\in\bbrn: |\xi|\approx 2^{-\mu+k+l}\} \q \text{ for }~l\in\bbn.$$
Then we write $ \widehat {\Phi_{k,{l}}^{\mu}:}=\widehat {\Gamma^{\mu}_{k,{l}}}\, \widehat {\varphi^{\mu}_k}$ and $ \widehat {\Psi_{k,{l}}^{\mu}:}=\widehat {\Gamma^{\mu}_{k,{l}}}\, \widehat {\psi^{\mu}_k}$  so that
 \begin{equation}
     \varphi^{\mu}_k= \sum_{{l}\in \N_0} \Phi_{k,{l}}^{\mu}\q \text{ and }\q  \psi^{\mu}_k= \sum_{{l}\in \N_0} \Psi_{k,{l}}^{\mu} .   \label{e0111}
 \end{equation}

Intuitively, as $\varphi_0^\mu$ is adapted to $Q=[-1,1]^{2n}$, $\wh {\varphi_0^\mu}$ is adapted to $Q$ as well. By the rapid decay of $\wh{\varphi_0^\mu}$,  we see that $2^{lM}\wh{\Phi_0^\mu}$ and $2^{lM}\Phi_0^\mu$ are adapted to $Q$ too.
To capture the decay in $l$ as we observed, we have the following estimates.

\begin{lemma}\label{005d}
Let $N$ be a sufficiently large integer greater than $n$ and choose two integers $0<M\le N$ and  $0<L<N-n$.
Then we have
\begin{equation}\label{lm26fstd}
\big\Vert \wh{\Phi_{k,l}^{\mu}}\big\Vert_{L^{\infty}(\bbrn)},\, \big\Vert \wh{\Psi_{k,l}^{\mu}}\big\Vert_{L^{\infty}(\bbrn)}\, \lesssim_{M,N} 2^{(\mu-l)M} 
\end{equation}
and
 \begin{equation}\label{partialphiklmud}
 \big| \partial^{\alpha}\Phi_{k,l}^{\mu}(x)\big|, \, \big| \partial^{\alpha}\Psi_{k,l}^{\mu}(x)\big|\, \lesssim_{M,L,N,\alpha} 2^{kn}2^{(\mu-l) (M-n)}
 2^{|\alpha|(k-\mu+l)}\frac{1}{(1+2^{k-\mu}|x|)^L}
 \end{equation}
 for any multi-indices $\alpha$.
\end{lemma}
 \begin{proof}
 Without loss of generality, we may assume $k=0$, as the other cases when $k\in\bbz$ simply follow from the dilation
  $$\Phi_{k,l}^{\mu}(x)=2^{kn}\Phi_{0,l}^{\mu}(2^kx)\q \text{ and }\q\Psi_{k,l}^{\mu}(x)=2^{kn}\Psi_{0,l}^{\mu}(2^kx).$$
As before, we use $\phi$ to indicate $\varphi$ or $\psi$.
Since $\phi^\mu$ are smooth bumps adapted to $[-2,2]^n$ (uniformly in $\mu$) by Lemma~\ref{compact}, the inequality
 \begin{equation*}
 \big| \partial^{\alpha}\phi^{\mu}(x)\big|\lesssim_{\alpha,N} (1+|x|)^{-N}
 \end{equation*}
holds for all $N>0$ and all multi-indices $\alpha$ with $|\alpha|\le M$. We set $N>n$ to be a sufficiently large number.
Then for any multi-indices $\beta$ with $|\beta|<N-n$, using the fact that $\spt \phi^{\mu}\subset [-2^{\mu+1}, 2^{\mu+1}]$, we have
 \begin{eqnarray}\label{delbetawhd}
 \big| \partial^{\beta}\wh{\phi^{\mu}}(\xi)\big|&\approx& \Big| \int_{\bbrn} x^{\beta}\phi^{\mu}(x) e^{2\pi i\langle x,\xi\rangle}~dx\Big|  \nonumber\\
 & \lesssim_{\beta,N,M} &  (1+|\xi|)^{-M} ~\Big| \int_{\bbrn} \frac{|x|^{|\beta|}}{(1+|x|)^N} ~dx\Big|\lesssim (1+|\xi|)^{-M}  \nonumber,
 \end{eqnarray}
 which follows from integration by parts since $N-|\beta|>n$. This implies that
 \begin{align}\label{parbetaxi}
 \big| \partial^{\beta}\big(\wh{\phi^{\mu}}\wh{\Gamma_{0,l}^{\mu}}\big)(\xi)\big|& \lesssim \sum_{\beta'\le \beta}\big|\partial^{\beta'}\wh{\phi^{\mu}}(\xi) \big| \, \big|\partial^{\beta-\beta'}\wh{\Gamma_{0,l}^{\mu}}(\xi) \big|\nonumber\\
 &\lesssim 2^{\mu |\beta|}(1+|\xi|)^{-M}\begin{cases}
 \chi_{|\xi|\approx 2^{-\mu+l}} & l\in\bbn\\
 \chi_{|\xi|\lesssim 2^{-\mu}} & l=0
 \end{cases}\notag\\
  &\lesssim 2^{\mu M}2^{\mu |\beta|}(1+2^{\mu}|\xi|)^{-M}\begin{cases}
 \chi_{|\xi|\approx 2^{-\mu+l}} & l\in\bbn\\
 \chi_{|\xi|\lesssim 2^{-\mu}} & l=0
 \end{cases}
 \end{align}

 Finally, we deduce \eqref{lm26fstd} from \eqref{parbetaxi}.

 Moreover, when $l\ge 1$, for any $0<L<N-n$ and $\mu\ge 0$,  we have
  \begin{align*}
  \big| \partial^{\alpha} (\phi^{\mu}{\ast} \Gamma^\mu_{0,l})(x)\big| &\lesssim \Big| \int_{|\xi| \approx 2^{-\mu+l}}  \xi^\alpha  \wh{\phi^{\mu}}(\xi) \wh{\Gamma^\mu_{0,l}}(\xi ) e^{2\pi i\langle x,\xi\rangle}\;dx\Big|\\
 &\lesssim |x|^{-L}\int_{|\xi|\approx 2^{-\mu+l}}\sum_{|\beta|+|\beta'|= L} \big|\partial^{\beta}\xi^\alpha \big| \big| \partial^{\beta'}_{\xi} \big(\wh{\phi^{\mu}} \wh{\Gamma^\mu_{0,l}}\big)(\xi ) \big|\; d\xi\\
 &\lesssim(2^{-\mu}|x|)^{-L} 2^{\mu M}2^{-\mu|\alpha|}
\int_{|\xi|\approx2^{-\mu+l}} \sum_{\beta\le\alpha}  (2^\mu|\xi|)^{|\alpha-\beta|}(1+2^\mu|\xi|)^{-M}\;d\xi\\
&\approx
(2^{-\mu}|x|)^{-L} 2^{\mu M}2^{-\mu|\alpha|} 2^{-(\mu-l)n}
\int_{|\eta|\approx1} \sum_{\beta\le\alpha}  (2^l|\eta|)^{|\alpha-\beta|}(1+2^l|\eta|)^{-M}\;d\eta
\\
 &\lesssim (2^{-\mu}|x|)^{-L} 2^{(\mu-l)(M-n)}2^{(-\mu+l)|\alpha|}
  \end{align*}
 where the third inequality follows from \eqref{parbetaxi}.
 Similarly, one can verify that
 $$
   \big| \partial^{\alpha} (\phi^{\mu}{\ast} \Gamma^\mu_{0,l})(x)\big|\lesssim 2^{(\mu-l)(M-n)}2^{(-\mu+l)|\alpha|}.
 $$
 The above two estimates yield \eqref{partialphiklmud} for $k=0$ and $l\ge 1$.
 One can easily modify this argument to verify the case when $k=0$ and $l=0$.
In summary, we obtain the inequality \eqref{partialphiklmud} when $k=0$.
 \end{proof}

According to Lemma \ref{005d},
$$\big| \Psi_{\mu,l}^{\mu}(x)\big|+\big| \nabla\Psi_{\mu,l}^{\mu}(x)\big|\lesssim 2^{\mu (M+1)}2^{-l(M-n-1)}\frac{1}{(1+|x|)^L}.$$
Moreover, $\Psi_{\mu,l}^{\mu}=\Gamma_{\mu,l}^{\mu}\ast \psi_{\mu}^{\mu}$ satisfies the mean value zero property, which is inherited from \eqref{meanzeropm}.
Therefore,  it follows from the Littlewood-Paley theroy \eqref{LPtheory} that if $M>0$ is sufficiently large, then
\begin{equation}\label{molpthoery1}
\big\Vert \big\{ (\Psi_{k,l}^{\mu}\ast f\big\}_{k\in\bbz}\big\Vert_{L^p(\ell^2)}\lesssim 2^{\mu M}2^{-l(M-n-2)}\Vert f\Vert_{L^p(\bbr^{n})}
\end{equation} holds for $f\in L^p(\bbrn)$,
and the bi-paramter version also holds due to \eqref{biLPtheory} so that
\begin{equation}\label{molpthoery2}
\big\Vert \big\{ (\Psi_{j,l_1}^{\mu_1}\otimes\Psi_{k,l_2}^{\mu_2})\ast F\big\}_{j,k\in\bbz}\big\Vert_{L^p(\ell^2)}\lesssim 2^{(\mu_1+\mu_2)M}2^{-(l_1+l_2)(M-n{-2})}\Vert F\Vert_{L^p(\bbr^{2n})}
\end{equation} holds
 for $F\in L^p(\bbr^{2n})$ as $\Psi_{k,l}^{\mu}(x)=2^{(k-\mu)n}\Psi_{\mu,l}^{\mu}(2^{k-\mu}x)$ for all $k\in\bbz$.\\

 Now let us prove lemma \ref{lmnoncompactbip}
\begin{proof}[Proof of lemma~\ref{lmnoncompactbip}]
We first consider the operator $\mathcal{SS}_{\vec M}^{(u),\mmu}$.
Using the decomposition \eqref{e0111}, we write
\begin{align*}
& \Big\|(\psi_j^{\mu_1}\otimes \psi_k^{\mu_2}) \ast F \Big({m_1+y_1 \over 2^{j-M_1}}, {m_2+y_2 \over 2^{k-M_2}}\Big)\Big\|_{L^u(|y|\lesssim 1)} \chi_{R^{\vec m}_{j,k,\vec M}}(x)\\
\lesssim&  \sum_{l_1,l_2\in\bbn_0}\Big\|(\Psi_{j,l_1}^{\mu_1}\otimes \Psi_{k,l_2}^{\mu_2}) \ast F \Big({m_1+y_1 \over 2^{j-M_1}}, {m_2+y_2 \over 2^{k-M_2}}\Big)\Big\|_{L^u(|y|\lesssim 1)} \chi_{R^{\vec m}_{j,k,\vec M}}(x)\\
\lesssim&\sum_{l_1,l_2\in\bbn_0}2^{n(M_1+M_2)(1/p_0-1/u)}       \big( \max\{2^{-n(\mu_1-l_1)},1\}\big)^{1/p_0-1/u}\\
&\q\times \big(\max\{2^{-n(\mu_2-l_2)},1\}\big)^{1/p_0-1/u} \mathcal{M}_{p_0}^{(n),\otimes}\big((\Psi_{j,l_1}^{\mu_1}\otimes \Psi_{k,l_2}^{\mu_2})\ast F \big)(x)\chi_{R^{\vec m}_{j,k,\vec M}}(x)
\end{align*}
where Lemma~\ref{t301} is applied in last step. Therefore, we have proved that
\begin{align*}
 &\big\Vert    \mathcal{SS}^{(u),\mmu}_{\vec{M}}(F)   \big\Vert_{L^p(\bbr^{2n})}\\
 &\lesssim 2^{n(M_1+M_2)(1/p_0-1/u)}  \sum_{l_1,l_2\in\bbn_0} \big( \max\{2^{-n(\mu_1-l_1)},1\}\big)^{1/p_0-1/u} \big(\max\{2^{-n(\mu_2-l_2)},1\}\big)^{1/p_0-1/u}\\
 &\qq\qq\qq\times   \Big\Vert \Big( \sum_{j,k\in\bbz}  \Big( \mathcal{M}_{p_0}^{(n),\otimes}\big((\Psi_{j,l_1}^{\mu_1}\otimes \Psi_{k,l_2}^{\mu_2})\ast F \big)   \Big)^2     \sum_{m_1,m_2\in\bbzn} \chi_{R^{\vec m}_{j,k,\vec M}}
   \Big)^{1/2}   \Big\Vert_{L^p(\bbr^{2n})}\\
    &\lesssim 2^{n(M_1+M_2)(1/p_0-1/u)}  \sum_{l_1,l_2\in\bbn_0} \big( \max\{2^{-n(\mu_1-l_1)},1\}\big)^{1/p_0-1/u}\\
    &\qq \qq\times \big(\max\{2^{-n(\mu_2-l_2)},1\}\big)^{1/p_0-1/u} \big\Vert \big\{      \mathcal{M}_{p_0}^{(n),\otimes}\big((\Psi_{j,l_1}^{\mu_1}\otimes \Psi_{k,l_2}^{\mu_2})\ast F \big)        \big\}_{j,k\in\bbz}\big\Vert_{L^p(\ell^2)}\\
    &\lesssim 2^{n(M_1+M_2)(1/p_0-1/u)}  \sum_{l_1,l_2\in\bbn_0} \big( \max\{2^{-n(\mu_1-l_1)},1\}\big)^{1/p_0-1/u}\\
    &\qq\qq \times \big(\max\{2^{-n(\mu_2-l_2)},1\}\big)^{1/p_0-1/u}\big\Vert \big\{    (\Psi_{j,l_1}^{\mu_1}\otimes \Psi_{k,l_2}^{\mu_2})\ast F         \big\}_{j,k\in\bbz}\big\Vert_{L^p(\ell^2)}\\
    &\lesssim 2^{n(M_1+M_2)(1/p_0-1/u)} \sum_{l_1,l_2\in\bbn_0}  \big( \max\{2^{-n(\mu_1-l_1)},1\}\big)^{1/p_0-1/u}\\
    &\qq\qq\times \big(\max\{2^{-n(\mu_2-l_2)},1\}\big)^{1/p_0-1/u}2^{(\mu_1+\mu_2)M}2^{-(l_1+l_2)(M-n-2)}\Vert F\Vert_{L^p(\bbr^{2n})}\\
    &= 2^{n(M_1+M_2)(1/p_0-1/u)}2^{(\mu_1+\mu_2)M}\Vert F\Vert_{L^p(\bbr^{2n})}
\end{align*}
for sufficiently large $M>0$, as desired, where the last inequality follows form \eqref{molpthoery2}.\\

Similarly, using Lemma~\ref{t301} and \eqref{e0111}, we have
\begin{align*}
& \Big\|(\psi_j^{\mu_1}\otimes \vp_k^{\mu_2}) \ast F \Big({m_1+y_1 \over 2^{j-M_1}}, {m_2+y_2 \over 2^{k-M_2}}\Big)\Big\|_{L^u(|y|\lesssim 1)} \chi_{R^{\vec m}_{j,k,\vec M}}(x)\\
\lesssim & \sum_{l_1,l_2\in\bbn_0}\Big\|(\Psi_{j,l_1}^{\mu_1}\otimes \Phi_{k,l_2}^{\mu_2}) \ast F \Big({m_1+y_1 \over 2^{j-M_1}}, {m_2+y_2 \over 2^{k-M_2}}\Big)\Big\|_{L^u(|y|\lesssim 1)} \chi_{R^{\vec m}_{j,k,\vec M}}(x)\\
\lesssim &2^{n(M_1+M_2)(1/p_0-1/u)}\sum_{l_1,l_2\in\bbn_0}\big( \max\{2^{-n(\mu_1-l_1)},1\}\big)^{1/p_0-1/u}\\
&\q\times  \big(\max\{2^{-n(\mu_2-l_2)},1\}\big)^{1/p_0-1/u}\mathcal{M}_{p_0}^{(n),\otimes}\big((\Psi_{j,l_1}^{\mu_1}\otimes \Phi_{k,l_2}^{\mu_2})\ast F \big)(x)\chi_{R_{j,k,\vec{M}}^{\vec m}}(x)
\end{align*}
and thus
\begin{align*}
 &\big\Vert    \mathcal{SM}^{(u),\mmu}_{\vec{M}}(F)   \big\Vert_{L^p(\bbr^{2n})}\lesssim 2^{n(M_1+M_2)(1/p_0-1/u)}  \sum_{l_1,l_2\in\bbn_0} \big( \max\{2^{-n(\mu_1-l_1)},1\}\big)^{1/p_0-1/u}\\
 &\q\times  \big(\max\{2^{-n(\mu_2-l_2)},1\}\big)^{1/p_0-1/u}  \Big\Vert \Big( \sum_{j\in\bbz}  \Big( \sup_{k\in\bbz}\mathcal{M}_{p_0}^{(n),\otimes}\big((\Psi_{j,l_1}^{\mu_1}\otimes \Phi_{k,l_2}^{\mu_2})\ast F \big)   \Big)^2       \Big)^{1/2}   \Big\Vert_{L^p(\bbr^{2n})}.
   \end{align*}
Now we see that
\begin{align}\label{phipsiastf}
 &\Big\Vert \Big( \sum_{j\in\bbz}  \Big( \sup_{k\in\bbz}\mathcal{M}_{p_0}^{(n),\otimes}\big((\Psi_{j,l_1}^{\mu_1}\otimes \Phi_{k,l_2}^{\mu_2})\ast F \big)   \Big)^2       \Big)^{1/2}   \Big\Vert_{L^p(\bbr^{2n})}\nonumber\\
 &\le   \Big\Vert \Big( \sum_{j\in\bbz}  \Big(\mathcal{M}_{p_0}^{(n),\otimes}\big(\sup_{k\in\bbz}|(\Psi_{j,l_1}^{\mu_1}\otimes \Phi_{k,l_2}^{\mu_2})\ast F| \big)\Big)^2       \Big)^{1/2}   \Big\Vert_{L^p(\bbr^{2n})}\nonumber\\
&\lesssim  \Big\Vert \Big( \sum_{j\in\bbz}  \Big(\sup_{k\in\bbz}|(\Psi_{j,l_1}^{\mu_1}\otimes \Phi_{k,l_2}^{\mu_2})\ast F| \Big)^2       \Big)^{1/2}   \Big\Vert_{L^p(\bbr^{2n})}
\end{align}
and in view of \eqref{partialphiklmud},
\begin{align*}
&\sup_{k\in\bbz}|(\Phi_{j,l_1}^{\mu_1}\otimes \Psi_{k,l_2}^{\mu_2})\ast F(x_1,x_2)|\lesssim 2^{(\mu_2-l_2)M}\mathcal{M}^{(n)}\big( \Psi_{j,l_1}^{\mu_1}\ast_1 F(x_1,\cdot)\big)(x_2)
\end{align*}
for sufficiently large $M>0$, where $\ast_1$ means the convolution in the first variable $x_1$ so that
$$\Psi_{j,l_1}^{\mu_1}\ast_1 F(x_1,x_2)=\int_{\bbrn}\Psi_{j,l_1}^{\mu_1}(x_1-y_1)F(y_1,x_2)~dy_1.$$
Now, applying \eqref{multimaximalin} and \eqref{molpthoery1}, the right-hand side of \eqref{phipsiastf} is dominated by a constant times
\begin{align*}
&2^{(\mu_2-l_2)M}\bigg\Vert\Big\Vert \Big( \sum_{j\in\bbz}  \Big|\Psi_{j,l_1}^{\mu_1}\ast_1 \big(F(\cdot,x_2)\big)(x_1) \Big|^2       \Big)^{1/2}   \Big\Vert_{L^p(x_2)}\bigg\Vert_{L^p(x_1)}\\
\lesssim& 2^{(\mu_1+\mu_2)M}2^{-(l_1+l_2)(M-n-2)}\Vert F\Vert_{L^p(\bbr^{2n})},
\end{align*}
which finally proves
\begin{align*}
&\big\Vert    \mathcal{SM}^{(u),\mmu}_{\vec{M}}(F)   \big\Vert_{L^p(\bbr^{2n})}\\
&\lesssim 2^{n(M_1+M_2)(1/p_0-1/u)} \sum_{l_1,l_2\in\bbn_0}  \big( \max\{2^{-n(\mu_1-l_1)},1\}\big)^{1/p_0-1/u}\\
    &\qq\times \big(\max\{2^{-n(\mu_2-l_2)},1\}\big)^{1/p_0-1/u}2^{(\mu_1+\mu_2)M}2^{-(l_1+l_2)(M-n-2)}\Vert F\Vert_{L^p(\bbr^{2n})}\\
    &= 2^{n(M_1+M_2)(1/p_0-1/u)}2^{(\mu_1+\mu_2)M}\Vert F\Vert_{L^p(\bbr^{2n})}.
\end{align*}

\hfill

As $    \mathcal{MS}^{(u),\mmu}_{\vec{M}}(F)   $ is controlled by the corresponding $    \mathcal{SM}^{(u),\mmu}_{\vec{M}}(F)   $ pointwisely, we have
\begin{equation*}
\big\Vert    \mathcal{MS}^{(u),\mmu}_{\vec{M}}(F)   \big\Vert_{L^p(\bbr^{2n})}\lesssim 2^{n(M_1+M_2)(1/p_0-1/u)}2^{(\mu_1+\mu_2)M}\Vert F\Vert_{L^p(\bbr^{2n})}.
\end{equation*}

\hfill

For the last one, we apply similar arguments to prove
\begin{align*}
&\big\Vert    \mathcal{MM}^{(u),\mmu}_{\vec{M}}(F)   \big\Vert_{L^p(\bbr^{2n})}\lesssim 2^{n(M_1+M_2)(1/p_0-1/u)}  \sum_{l_1,l_2\in\bbn_0} \big( \max\{2^{-n(\mu_1-l_1)},1\}\big)^{1/p_0-1/u}\\
 &\q\times  \big(\max\{2^{-n(\mu_2-l_2)},1\}\big)^{1/p_0-1/u}  \Big\Vert  \sup_{j,k\in\bbz}\mathcal{M}_r^{(n),\otimes}\big((\Phi_{j,l_1}^{\mu_1}\otimes \Phi_{k,l_2}^{\mu_2})\ast F \big)     \Big\Vert_{L^p(\bbr^{2n})}.
\end{align*}
Simply using \eqref{simplemaximalin} and \eqref{partialphiklmud}, we can bound the $L^p$ norm in the last expression by
$$2^{(\mu_1+\mu_2)M}2^{-(l_1+l_2)(M-n-2)}\Vert F\Vert_{L^p(\bbr^{2n})},$$
which completes the proof of the estimate for $ \mathcal{MM}^{(u),\mmu}_{\vec{M}}$.
\end{proof}

\hfill

 \section{The main reduction}\label{reductionsection}
 In this section, we will proceed the main reduction of the operator \eqref{op}, in  preparation for proving Theorem \ref{bl}.
We will follow the setups in  \cite{MPTT1, MPTT2}.

Let us take a function $\psi \in \sch$ satisfying $\spt \wh{\psi} \subseteq  \{\xi\in\n:\ 2^{-1}\leq \left|\xi \right| \leq 2\}$, and
  $$\sum_j \widehat {\psi_j}({\xi})=1, \ \ \  \xi \neq 0, $$
  where $\psi_j:=2^{jn} \psi(2^j \cdot)$.
Then $T_m(f,g)(x)$ can be decomposed as
\begin{equation}
  \sum_{j_1,k_1\in \z}\sum_{j_2,k_2\in \z}\int_{\bn\times \bn} m(\xi,\eta)  \wh{\psi_{j_1}}({\xi_1 }) \wh{\psi_{j_2}}({\xi_2}) \wh{\psi_{k_1}}({\eta_1}) \wh{\psi_{k_2}}({\eta_2 }) \wh{f}(\xi) \wh{g}(\eta) e^{2\pi i \langle x,\xi+\eta\rangle}\; d\xi d\eta.  \label{Tdecompose}
\end{equation}
We express
\begin{align}
 &  \sum_{j_1,k_1 \in \z} \wh{\psi_{j_1}}({\xi_1})\wh{\psi_{k_1}}({\eta_1 }) \nonumber\\
 &= \sum_{j_1\ll k_1} \wh{\psi_{j_1}}({\xi_1 })\wh{\psi_{k_1}}({\eta_1})+\sum_{j_1\approx k_1} \wh{\psi_{j_1}}({\xi_1})\wh{\psi_{k_1}}({\eta_1})+\sum_{j_1\gg k_1}  \wh{\psi_{j_1}}({\xi_1 })\wh{\psi_{k_1}}({\eta_1})\nonumber \\
 &=
 \sum_{ j \in \z} \wh{\varphi_{j}}({\xi_1 })\wh{\psi_{j}}({\eta_1 })+\sum_{j\in \z} \wh{\psi_j}({\xi_1 })\wh{\psi_j}({\eta_1})+\sum_{j\in \z} \wh{\psi_j}({\xi_1 })\wh{\varphi_j}({\eta_1 })
 \label{phipsi0}
\end{align}
where  $j\ll k$ means $j-k<-100$ and $\wh{\varphi_{j}} :=\sum_{k\ll j}\wh{\psi_k} $ so that $\spt \wh{\varphi_j}\subset \{\xi\in \bbrn: |\xi|\le 2^{-100+j}\}$. Clearly, $\psi$ and $\varphi$ are $\Psi$-type and $\Phi$-type functions we discussed in the previous section, respectively.
A similar argument can be made for  other sums to write
\begin{equation*}
 \sum_{j_2,k_2 \in \z} \wh{\psi_{j_2}}({\xi_2})\wh{\psi_{k_2}}({\eta_2 }) = \sum_{ k \in \z} \wh{\varphi_{k}}({\xi_2 })\wh{\psi_{k}}({\eta_2 })+\sum_{k\in \z} \wh{\psi_k}({\xi_2 })\wh{\psi_k}({\eta_2})+\sum_{k\in \z} \wh{\psi_k}({\xi_2 })\wh{\varphi_k}({\eta_2 })
\end{equation*}
 and thus
\begin{align*}
  & \sum_{j_1,k_1 \in \z} \sum_{j_2,k_2 \in \z} \wh{\psi_{j_1}}({\xi_1 })\wh{\psi_{k_1}}({\eta_1})\wh{\psi_{j_2}}({\xi_2 })\wh{\psi_{k_2}}({\eta_1}) \\
  &=\Big( \sum_{ j \in \z} \wh{\varphi_{j}}({\xi_1 })\wh{\psi_{j}}({\eta_1 })+\sum_{j\in \z} \wh{\psi_j}({\xi_1 })\wh{\psi_j}({\eta_1})+\sum_{j\in \z} \wh{\psi_j}({\xi_1 })\wh{\varphi_j}({\eta_1 })\Big) \\
  & \q\q \times \Big( \sum_{ k \in \z} \wh{\varphi_{k}}({\xi_2 })\wh{\psi_{k}}({\eta_2 })+\sum_{k\in \z} \wh{\psi_k}({\xi_2 })\wh{\psi_k}({\eta_2})+\sum_{k\in \z} \wh{\psi_k}({\xi_2 })\wh{\varphi_k}({\eta_2 })\Big) \\
 &=:\big(A_{12}+A_{22}+A_{21}\big)(\xi_1,\eta_1)\times  \big( B_{12}+B_{22}+B_{21}\big)(\xi_2,\eta_2).
\end{align*}
Based on the decomposition, $T_m(f,g)$ can be separated into $9$ pieces, and by symmetry it suffices to consider just the following four terms;
\begin{equation}\label{reduced432}
\begin{aligned}
  A_{12}B_{12}&=&\sum_{j,k\in\z}\widehat{\varphi_j}(\xi_1) \widehat{\varphi_k}(\xi_2) \widehat{\psi_j}(\eta_1) \widehat{\psi_k}(\eta_2), \\
   A_{12}B_{22}&=&\sum_{j,k\in\z}\widehat{\varphi_j}(\xi_1) \widehat{\psi_k}(\xi_2) \widehat{\psi_j}(\eta_1) \widehat{\psi_k}(\eta_2), \\
   A_{12}B_{21}&=&\sum_{j,k\in\z}\widehat{\varphi_j}(\xi_1) \widehat{\psi_k}(\xi_2) \widehat{\psi_j}(\eta_1) \widehat{\varphi_k}(\eta_2),  \\
   A_{22}B_{22}&=& \sum_{j,k\in\z}\widehat{\psi_j}(\xi_1) \widehat{\psi_k}(\xi_2) \widehat{\psi_j}(\eta_1) \widehat{\psi_k}(\eta_2).
   \end{aligned}
\end{equation}
Next we insert some auxiliary $\Psi$-type and $\Phi$-type functions (without changing the decomposition) to share the support information of the decomposition with $m$ and $\widehat{f}(\xi)\wh{g}(\eta) $.
Actually by choosing a function $\widehat{\wt{\Theta}}(\xi_1,\eta_1)$, defined on $\bbrn\times\bbrn$, supported in $\{(\xi_1,\eta_1)\in\mathbb R^{2n}: {100}^{-1}\le |(\xi_1,\eta_1)|\le 100\}$ such that  $\widehat{\wt{\Theta}}(\xi_1,\eta_1)=1$ in the support of $\widehat\vp\otimes \widehat\psi$ and letting $\wt{\Theta_j}(x_1,y_1):=2^{2jn}\wt{\Theta}(2^jx_1,2^jy_1)$, we obtain
\begin{equation} \label{phipsi}
  \sum_{ j \in \z} \wh{\varphi_j}({\xi_1})\wh{\psi_j}({\eta_1 })=\sum_{ j \in \z} \wh{\varphi_j}({\xi_1 }) 
  \wh{\psi_j}({\eta_1})\wh{\wt{\psi_j}}(\xi_1+\eta_1) \wh{\wt {\Theta_j}}(\xi_1,\eta_1),
\end{equation}
where  $\wt{\psi_j}$ is  completed in a way that $ \widehat{{\wt{\psi}}}({\xi_1})=1$ for  $2^{-2}\le |\xi_1|\le 2^2$, and also supported in $\{\xi_1\in \n:\ 2^{-3}\le |\xi|\le 2^3\}$, and $\wt{\psi_j}:=2^{jn}\wt{\psi}(2^j\cdot)$.
Note that $\wt{\psi}$ is of $\Psi$-type as well.
One will see that the ``type" of these functions matters in our proof, thus without lose of generality, we will not distinguish between 
 $\psi$ and $\wt{\psi}$
 in notations, and simplify \eqref{phipsi} as
\begin{equation}
\label{phipsi2}
  \sum_{ j \in \z} \wh{\varphi_j}({\xi_1})\wh{\psi_j}({\eta_1 })=\sum_{ j \in \z} \wh{\varphi_j}({\xi_1 })
  \wh{\psi_j}({\eta_1})\wh{{\psi_j}}(\xi_1+\eta_1) \wh{\wt {\Theta_j}}(\xi_1,\eta_1).
\end{equation}
Similarly, we can write
\begin{align*}
  \sum_{j \in \z} \wh{\psi_j}({\xi_1})\wh{\psi_j}({\eta_1})
  &=\sum_{ j \in \z} \widehat {\psi_{j}}({\xi_1 })   \widehat {\psi_{j}} ({\eta_1 })  \widehat {{{\varphi}}_{j}}({\xi_1+\eta_1})  \widehat{ \wt{\Theta_{j}}}({\xi_1,\eta_1}),\\
  \sum_{j \in \z} \wh{\psi_j}({\xi_1 })\wh{\varphi_j}({\eta_1 })
  &=\sum_{ j \in \z}  \widehat {\psi_{j}} ({\xi_1 })  \widehat {\varphi_{j}}({\eta_1 }) \widehat{{\psi}_{j}}({\xi_1+\eta_1 }) \widehat{ \wt{\Theta_{j}}}({\xi_1,\eta_1}).
\end{align*}
Repeating the same argument for $\xi_2,\eta_2$, the four cases in \eqref{reduced432} can be completed as
\begin{align}
  A_{12}B_{12}&=\sum_{j,k}\widehat{\varphi_j}(\xi_1)   \widehat{\varphi_k}(\xi_2)    \widehat{\psi_j}(\eta_1)\widehat{\psi_k}(\eta_2) \widehat{\psi_j}(\xi_1+\eta_1)  \widehat{\psi_k}(\xi_2+\eta_2) \widehat{\wt{\Theta_j}}(\xi_1,\eta_1)\widehat{\wt{\Theta_k}}(\xi_2,\eta_2),  \nonumber \\
   A_{12}B_{22}&=\sum_{j,k}\widehat{\varphi_j}(\xi_1)  \widehat{\psi_k}(\xi_2) \widehat{\psi_j}(\eta_1) \widehat{\psi_k}(\eta_2)
\widehat{\psi_j}(\xi_1+\eta_1)\widehat{\varphi_k}(\xi_2+\eta_2) \widehat{\wt{\Theta_j}}(\xi_1,\eta_1)\widehat{\wt{\Theta_k}}(\xi_2,\eta_2),  \nonumber \\
   A_{12}B_{21}&=\sum_{j,k}\widehat{\varphi_j}(\xi_1) \widehat{\psi_k}(\xi_2)  \widehat{\psi_j}(\eta_1) \widehat{\varphi_k}(\eta_2) \widehat{\psi_j}(\xi_1+\eta_1)\widehat{\psi_k}(\xi_2+\eta_2)\widehat{\wt{\Theta_j}}(\xi_1,\eta_1)\widehat{\wt{\Theta_k}}(\xi_2,\eta_2),  \nonumber  \\
   A_{22}B_{22}&= \sum_{j,k}\widehat{\psi_j}(\xi_1)  \widehat{\psi_k}(\xi_2) \widehat{\psi_j}(\eta_1)  \widehat{\psi_k}(\eta_2)\widehat{\varphi_j}(\xi_1+\eta_1) \widehat{\varphi_k}(\xi_2+\eta_2) \widehat{\wt{\Theta_j}}(\xi_1,\eta_1)\widehat{\wt{\Theta_k}}(\xi_2,\eta_2). \label{reduced42}
\end{align}

Finally, from \eqref{Tdecompose}, $T_m(f,g)$ can be decomposed as
\begin{align*}
    & \sum_{j,k\in \z}\int_{\bn\times \bn}\Big( m(\xi,\eta) \widehat{\wt{\Theta_j}}(\xi_1,\eta_1)\widehat{\wt{\Theta_k}}(\xi_2,\eta_2)\Big) \Big(  \widehat{ {\phi}_j}({\xi_1})\widehat{{\phi}_{k}}({\xi_2 }) \wh{f}(\xi) \Big) \\
   & \qq  \times\Big(  \widehat{ {\phi}_{j}}({\eta_1}) \widehat{{\phi}_{k}}({\eta_2 }) \widehat{g}(\eta)\Big) \widehat {{{\phi}}_{j}}({\xi_1+\eta_1}) \widehat {{\phi}}_{k}({\xi_2+\eta_2})  e^{2\pi i \langle x,\xi+\eta\rangle}  d\xi d\eta\\
   &=\sum_{j,k\in \z}\int_{\mathbb{R}^{4n}}\int_{\mathbb{R}^{4n}}K_{j,k}(y_1,y_2,z_1,z_2)   \Big(  \widehat{ {\phi}_j}({\xi_1})\widehat{{\phi}_{k}}({\xi_2 }) \wh{f}(\xi) \Big)                \Big(  \widehat{ {\phi}_{j}}({\eta_1}) \widehat{{\phi}_{k}}({\eta_2 }) \wt{g}(\eta)\Big)  \\
   &  \qq \times  \widehat {{{\phi}}_{j}}({\xi_1+\eta_1}) \widehat {{\phi}}_{k}({\xi_2+\eta_2})  e^{2\pi i (\langle x_1-\frac{y_1}{2^j},\xi_1\rangle+\langle x_2-\frac{y_2}{2^k},\xi_2\rangle+\langle x_1-\frac{z_1}{2^j},\eta_1\rangle+\langle x_2-\frac{z_2}{2^k},\eta_2\rangle )} \;dy dz d\xi d\eta,
  \end{align*}
  where $\phi$ is either a $\varphi$-function or $\psi$-function, taken from the four pieces in \eqref{reduced42}, and
  $\wh{K_{j,k}} (\xi_1,\xi_2,\eta_1,\eta_2):=m(2^j\xi_1,2^k \xi_2,2^j\eta_1,2^k \eta_2) \widehat{\wt{\Theta}}(\xi_1,\eta_1)\widehat{\wt{\Theta}}(\xi_2,\eta_2)$.
 Here, we note that for each pair of $j,k$,
 \begin{align}
  &\Big\|\big(1+4\pi^2(|\cdot_1|^2+|\cdot_3|^2)\big)^{s_1/ 2}\big(1+4\pi^2(|\cdot_2|^2+|\cdot_4|^2)\big)^{s_2/ 2} K_{j,k}\Big\|_{L^{u'}((\bbrn)^4)} \nonumber\\
 & =\Big\|m(2^{j}\cdot_1,2^{k}\cdot_2,2^{j}\cdot_3,2^{k}\cdot_4)\widehat{\wt{\Theta}}(\cdot_1,\cdot_3)\widehat{\wt{\Theta}}(\cdot_2,\cdot_4)\Big\|_{L^u_{s_1,s_2}(\bbr^{2n}\times\bbr^{2n})}\lesssim \mathcal{A}_{s_1,s_2}^{u}[m]\label{estamm}
\end{align}
by the Hausdorff-Young inequality.
 Moreover, an important fact is that in $\widehat{\phi_j}(\xi_1)\widehat{\phi_j}(\eta_1)\widehat{\phi_j}(\xi_1+\eta_1)$, at least two of them are $\Psi$-type functions, and similarly for $\widehat{\phi_k}(\xi_2)\widehat{\phi_k}(\eta_2)\widehat{\phi_k}(\xi_2+\eta_2)$. This will be crucial in the proof, even for the multi-linear cases, which can be seen from the later proof.

Let's first study $A_{12}B_{21}$ part to which we denote by $T_{A_{12}B_{21}}$ the corresponding operator.
Then
\begin{align*}
 & T_{A_{12}B_{21}}\big(f,g\big)(x) \\
   &=\sum_{j,k\in \z}\int_{\mathbb{R}^{4n}}\int_{\mathbb{R}^{4n}}K_{j,k}(y_1,y_2,z_1,z_2) \Big(  \widehat{ {\varphi}_{j}}({\xi_1})\widehat{ {\psi}_{k}}({\xi_2})  \wh{f}(\xi) \Big) \Big(  \widehat{{\psi}_{j}}({\eta_1 })\widehat{{\varphi}_{k}}({\eta_2 }) \wh{g}(\eta)\Big)  \\
   & \q\times   \widehat {{{\psi}}_{j}}({\xi_1+\eta_1}) \widehat {{\psi}}_{k}({\xi_2+\eta_2})  e^{2\pi i (\langle x_1-\frac{y_1}{2^j},\xi_1\rangle+\langle x_2-\frac{y_2}{2^k},\xi_2\rangle+\langle x_1-\frac{z_1}{2^j},\eta_1\rangle+\langle x_2-\frac{z_2}{2^k},\eta_2\rangle )} \; dy dz d\xi d\eta.
\end{align*}
Let
\begin{equation}\label{defd0}
  D_{0}:=\{v\in \bbr^{2n}: |v|\leq 1\}
  \end{equation}
  and
  \begin{equation}\label{defdm}
   D_{M}:=\{v \in \bbr^{2n}: 2^{M-1}< |v|\leq 2^{M}\} \quad\text{ for }~ M\in\bbn.
\end{equation}
Then we decompose $T_{A_{12}B_{21}}$ as
\begin{equation*}
 T_{A_{12}B_{21}}(f,g)= \sum_{\vec{M} \in \mathbb{N}_0\times\bbn_0}T_{A_{12}B_{21}}^{\vec{M}}(f,g)
\end{equation*}
where
\begin{align*}
&T_{A_{12}B_{21}}^{\vec{M}}(f,g)(x_1,x_2)\\
&:=\sum_{j,k\in \z}\int_{\mathbb{R}^{4n}}\int_{\substack{(y_1,z_1)\in D_{M_1}\\(y_2,z_2)\in D_{M_2}     }}K_{j,k}(y_1,y_2,z_1,z_2) \Big(  \widehat{ {\varphi}_{j}}({\xi_1})\widehat{ {\psi}_{k}}({\xi_2})  \wh{f}(\xi) \Big) \Big(  \widehat{{\psi}_{j}}({\eta_1 })\widehat{{\varphi}_{k}}({\eta_2 }) \wh{g}(\eta)\Big)  \\
   &   ~\times \widehat {{{\psi}}_{j}}({\xi_1+\eta_1}) \widehat {{\psi}}_{k}({\xi_2+\eta_2})  e^{2\pi i (\langle x_1-\frac{y_1}{2^j},\xi_1\rangle+\langle x_2-\frac{y_2}{2^k},\xi_2\rangle+\langle x_1-\frac{z_1}{2^j},\eta_1\rangle+\langle x_2-\frac{z_2}{2^k},\eta_2\rangle )}  \;dy dz d\xi d\eta
\end{align*} for $\vec{M}=(M_1,M_2)\in \bbn_0\times \bbn_0$.

In order to prove Theorem~\ref{bl} we will actually prove
\begin{proposition}\label{pp42}
Assume $1<u\le 2$, $0<r<\infty$,  $1<p,q< \infty$ with $\frac{1}{r}=\frac{1}{p}+\frac{1}{q}$. Suppose that
\begin{equation}\label{s1s2condition}
s_1,s_2>\max{\Big(\tf{2n}u,\tfrac{n}{r}-\tfrac{n}{u'}\Big)}.
\end{equation}
   Let $1\leq p_0<\min{(p,u)}$, $1\leq q_0<\min{(q,u)}$. Then for any $\epsilon>0$ there exists a constant $C_\epsilon>0$ such that
  \begin{equation}
  \big\|T_{A_{12}B_{21}}^{\vec{M}}(f,g)\big\|_{L^{r} } \le  C_\epsilon   \mathcal{A}_{s_1,s_2}^{u}[m]2^{-M_1(s_1-{n\over u-\epsilon}-{n\over p_0})} 2^{-M_2(s_2-{n\over u-\epsilon}-{n\over p_0})}
 \|f\|_{L^p} \|g\|_{L^q} \label{g1}
  \end{equation}
  and
  \begin{equation}
  \big\|T_{A_{12}B_{21}}^{\vec{M}}(f,g)\big\|_{L^{r} } \le  C_\epsilon   \mathcal{A}_{s_1,s_2}^u[m]2^{-M_1(s_1-{n\over u-\epsilon}-{n\over q_0})} 2^{-M_2(s_2-{n\over u-\epsilon}-{n\over q_0})}
 \|f\|_{L^p} \|g\|_{L^q}. \label{g2}
  \end{equation}
\end{proposition}
The proof of the proposition will be given in  next section. Now taking the proposition temporarily for granted, let us complete the proof of Theorem \ref{bl}.
\begin{proof}[Proof of Theorem~\ref{bl}]
A similar argument of this reduction for the general multi-linear case can be found in
\cite[Section~4]{LHH}. We will discuss the bilinear case briefly for the sake of completeness.

By the Sobolev embedding theorem, we can suppose $s\leq n+\frac{n}{u}$. The range of $(\tf1p,\tf1q)$ described by \eqref{e001} with $1/r=1/p+1/q$ is the interior of the convex set with vertices
$(0,0)$, $(0,1)$, $(1,0)$, $(1,\tf{s}n-\tf1u)$, and $(\tf{s}n-\tf1u,1)$, where $s=\min\{s_1,s_2\}$. By the multi-linear interpolation, it suffices to verify \eqref{e002} when $(\tf1p,\tf1q)$ are close to these vertices.

To check the boundedness of $T_m$ when $(\tf1p,\tf1q)$ is close to  $(\tf{s}n-\tf1u,1)$, we notice that $\tf sn=\tf sn-\tf1u +\tf1u>\tf1p+\tf1u$.
As we have $p<u$ when $\tf1p$ is close to $\tf{s}n-\tf1u>\tf1{u}$, we can choose $\epsilon$ and $p_0<p$ appropriately so that $\tf sn>\tf1{u-\epsilon}+\tf1{p_0}$. Applying \eqref{g1}, we obtained the decay as $M_1$ and $M_2$ growth, which implies the desired boundedness of $T_m$ when $(\tf1p,\tf1q)$ is close to  $(\tf{s}n-\tf1u,1)$. The other four vertices can be handled similarly, and an interpolation argument, as we mentioned above, yields \eqref{e002} when $p,q,s_1,s_2$ satisfy the condition \eqref{e001}.
\end{proof}

\section{Proof of Proposition ~\ref{pp42}}\label{pop}

We fix $1<u\le 2$ and $s_1,s_2>0$ satisfying \eqref{s1s2condition}, and simply write
$$\mathcal{A}[m]=\mathcal{A}_{s_1,s_2}^u[m]$$
for a notational convenience in this section.
It suffices to establish \eqref{g2} as the other case \eqref{g1} would follow by symmetry.
We shall take advantage of the dual form $$\big\langle T_{A_{12}B_{21}}^{\vec{M}}(f,g), h \big\rangle$$
where we take $h$ from $L^{r'}$-function when $r>1$ and a characteristic function of a certain exceptional set $E'$ when $0<r\leq 1$.
 Before looking at the two cases separately, we first proceed a further reduction of the dual form, which will be useful for both.
  \begin{align}
   &  \big|\big\langle  T_{A_{12}B_{21}}^{\vec{M}}(f,g), h\big\rangle\big|\nonumber\\
   &\le \sum_{j,k\in \z}\Big|\int_{\mathbb{R}^{6n}} \int_{\substack{(y_1,z_1)\in D_{M_1}\\ (y_2,z_2)\in D_{M_2}}} K_{j,k}(y_1,y_2,z_1,z_2) \Big(  \widehat{ {\varphi}_{j}}({\xi_1}) \widehat{ {\psi}_{k}}({\xi_2}) \widehat f(\xi) \Big) \Big(  \widehat{{\psi}_{j}}({\eta_1 })\widehat{{\varphi}_{k}}({\eta_2 }) \widehat g(\eta)\Big) \nonumber \\
   &   ~\times \widehat {{{\psi}}_{j}}({\xi_1+\eta_1}) \widehat {{\psi}}_{k}({\xi_2+\eta_2})  e^{2\pi i (\langle x_1-\frac{y_1}{2^j},\xi_1\rangle+\langle x_2-\frac{y_2}{2^k},\xi_2\rangle+\langle x_1-\frac{z_1}{2^j},\eta_1\rangle+\langle x_2-\frac{z_2}{2^k},\eta_2\rangle )}   h(x)\; dy dz d\xi d\eta dx \Big| \nonumber \\
   &\le \sum_{j,k\in \z}\int_{\mathbb{R}^{2n}} \int_{\substack{(y_1,z_1)\in D_{M_1}\\ (y_2,z_2)\in D_{M_2}}} \big| K_{j,k}(y_1,y_2,z_1,z_2)\big|\, \Big|  \big({ {\varphi}_{j}} \otimes { {\psi}_{k}} \big) \ast f \Big(x_1-{y_1\over 2^j},x_2-{y_2\over 2^k}\Big)\Big|  \label{midmainterm} \\
   &  ~\times  \Big|  \big({ {\psi}_{k}}\otimes  { {\varphi}_{k}} \big) \ast g \Big(x_1-{z_1\over 2^j},x_2-{z_2\over 2^k}\Big)\Big|   \big|\big({{{\psi}}_{j}}\otimes {{\psi}}_{k}\big)\ast h(x_1,x_2)\big|\;  dy dz  dx  \nonumber
      \end{align}
since
$$
\widehat\psi_{j}({\xi_1+\eta_1}) \widehat \psi_{k}({\xi_2+\eta_2}) \int h(x)e^{2\pi ix\cdot (\xi+\eta)}dx=\int (\psi_j\otimes\psi_k)*h(x)e^{2\pi ix\cdot(\xi+\eta)}\;dx.
$$


As in \eqref{defd0} and \eqref{defdm},  the definitions of $D_0$ and $D_M$ for $M\ge 1$ are a little bit different; $D_0$ stands for the unit ball in $\bbr^{2n}$ and $D_M$ does an annulus of radius $2^M$. For convenience of notation in the presentation of the proof we will be only concerned with the case when $M_1, M_2\ge 1$ in which $(y_1,z_1)\in D_{M_1}$ indicates $|(y_1,z_1)|\approx 2^{M_1}$ and $(y_2,z_2)\in D_{M_2}$ does $|(y_2,z_2)|\approx 2^{M_2}$, as a similar and simpler procedure is applicable to the case when $M_1=0$ or $M_2=0$, in which \eqref{am22mm} below will be replaced by
 \begin{align*}
 \Big( \int_{\substack{|(y_1,z_1)|\lesssim 1\\|(y_2,z_2)|\approx 1}}     2^{2M_2n}\big|K_{j,k}(y_1,2^{M_2}y_2,z_1,2^{M_2}z_2)\big|^{u'} ~dy dz        \Big)^{1/u'}&\lesssim \mathcal{A}[m]2^{-M_2(s_2-\frac{2n}u)},\\
  \Big( \int_{\substack{|(y_1,z_1)|\approx 1\\|(y_2,z_2)|\lesssim 1}}     2^{2M_1n}
  \big|K_{j,k}(2^{M_1}y_1,y_2,2^{M_1}z_1,z_2)\big|^{u'}  ~dy dz        \Big)^{1/{u'} }&\lesssim \mathcal{A}[m]2^{-M_1(s_1-\frac{2n}u)},
 \end{align*}
 or
 \begin{equation*}
   \Big( \int_{\substack{|(y_1,z_1)|\lesssim 1\\|(y_2,z_2)|\lesssim 1}}     \big|K_{j,k}(y_1,y_2,z_1,z_2)\big|^{u'}  ~dy dz        \Big)^{1/{u'} }\lesssim \mathcal{A}[m].
 \end{equation*}

 Now we assume that both $M_1,M_2\ge 1$. Then by applying a change of variables $(y_1,z_1)\to (2^{M_1}y_1,2^{M_1}z_1)$, $(y_2,z_2)\to (2^{M_2}y_2,2^{M_2}z_2)$ and $(x_1,x_2)\to (2^{M_1-j}x_1,2^{M_2-k}x_2)$, we write \eqref{midmainterm}  as
   \begin{align*}
       & \sum_{j,k\in \z}\int_{\mathbb{R}^{2n}} \int_{\substack{1/2<|(y_1,z_1)|\le 1\\ 1/2<|(y_2,z_2)|\le 1}} 2^{2M_1n}  2^{2M_2n} | K_{jk}(2^{M_1}y_1,2^{M_2}y_2,2^{M_1}z_1,2^{M_2}z_2)| \\
   & \q\q\times  2^{(-j+M_1)n} 2^{(-k+M_2)n} \Big|  \big({ {\varphi}_{j}} \otimes { {\psi}_{k}} \big) \ast f\Big({x_1-y_1\over 2^{j-M_1}},{x_2-y_2\over 2^{k-M_2}}\Big) \Big|     \\
    &\q\q\q \times \Big|  \big({ {\psi}_{k}}\otimes  { {\varphi}_{k}} \big) \ast g \Big({x_1-z_1\over 2^{j-M_1}},{x_2-z_2\over 2^{k-M_2}}\Big) \Big|\,\Big|\big({{{\psi}}_{j}}\otimes {{\psi}}_{k}\big)\ast h \Big({x_1\over 2^{j-M_1}},{x_2\over 2^{k-M_2}}\Big) \Big| \;dy dz  dx.
    \end{align*}
   Decomposing  $\mathbb{R}^{2n}$ as the product of union of cubes
   $$\mathbb{R}^{2n}=\bigcup_{m_1,m_2\in \z^n} \big(m_1+[0,1)^n\big) \times  \big( m_2+ [0,1)^n\big), $$
   we finally obtain
  \begin{align}
    &\big|\big\langle T_{A_{12}B_{21}}^{\vec{M}}(f,g),h\big\rangle\big|\nonumber \\
    &\lesssim\sum_{j,k\in \z} \sum_{m_1,m_2\in \z^n}\int_{[0,1]^n\times[0,1]^n } \int_{\substack{1/2<|(y_1,z_1)|\le 1\\ 1/2<|(y_2,z_2)|\le 1}} 2^{2M_1n}  2^{2M_2n} | K_{jk}(2^{M_1}y_1,2^{M_2}y_2,2^{M_1}z_1,2^{M_2}z_2)| \nonumber \\
    &\q\times   2^{(-j+M_1)n} 2^{(-k+M_2)n} \Big|  \big({ {\varphi}_{j}} \otimes { {\psi}_{k}} \big) \ast f\Big({m_1+\beta_1-y_1\over 2^{j-M_1}},{m_2+\beta_2-y_2\over 2^{k-M_2}}\Big) \Big|  \nonumber \\
     &\q\q\times\Big|  \big({ {\psi}_{k}}\otimes  { {\varphi}_{k}} \big) \ast g \Big({m_1+\beta_1-z_1\over 2^{j-M_1}},{m_2+\beta_2-z_2\over 2^{k-M_2}}\Big) \Big| \nonumber \\
    &\q\q\q\times\Big|\big({{{\psi}}_{j}}\otimes {{\psi}}_{k}\big)\ast h \Big({m_1+\beta_1\over 2^{j-M_1}},{m_2+\beta_2\over 2^{k-M_2}}\Big) \Big|\; dy dz  d\beta \label{rdend}
   \end{align}

Now we are ready to prove \eqref{g2}.

\hfill

\subsection{The case when $r>1$ and $p>2$}\label{pg2}~

Applying H\"older's inequality for the integration w.r.t $y,z$ in \eqref{rdend}, we can control $ |\langle  T_{A_{12}B_{21}}^{\vec{M}}(f,g), h\rangle|$  by
  \begin{align}
    &\sup_{j,k \in \z} \Big( \int_{\substack{|(y_1,z_1)|\approx 1\\|(y_2,z_2)|\approx 1}}\big| 2^{2M_1n}2^{2M_2n}K_{j,k}(2^{M_1}y_1,2^{M_2}y_2,2^{M_1}z_1,2^{M_2}z_2)\big|^{u'} ~dy dz\Big)^{1/u'}\nonumber\\
    &\q \times \sum_{j,k\in \z} \sum_{m_1,m_2\in \z^n}  2^{(-j+M_1)n} 2^{(-k+M_2)n}  \int_{|\beta_1|,|\beta_2|\le \sqrt{n}} \Big|\big({{{\psi}}_{j}}\otimes {{\psi}}_{k}\big)\ast h \Big({m_1+\beta_1\over 2^{j-M_1}},{m_2+\beta_2\over 2^{k-M_2}}\Big) \Big| \nonumber\\
      & \q\q\times \Big(\int_{{|y_1|,|y_2|\le 1} }\ \Big|  \big({ {\varphi}_{j}} \otimes { {\psi}_{k}} \big) \ast f\Big({m_1+\beta_1-y_1\over 2^{j-M_1}},{m_2+\beta_2-y_2\over 2^{k-M_2}}\Big) \Big|^u dy \Big)^{1/ u}  \nonumber\\
    & \q\q\q\times \Big(\int_{{|z_1|, |z_2|\le 1 } } \Big|  \big({ {\psi}_{k}}\otimes  { {\varphi}_{k}} \big) \ast g \Big({m_1+\beta_1-z_1\over 2^{j-M_1}},{m_2+\beta_2-z_2\over 2^{k-M_2}}\Big) \Big|^u dz\Big)^{1/u}~d\beta. \label{eqpg2}
   \end{align}
We first see that
    \begin{align}\label{am22mm}
    &  \sup_{j,k \in \z} \Big( \int_{\substack{|(y_1,z_1)|\approx 1\\|(y_2,z_2)|\approx 1}}\big| 2^{2M_1n}2^{2M_2n}K_{j,k}(2^{M_1}y_1,2^{M_2}y_2,2^{M_1}z_1,2^{M_2}z_2)\big|^{u'} \;dy dz\Big)^{1/u'}\nonumber\\
    &=2^{-M_1(s_1-\frac{2n}u)}2^{-M_2(s_2-\frac{2n}u)} \sup_{j,k \in \z} \Big( \int_{\substack{|(y_1,z_1)|\approx 2^{M_1}\\|(y_2,z_2)|\approx 2^{M_2}}} \big|2^{M_1s_1}2^{M_2s_2}K_{j,k}(y_1,y_2,z_1,z_2) \big|^{u'}\;dydz\Big)^{1/u'}\nonumber\\
      &\lesssim \mathcal{A}[m] 2^{-M_1(s_1-\frac{2n}u)}2^{-M_2(s_2-\frac{2n}u)}
   \end{align}
   where \eqref{estamm} is applied in the inequality.

   Moreover, applying a simple change of variables, the remaining term is bounded by
\begin{align*}
     &  \sum_{j,k\in \z} \sum_{m_1,m_2\in \z^n}     2^{(-j+M_1)n} 2^{(-k+M_2)n}   \int_{{ |\beta_1|, |\beta_2|\le \sqrt{n} }}\Big|\big({{{\psi}}_{j}}\otimes {{\psi}}_{k}\big)\ast h \Big({m_1+\beta_1\over 2^{j-M_1}},{m_2+\beta_2\over 2^{k-M_2}}\Big) \Big|  d\beta \\
     & \q \times  \Big(\int_{{|y_1|, |y_2|\le 2\sqrt{n}} }\ \Big|  \big({ {\varphi}_{j}} \otimes { {\psi}_{k}} \big) \ast f\Big({m_1+y_1\over 2^{j-M_1}},{m_2+y_2\over 2^{k-M_2}}\Big) \Big|^u ~dy\Big)^{1/ u} \\
   & \q\q\times \Big(\int_{{|z_1|, |z_2|\le 2\sqrt{n} } } \Big|  \big({ {\psi}_{k}}\otimes  { {\varphi}_{k}} \big) \ast g \Big({m_1+z_1\over 2^{j-M_1}},{m_2+z_2\over 2^{k-M_2}}\Big) \Big|^u ~dz\Big)^{1/ u} .
\end{align*}
Therefore, $\big| \big\langle T_{A_{12}B_{21}}^{\vec{M}}(f,g),h\big\rangle\big|$ is bounded by a constant multiple of $\mathcal{A}[m]2^{-M_1(s_1-n)}2^{-M_2(s_2-n)} $ times
\begin{equation}
\sum_{j,k\in \z} \sum_{m_1,m_2\in \z^n} |I^{m_1}_{j-M_1}|\, |J^{m_2}_{k-M_2}| \tilde {\tilde f}_{jk}^{\vec{M}}(\vec{m}) \, \tilde{\tilde g}_{jk}^{\vec{M}}(\vec{m}) \, { \tilde h}_{jk}^{\vec{M}}(\vec{m}) , \label{eq2pg2}
\end{equation}
where we recall
\begin{align*}
I^{m_1}_{j-M_1}&=2^{-(j-M_1)}m_1+[0,2^{-(j-M_1)})^n\\
I^{m_2}_{k-M_2}&=2^{-(k-M_2)}m_2+[0,2^{-(k-M_2)})^n
\end{align*}
and define
\begin{align*}
\tilde {\tilde f}_{jk}^{\vec{M}}(\vec{m}) &:=\Big\Vert      \big({ {\varphi}_{j}} \otimes { {\psi}_{k}} \big) \ast f\Big({m_1+y_1\over 2^{j-M_1}},{m_2+y_2\over 2^{k-M_2}}\Big)   \Big\Vert_{L^u(|y|\le 6\sqrt{n})},\\
\tilde{\tilde g}_{jk}^{\vec{M}}(\vec{m}) &:=\Big\Vert  \big({ {\psi}_{j}}\otimes  { {\varphi}_{k}} \big) \ast g \Big({m_1+y_1\over 2^{j-M_1}},{m_2+y_2\over 2^{k-M_2}}\Big)\Big\Vert_{L^u(|y|\le 6\sqrt{n})},\\
{ \tilde h}_{jk}^{\vec{M}}(\vec{m})&:=\Big\Vert   \big({{{\psi}}_{j}}\otimes {{\psi}}_{k}\big)\ast h \Big({m_1+y_1\over 2^{j-M_1}},{m_2+y_2\over 2^{k-M_2}}\Big)     \Big\Vert_{L^1(|y|\le 6\sqrt{n})}.
\end{align*}

Writing $R^{\vec m}_{j,k,\vec M}=I^{m_1}_{j-M_1} \times J^{m_2}_{k-M_2}$ as before,  \eqref{eq2pg2} can be written as
\begin{align}
 &  \int_{\bn}\sum_{\substack{j,k\in \z\\ m_1,m_2\in\bbzn}} \tilde {\tilde f}_{jk}^{\vec{M}}(\vec{m}) \, \tilde {\tilde g}_{jk}^{\vec{M}}(\vec{m}) \,{\tilde h}_{jk}^{\vec{M}}(\vec{m}) \chi_{R^{\vec m}_{j,k,\vec M}}(x) \;dx  \nonumber\\
 &\lesssim   \int_{\bn}{ \MS}_{\vec M}^{(u)}(f)(x){ \SM}_{\vec M}^{(u)}(g)(x){\DS}_{\vec M}^{(1)}(h)(x) \;dx,
 \label{smdecomposepg2}
 \end{align}
where the hybrid operators are defined as in \eqref{hybridoperator}
with a specific range of $y$ in the $L^1$ or $L^u$ norm rather than simply writing $|y|\lesssim 1$.

Finally, using Lemmas \ref{lmmsl1} and \ref{lemmassmm}, we obtain
\begin{align*}
  &  \big|\big\langle  T_{A_{12}B_{21}}^{\vec{M}}(f,g), h\big\rangle \big|\\
  &\lesssim \mathcal{A}[m]2^{-M_1(s_1-\tf{2n}u)}2^{-M_2(s_2-\tf{2n}u)}    \big\| { \MS}_{\vec M}^{(u)}(f)\big\|_{L^p(\bbr^{2n})}  \big\|{ \SM}_{\vec M}^{(u)}(g)\big\|_{L^q(\bbr^{2n})} \big\|{\DS}_{\vec M}^{(1)}(h)\big\|_{L^{r'}(\bbr^{2n})} \\
        &\lesssim \mathcal{A}[m]  2^{-M_1(s_1-{n\over p_0}-{n\over q_0})} 2^{-M_2(s_2-{n\over p_0}-{n\over q_0})}   \|f\|_{L^p(\bbr^{2n})} \|g\|_{L^q(\bbr^{2n})} \|h\|_{L^{r'}(\bbr^{2n})}
\end{align*}
for all $h\in L^{r'}(\bbr^{2n})$.
Then  \eqref{g2} follows from taking $p_0=u-\epsilon$ for any $\epsilon>0$ due to $p>2\ge u$.

\hfill

\subsection{The case when $r>1$ and $p\le 2$}
\label{pl2}
~\\

Moving back to the reduction \eqref{rdend},  we first make  a change of variables $\beta_1-y_1\to \beta_1$, $\beta_2-y_2\to \beta_2$ to bound $\big| \big\langle T_{A_{12}B_{21}}^{M_1,M_2}(f,g),h\rangle\big|$ by
 \begin{align*}
    &\sum_{j,k\in \z} \sum_{m_1,m_2\in \z^n}\int_{\beta_1,\beta_2\in [-2,2]^n } \int_{\substack{1/2<|(y_1,z_1)|\le 1\\ 1/2<|(y_2,z_2)|\le 1}} 2^{2M_1n}  2^{2M_2n} \big| K_{j,k}(2^{M_1}y_1,2^{M_2}y_2,2^{M_1}z_1,2^{M_2}z_2)\big|\\
    &\q\times   2^{(-j+M_1)n} 2^{(-k+M_2)n} \Big|  \big({ {\varphi}_{j}} \otimes { {\psi}_{k}} \big) \ast f\Big({m_1+\beta_1\over 2^{j-M_1}},{m_2+\beta_2\over 2^{k-M_2}}\Big) \Big|  \\
      &\q\q\times\Big|  \big({ {\psi}_{k}}\otimes  { {\varphi}_{k}} \big) \ast g \Big({m_1+\beta_1+y_1-z_1\over 2^{j-M_1}},{m_2+\beta_2+y_2-z_2\over 2^{k-M_2}}\Big) \Big|  \\
    &\q\q\q\times \Big|\big({{{\psi}}_{j}}\otimes {{\psi}}_{k}\big)\ast h \Big({m_1+\beta_1+y_1\over 2^{j-M_1}},{m_2+\beta_2+y_2\over 2^{k-M_2}}\Big) \Big| ~dy dz  d\beta.
   \end{align*}
 Similar to the previous case,  using H\"older's inequality for the integral over $y,z$, the estimate \eqref{am22mm}, and a change of variables, the above expression is dominated by a constant multiple of $\mathcal{A}[m] 2^{-M_1(s_1-\tf{2n}u)}2^{-M_2(s_2-\tf{2n}u)}$ times
   \begin{equation}\label{eq2}
 \sum_{j,k\in \z}\sum_{m_1,m_2\in \z^n} |I^{m_1}_{j-M_1}|\, |J^{m_2}_{k-M_2}| \tilde f_{jk}^{\vec{M}}(\vec{m}) \, \tilde{\tilde g}_{jk}^{\vec{M}}(\vec{m}) \, \tilde{ \tilde h}_{jk}^{\vec{M}}(\vec{m})
   \end{equation}
   where
   \begin{align}
   \tilde f_{jk}^{\vec{M}}(\vec{m})&:=\Big\Vert     \big({ {\varphi}_{j}} \otimes { {\psi}_{k}} \big) \ast f\Big({m_1+y_1\over 2^{j-M_1}},{m_2+y_2\over 2^{k-M_2}}\Big)     \Big\Vert_{L^1(|y|\le 6\sqrt{n})}, \label{ftildejk}\\
    \tilde{\tilde g}_{jk}^{\vec{M}}(\vec{m})&:= \Big\Vert       \big({ {\psi}_{k}}\otimes  { {\varphi}_{k}} \big) \ast g \Big({m_1+y_1\over 2^{j-M_1}},{m_2+y_2\over 2^{k-M_2}}\Big)   \Big\Vert_{L^u(|y|\le 6\sqrt{n})}, \label{gdtildejk}\\
    \tilde{ \tilde h}_{jk}^{\vec{M}}(\vec{m})&:=\Big\Vert   \big({{{\psi}}_{j}}\otimes {{\psi}}_{k}\big)\ast h \Big({m_1+y_1\over 2^{j-M_1}},{m_2+y_2\over 2^{k-M_2}}\Big)   \Big\Vert_{L^u(y\in [-2,2]^{2n})}. \label{hjkmm}
   \end{align}

Setting $R^{\vec m}_{j,k,\vec M}:=I^{m_1}_{j-M_1} \times J^{m_2}_{k-M_2}$ as above,  \eqref{eq2} can be written as
\begin{align}\label{smdecompose}
 &  \int_{\bn}\sum_{\substack{j,k\in \z\\ m_1,m_2\in\bbzn}} \tilde f_{jk}^{\vec{M}}(\vec{m}) \, \tilde {\tilde g}_{jk}^{\vec{M}}(\vec{m}) \, \tilde {\tilde h}_{jk}^{\vec{M}}(\vec{m}) \chi_{R^{\vec m}_{j,k,\vec M}}(x)~ dx  \nonumber\\
 &\lesssim   \int_{\bn} {\MS}_{\vec M}^{(1)}(f)(x) \, {\SM}_{\vec M}^{(u)}(g)(x)\, { \DS}_{\vec M}^{(u)}(h)(x)~ dx,
 \end{align}

Combining all together and using Lemmas \ref{lmmsl1} and \ref{lemmassmm}, 
we finally have
\begin{align*}
   & \big|\big\langle  T_{A_{12}B_{21}}^{\vec{M}}(f,g), h\big\rangle\big|\\
  &\lesssim \mathcal{A}[m] 2^{-M_1(s_1-\tf{2n}u)}2^{-M_2(s_2-\tf{2n}u)}   \big\| {\MS}_{\vec M}^{(1)}(f)\big\|_{L^p(\bbr^{2n})}  \big\|{\SM}_{\vec M}^{(u)}(g)\big\|_{L^q(\bbr^{2n})} \big\|{ \DS}_{\vec M}^{(u)}(h)\big\|_{L^{r'}(\bbr^{2n})} \\
     &\lesssim \mathcal{A}[m] 2^{-M_1(s_1-\tf{2n}u)}2^{-M_2(s_2-\tf{2n}u)}    2^{(M_1+M_2)n({1\over q_0}-{1\over u})} 2^{(M_1+M_2)n({1\over r'_0}-{1\over u})} \\
     &\qq\qq\qq\qq\qq\qq\qq\qq\times \|f\|_{L^p(\bbr^{2n})} \|g\|_{L^q(\bbr^{2n})} \|h\|_{L^{r'}(\bbr^{2n})} \\
        &= \mathcal{A}[m]2^{-M_1(s_1-{n\over q_0}-{n\over r'_0})} 2^{-M_2(s_2-{n\over q_0}-{n\over r'_0})}   \|f\|_{L^p(\bbr^{2n})} \|g\|_{L^q(\bbr^{2n})} \|h\|_{L^{r'}(\bbr^{2n})}
\end{align*}
for $h\in L^{r'}(\bbr^{2n})$ and all $1<r'_0<\min{(u,r')}$.

Since $1<r<p\le 2$, we have $r'>2\ge u$ and thus \eqref{g2} follows from choosing $r'_0=u-\epsilon$.

\hfill

\subsection{The case when $0<r\leq 1$}
~\\
In order to verify \eqref{g2} when $0<r\leq 1$, we will actually prove the following estimate:\\
Let $\epsilon>0$ be arbitrary. Then
\begin{equation}\label{e020}
\big\|T_{A_{12}B_{21}}^{\vec{M}}(f,g)\big\|_{L^{r,\infty}(\bbr^{2n}) } \lesssim_{\epsilon}  \mathcal{A}[m]2^{-M_1(s_1-{n\over u-\epsilon}-{n\over q_0})} 2^{-M_2(s_2-{n\over u-\epsilon}-{n\over q_0}) } \|f\|_{L^p(\bbr^{2n})} \|g\|_{L^q(\bbr^{2n})}
\end{equation}
for 
$1\leq q_0<\min{(q,u)}$, where we recall that  $L^{r,\infty}(\bbr^{2n})$ ($0<r<\infty$) is the weak-$L^r$ space with norm
$$\|f\|_{L^{r,\infty}(\bbr^{2n})}=\sup_{t>0} t \big| \big\{x\in \bbr^{2n}: |f(x)|>t \big\}\big|^{{1}/{r}}.$$
The quasi-Banach space $L^{r,\infty}$ does not have a good dual space, but we have the following substitute which is enough for our purpose. One may refer to
\cite[Lemma 2.6]{muscalu2013classical} for its proof.
\begin{lemma}\label{weak}
  Let $0<r\leq 1$ and $A>0$. The following two statements are equivalent:

(i)  $\|f\|_{L^{r,\infty}(\bbr^{2n})}\leq A$,

(ii) for every set $E$ in $\bbr^{2n}$ with $0<|E|<\infty$, there exists a measurable subset $E'\subseteq E$ with $|E'|\approx |E|$ and $|\langle f,\chi_{E'}\rangle|\lesssim A |E|^{1/r'}$, where $1/r+1/r'=1$.
\end{lemma}

Let $E$ be any subset of  $\bbr^{2n}$. Then, by Lemma \ref{weak} and scaling invariance, it suffices to prove that for $\|f\|_{L^p(\bbr^{2n})}=\|g\|_{L^q(\bbr^{2n})}=1=|E|$, there exists a subset $E'\subseteq E$ with $|E'|\approx 1$, such that
  \begin{equation}\label{wkest}
\big|\big\langle  T_{A_{12}B_{21}}^{\vec{M}}(f,g), \chi_{E'}\big\rangle \big| \lesssim  \mathcal{A}[m]2^{-M_1(s_1-{n\over u-\epsilon}-{n\over q_0})} 2^{-M_2(s_2-{n\over u-\epsilon}-{n\over q_0})}. 
  \end{equation}
Letting $h=\chi_{E'}$ with $E'$ to be determined later,
in order to control $|\langle  T_{A_{12}B_{21}}^{\vec{M}}(f,g), h\rangle| $, we will take advantage of Lemma~\ref{compact} to write
\begin{align*}
 &  \big({ {\psi}_{j}} \otimes { {\psi}_{k}} \big) \ast h\Big({m_1+y_1\over 2^{j-M_1}},{m_2+y_2\over 2^{k-M_2}}\Big) \nonumber \\
 &=\sum_{ \mu_1, \mu_2 \in \mathbb{N}} 2^{-1000n(\mu_1+\mu_2)} \big({ \psi^{\mu_1}_{j}} \otimes { \psi^{\mu_2}_{k}} \big) \ast h\Big({m_1+y_1\over 2^{j-M_1}},{m_2+y_2\over 2^{k-M_2}}\Big).  \nonumber
\end{align*}

As a result, $\big| \big\langle  T_{A_{12}B_{21}}^{\vec{M}}(f,g), h\big\rangle\big|$,
via \eqref{eq2}, can be estimated by
\begin{equation}\label{mudecompose}
\begin{aligned}
 & \mathcal{A}[m] 2^{-M_1(s_1-\tf{2n}u)}2^{-M_2(s_2-\tf{2n}u)}\sum_{\mu_1,\mu_2 \in \mathbb{N}} 2^{-1000n( \mu_1+\mu_2 ) }      \sum_{j,k\in \z} \sum_{m_1,m_2\in \z^n}  \\
 &  \qq\qq \times \int_{\bn}  \tilde f_{jk}^{\vec{M}}(\vec{m}) \, \tilde {\tilde g}_{jk}^{\vec{M}}(\vec{m})   \tilde{\tilde h}_{jk}^{\vec{M},\vec \mu}(\vec{m})  \chi_{R^{\vec m}_{j,k,\vec M}}(x)~ dx
 \end{aligned}
  \end{equation}
  where $ \tilde f_{jk}^{\vec{M}}(\vec{m})$ and $\tilde {\tilde g}_{jk}^{\vec{M}}(\vec{m})$ are defined as in \eqref{ftildejk} and \eqref{gdtildejk}, and
  $$\tilde{\tilde h}_{jk}^{\vec{M},\vec \mu}(\vec{m}) :=\Big\Vert  \big({{{\psi}}_{j}}^{\mu_1}\otimes {{\psi}}_{k}^{\mu_2}\big)\ast h \Big({m_1+y_1\over 2^{j-M_1}},{m_2+y_2\over 2^{k-M_2}}\Big)    \Big\Vert_{L^u(y\in [-2,2]^{2n})}$$
  for $\vec\mu=(\mu_1,\mu_2)$, similar to \eqref{hjkmm}.
  The sum over $j,k\in\bbz$ and $m_1,m_2\in\bbzn$ can be thought of as the sum over any dyadic rectangles $R_{j,k,\vec{M}}^{\vec m}$. Let $\RR$ be any finite collections of dyadic rectangles.

  We now claim that there exists a constant $N<1000n$ such that for any fixed $\mu_1,\mu_2\in\bbn$,

\begin{align}
&\sum_{R_{j,k,\vec{M}}^{\vec m}\in\RR} \int_{\bn}\tilde f_{jk}^{\vec{M}}(\vec{m}) \, \tilde {\tilde g}_{jk}^{\vec{M}}(\vec{m})   \tilde{\tilde h}_{jk}^{\vec{M},\vec \mu}(\vec{m})  \chi_{R^{\vec m}_{j,k,\vec M}}(x)~dx\nonumber\\
&\lesssim 2^{N(\mu_1+\mu_2)}2^{(M_1+M_2)n({1\over q_0}-{1\over u})}  2^{(M_1+M_2)n({1\over u-\epsilon}-{1\over u})}\label{e021}
\end{align}
uniformly in $\RR$, and then \eqref{wkest} follows immediately.\\

To verify \eqref{e021}, we will follow the strategy in \cite{MPTT1,MPTT2} to define $E'$ carefully with the help of the boundedness of hybrid operators established in Lemmas~\ref{lmmsl1}, ~\ref{lemmassmm}, and ~\ref{lmnoncompactbip}, which allows us to split the left-hand side of \eqref{e021} into two terms. The main term can be estimated using a standard Calder\'on-Zygmund argument, while the remaining term, saying ``tail part", is relatively  simpler. In fact, the tail part vanishes due to the definition of $E'$. Let us describe the details below.

Let $C>1$ be a sufficiently large number to be chosen shortly and we define the exceptional sets
\begin{align*}
 \Omega_{\mmu,0}^f&:=\Big\{x\in \bn:\MS_{\vec{M}}^{(1)}(f)(x)>C 2^{10n(\mu_1+\mu_2)}\Big\} \\
 \Omega_{\mmu,0}^g&:=\Big\{x\in \bn: \SM_{\vec{M}}^{(u)}(g)(x)>C 2^{(M_1+M_2)n({1\over q_0}-{1\over u})} 2^{10n ( \mu_1+\mu_2)}\Big\}\\
  \Omega_{\vec{\mu}}&:=\Omega_{\mmu,0}^f\cup \Omega_{\mmu,0}^g\\
  \widetilde{\Omega}_{\vec  \mu}&:=\Big\{x\in \bn: \mathcal{M}^{(n),\otimes}\big( \chi_{\Omega_{{\vec \mu}}}\big)(x)>{1 \over  100}\Big\}, \\
  \widetilde{\widetilde{\Omega}}_{\vec \mu}&:=\Big\{x\in \bn: \mathcal{M}^{(n),\otimes}\big( \chi_{\widetilde{\Omega}_{{\vec \mu}}}\big)(x)>{ 2^{-2n(\mu_1+\mu_2)-8n}}\Big\},
\end{align*}
and
$$  \Omega:= \bigcup_{\vec \mu \in \mathbb{N}^2}\widetilde{\widetilde{\Omega}}_{\vec \mu}.$$
Choose $1<t<\min\{p,q\}$ and then we see that
\begin{align*}
|\Omega|&\le \sum_{\mu_1,\mu_2\in\bbn}\big| \widetilde{\widetilde{\Omega}}_{\vec \mu} \big|\le \sum_{\mu_1,\mu_2\in\bbn} 2^{(2n(\mu_1+\mu_2)+8n)t}\big|\widetilde{\Omega}_{\vec  \mu} \big|\lesssim \sum_{\mu_1,\mu_2\in\bbn} 2^{2nt(\mu_1+\mu_2)}\big| \Omega_{\vec{\mu}}\big|
\end{align*}
by using Chebyshev's inequallity and \eqref{fsinq0}.
Moreover, using  Lemmas \ref{lmmsl1} and \ref{lemmassmm},
\begin{align*}
\big| \Omega_{\vec{\mu}}\big|&\le \Big| \Big\{x\in \bn:\MS_{\vec{M}}^{(1)}(f)(x)>C 2^{10n(\mu_1+\mu_2)}\Big\}\Big| \\
   &\qq +  \Big|\Big\{x\in \bn: \SM_{\vec{M}}^{(u)}(g)(x)>C 2^{(M_1+M_2)n({1\over q_0}-{1\over u})} 2^{10n ( \mu_1+\mu_2)}\Big\}\Big|\\
   &\le C^{-p}2^{-10n(\mu_1+\mu_2)p}\big\Vert \MS_{\vec{M}}^{(1)}(f) \big\Vert_{L^p(\bbr^{2n})}^p\\
   &\qq +C^{-q}2^{-q(M_1+M_2)n(1/q_0-1/u)}2^{-10n(\mu_1+\mu_2)q}\big\Vert  \SM_{\vec{M}}^{(u)}(g)    \big\Vert_{L^q(\bbr^{2n})}^q\\
   &\lesssim C^{-t}2^{-10n(\mu_1+\mu_2)t},
\end{align*}
which shows that
$$|\Omega|\lesssim C^{-t}\sum_{\mu_1,\mu_2\in\bbn}2^{-8n(\mu_1+\mu_2)t}\lesssim C^{-t}.$$
Therefore, we may pick $C>1$ so large that $|\Omega|\le \frac{1}{10}$ and then we define
$$E':=E\setminus \Omega.$$
Clearly, \begin{equation}\label{este'}
\frac{9}{10}\le |E'|\le 1
\end{equation}
and thus it remains to verify \eqref{e021} with $h:=\chi_{E'}$.

\hfill

\begin{proof}[Proof of {\eqref{e021}}]
We first claim that the sum in the left-hand side of \eqref{e021} can be replaced by the sum over $R_{j,k,\vec{M}}^{\vec m}\in \RR$ satisfying
\begin{equation*}
R_{j,k,\vec{M}}^{\vec m}\cap (\wt{\Omega}_{\vec{\mu}})^c\not= \emptyset.
\end{equation*}
To verify this, suppose that
\begin{equation}\label{assrjkmm}
R_{j,k,\vec{M}}^{\vec m}\cap (\wt{\Omega}_{\vec{\mu}})^c= \emptyset,
\end{equation}
which is equivalent to
$$R_{j,k,\vec{M}}^{\vec m}\subset \wt{\Omega}_{\vec{\mu}}.$$
If $x\in \big( 2^{\mu_1+4}I_{j-M_1}^{m_1}\big)\times \big( 2^{\mu_2+4}I_{k-M_2}^{m_2}\big)$, then
\begin{align*}
\mathcal{M}^{(n),\otimes}\big(\chi_{\wt{\Omega}_{\vec \mu}}\big)(x)&\ge \frac{1}{|2^{\mu_1+4}I_{j-M_1}^{m_1}|\,|2^{\mu_2+4}I_{k-M_2}^{m_2}|}\int_{( 2^{\mu_1+4}I_{j-M_1}^{m_1})\times ( 2^{\mu_2+4}I_{k-M_2}^{m_2})} \chi_{\wt{\Omega}_{\mmu}}(x)~dx\\
&\ge 2^{-(\mu_1+\mu_2)n-8n}\frac{1}{|I_{j-M_1}^{m_1}|\,|J_{k-M_2}^{m_2}|}\big| R_{j,k,\vec{M}}^{\vec{m}}\cap \wt{\Omega}_{\mmu} \big|=2^{-(\mu_1+\mu_2)n-8n},
\end{align*}
which implies $x\in \widetilde{\widetilde{\Omega}}_{\vec \mu}$. That is, by the definition of $\widetilde{\widetilde{\Omega}}_{\vec \mu}$,
$$ \big( 2^{\mu_1+4}I_{j-M_1}^{m_1}\big)\times \big( 2^{\mu_2+4}I_{k-M_2}^{m_2}\big)\subset \widetilde{\widetilde{\Omega}}_{\vec \mu}.$$
On the other hand,
\begin{align}
&\Big|  \big({{{\psi}}_{j}}^{\mu_1}\otimes {{\psi}}_{k}^{\mu_2}\big)\ast h \Big({m_1+y_1\over 2^{j-M_1}},{m_2+y_2\over 2^{k-M_2}}\Big)  \Big|\nonumber\\
&\le \int_{E'}\Big|\psi_j^{\mu_1}\Big(\frac{m_1+y_1}{2^{j-M_1}}-z_1 \Big) \Big|\, \Big|\psi_k^{\mu_2}\Big( \frac{m_2+y_2}{2^{k-M_2}}-z_2\Big) \Big| |h(z_1,z_2)|~dz_1 dz_2 \nonumber\\
&\le \int_{E\setminus \wt{\wt{\Omega}}_{\mmu} }\Big|\psi_j^{\mu_1}\Big(\frac{m_1+y_1}{2^{j-M_1}}-z_1 \Big) \Big|\, \Big|\psi_k^{\mu_2}\Big( \frac{m_2+y_2}{2^{k-M_2}}-z_2\Big) \Big| |h(z_1,z_2)|~dz_1 dz_2.\label{intewtwtom}
\end{align}
We note that the function $\psi$, which was taken at the beginning of Section \ref{reductionsection}, is adapted to $[-1,1]^n$ and thus
\begin{equation}\label{myzmuj}
\frac{m_1+y_1}{2^{j-M_1}}-z_1\in 2^{\mu_1-j}[-1,1]^n\q \text{ and }\q \frac{m_2+y_2}{2^{k-M_2}}-z_2\in 2^{\mu_2-k}[-1,1]^n
\end{equation}
as $Q=[-1,1]^n$ in \eqref{suppphi_k}. Then \eqref{myzmuj} implies that
$$z_1\in 2^{-j+M_1}m_1+2^{-j+\mu_1+M_1}[-3,3]^n\subset 2^{\mu_1+4}I_{j-M_1}^{m_1},$$
$$z_2\in 2^{-k+M_2}m_2+2^{-k+\mu_2+M_2}[-3,3]^n\subset 2^{\mu_2+4}J_{k-M_2}^{m_2}$$
and thus $z\in \wt{\wt{\Omega}}_{\mmu}$
 for $y=(y_1,y_2)\in [-2,2]^{n}\times [-2,2]^{n}$. This yields that the right-hand side of \eqref{intewtwtom} vanishes as the integral is taken over $z\in E\setminus \wt{\wt{\Omega}}_{\mmu}$. This means that
 $$\tilde{\tilde h}_{jk}^{\vec{M},\vec \mu}(\vec{m})=0 \q \text{if \eqref{assrjkmm} holds},$$
as desired.\\

As a result, the proof of \eqref{e021} can be reduced to
\begin{align}
&\sum_{\substack{R_{j,k,\vec{M}}^{\vec m}\in\RR\\  R_{j,k,\vec{M}}^{\vec m}\cap (\wt{\Omega}_{\vec{\mu}})^c\not= \emptyset    }} \int_{\bn}\tilde f_{jk}^{\vec{M}}(\vec{m}) \, \tilde {\tilde g}_{jk}^{\vec{M}}(\vec{m})   \tilde{\tilde h}_{jk}^{\vec{M},\vec \mu}(\vec{m})  \chi_{R^{\vec m}_{j,k,\vec M}}(x)~dx\nonumber\\
&\lesssim 2^{(N+20n)(\mu_1+\mu_2)}2^{(M_1+M_2)n({1\over q_0}-{1\over u})}  2^{(M_1+M_2)n({1\over u-\epsilon}-{1\over u})}\label{e021reduction}
\end{align}
uniformly in $\RR$.
If $x\in R\cap (\wt{\Omega}_{\vec{\mu}})^c\not= \emptyset $, then, by the definition of $\wt{\Omega}_{\mmu}$, we have
\begin{equation*}
\frac{| R\cap \Omega_{\mmu}|}{|R|}\le \mathcal{M}^{(n),\otimes}\big(\chi_{\Omega_{\mmu}}\big)(x)< \frac{1}{100},
\end{equation*}
which yields
$$| R\cap \Omega_{\mmu}|< \frac{1}{100}|R|.$$
This further implies that
\begin{equation}\label{rcapomegaest}
|R\cap \Omega_{\mmu,0}^f|<\frac{1}{100}|R| \q\text{ and }\q |R\cap \Omega_{\mmu,0}^g|<\frac{1}{100}|R|  \q \text{ if }~ R\cap (\wt{\Omega}_{\vec{\mu}})^c\not= \emptyset.
\end{equation}

Let $C>1$ be the constant that appeared in the definition of $\Omega$ so that \eqref{este'} holds.
We define
\begin{equation*}
\Omega_{\mmu,\kappa}^f:=\Big\{x\in \bbr^{2n}: \MS_{\vec M}^{(1)}(f)(x)>2^{-\kappa}C 2^{10n (\mu_1+\mu_2)} \Big\}.
\end{equation*}
Note that $ \Omega_{\mmu,1}^f\subset \Omega_{\mmu,2}^f\subset \cdots $ and
$$\bigcup_{\kappa\ge 1}\Omega_{\mmu,\kappa}^f=\big\{x\in\bbr^{2n}: \MS_{\vec M}^{(1)}(f)(x)\not=0\big\}.$$
Now we define disjoint partitions of $\RR$ as
\begin{equation*}
\RR_{\mmu,\kappa}^{f}:=\begin{cases}
\Big\{ R\in\RR: \big| R\cap \Omega_{\mmu,1}^{f}\big|>\frac{1}{100}|R|\Big\}, & \kappa=1\\
\Big\{ R\in\RR\setminus \cup_{j=1}^{\kappa-1}\RR_{\mmu,j}^{f}:   \big| R\cap \Omega_{\mmu,\kappa}^{f}\big|>\frac{1}{100}|R|       \Big\}, & \kappa\ge 2
\end{cases}
\end{equation*}
so that
\begin{equation}\label{rpart1}
\RR= \dot{\bigcup_{\kappa\ge 1}}\RR_{\mmu,\kappa}^{f}.
\end{equation}

Similarly, for each $\kappa\in \bbn$,
\begin{equation*}
\Omega_{\mmu,\kappa}^g:=\Big\{x\in \bbr^{2n}: \SM_{\vec M}^{(u)}(g)(x)>2^{-\kappa}C 2^{10n (\mu_1+\mu_2)} 2^{(M_1+M_2)n(\frac{1}{q_0}-\frac{1}{u})}\Big\}
\end{equation*}
and
\begin{equation*}
\RR_{\mmu,\kappa}^{g}:=\begin{cases}
\Big\{ R\in\RR: \big| R\cap \Omega_{\mmu,1}^{g}\big|>\frac{1}{100}|R|\Big\}, & \kappa=1\\
\Big\{ R\in\RR\setminus \cup_{j=1}^{\kappa-1}\RR_{\mmu,j}^{g}:   \big| R\cap \Omega_{\mmu,\kappa}^{g}\big|>\frac{1}{100}|R|       \Big\}, & \kappa\ge 2
\end{cases}.
\end{equation*}
Due to the structure similar to $\Omega_{\mmu,\kappa}^f$ and $\RR_{\mmu,\kappa}^{f}$, it follows that
$$\bigcup_{\kappa\ge 1}\Omega_{\mmu,\kappa}^g=\big\{x\in\bbr^{2n}: \SM_{\vec M}^{(u)}(g)(x)\not=0\big\}.$$
and
\begin{equation}\label{rpart2}
\RR= \dot{\bigcup_{\kappa\ge 1}}\RR_{\mmu,\kappa}^{g}.
\end{equation}

Let $N<1000n$ and we additionally define
\begin{equation*}
\Omega_{\mmu,\kappa}^{h}:=\Big\{x\in\bbr^{2n} : \DS_{\vec{M}}^{(u),\mmu}(h)(x)>2^{-\kappa}C 2^{(M_1+M_2)n(\frac{1}{u-\epsilon}-\frac{1}{u})}2^{N(\mu_1+\mu_2)} \Big\}
\end{equation*}
for $\kappa\in\bbz$.
Then it follows from Chebyshev's inequality and Lemma \ref{lmnoncompactbip} that
\begin{align*}
\big|\Omega_{\mmu,\kappa}^{h} \big|&\le \big( 2^{-\kappa}C 2^{(M_1+M_2)n(\frac{1}{u-\epsilon}-\frac{1}{u})} \big)^{-u}\big\Vert \DS_{\vec{M}}^{(u),\mmu}(h)\big\Vert_{L^u(\bbr^{2n})}^u\\
&\lesssim 2^{u\kappa}C^{-u}\Vert h\Vert_{L^u(\bbr^{2n})}^u\\
&\lesssim 2^{u\kappa}C^{-u}\to 0\q \text{as }~\kappa\to -\infty.
\end{align*}
Since $\RR$ is a finite set, we may choose a large integer $\NN=\NN_{\RR}$, depending on $\RR$, such that
\begin{equation}\label{rinteroest}
\big| R\cap \Omega_{\mmu,-\NN}^h\big|<\frac{1}{100}|R|\q \text{ for }~ R\in\RR.
\end{equation}
We now construct disjoint partitions of $\RR$ as follows:
\begin{equation*}
\RR_{\mmu,\kappa}^{h}:=\begin{cases}
\Big\{ R\in\RR: \big| R\cap \Omega_{\mmu,-\NN+1}^{h}\big|>\frac{1}{100}|R|\Big\}, & \kappa=-\NN+1\\
\Big\{ R\in\RR\setminus \cup_{j=-\NN+1}^{\kappa-1}\RR_{\mmu,j}^{h}:   \big| R\cap \Omega_{\mmu,\kappa}^{h}\big|>\frac{1}{100}|R|       \Big\}, & \kappa\ge -\NN+2
\end{cases}.
\end{equation*}
for $\kappa\ge -\NN+1$.
Clearly, we also have
$$\bigcup_{\kappa\ge -\NN+1}\Omega_{\mmu,\kappa}^h=\big\{x\in\bbr^{2n}: \DS_{\vec M}^{(u),\mmu}(h)(x)\not=0\big\}.$$
and
\begin{equation}\label{rpart3}
\RR= \dot{\bigcup_{\kappa\ge -\NN+1}}\RR_{\mmu,\kappa}^{h}.
\end{equation}

We observe that
\begin{equation}\label{omegainter1}
\big| R\cap \Omega_{\mmu,\kappa-1}^{f}\big|<\frac{1}{100}|R| \q \text{ if }~ R\cap (\wt{\Omega}_{\mmu})^c\not=\emptyset \text{ and } R\in \RR_{\mmu,\kappa}^f ~\text{ for }~\kappa\ge 1,
\end{equation}
which follows from \eqref{rcapomegaest} for $\kappa=1$ and from the definition of $\RR_{\mmu,\kappa}^f$ for $\kappa\ge 2$.
Similarly, due to \eqref{rcapomegaest} and the definition of $\RR_{\mmu,\kappa}^g$,
\begin{equation}\label{omegainter2}
\big| R\cap \Omega_{\mmu,\kappa-1}^{g}\big|<\frac{1}{100}|R| \q \text{ if }~ R\cap (\wt{\Omega}_{\mmu})^c\not=\emptyset \text{ and }  R\in \RR_{\mmu,\kappa}^g ~\text{ for }~\kappa\ge 1,
\end{equation}
 and lastly,
 \begin{equation}\label{omegainter3}
\big| R\cap \Omega_{\mmu,\kappa-1}^h\big|<\frac{1}{100}|R|\q \text{ if }~ R\in \RR_{\mmu,\kappa}^{h} ~\text{ for }~\kappa\ge -\NN+1
\end{equation}
 according to \eqref{rinteroest} and the definition of $\RR_{\mmu,\kappa}^h$.

\hfill

In view of \eqref{rpart1}, \eqref{rpart2}, and \eqref{rpart3}, letting
$$\mathcal{R}_{\kappa_1,\kappa_2,\kappa_3}:=\RR_{\mmu,\kappa_1}^{f}\cap \RR_{\mmu,\kappa_2}^{g}\cap \RR_{\mmu,\kappa_3}^{h},$$
the left-hand side of \eqref{e021reduction} can be written as
\begin{equation}
  \sum_{\substack{\kappa_1,\kappa_2\ge 1\\ \kappa_3\ge-\NN+1}} \sum_{\substack{R_{j,k,\vec{M}}^{\vec m}\in \mathcal{R}_{\kappa_1,\kappa_2,\kappa_3}\\R_{j,k,\vec{M}}^{\vec m}\cap(\wt{\Omega}_{\mmu})^c\not=\emptyset }} |R^{\vec m}_{j,k,\vec M}|\tilde f_{jk}^{M_1,M_2}(\vec m) \, \tilde {\tilde g}_{jk}^{M_1,M_2}(\vec m)   \tilde{\tilde h}_{jk}^{M_1,M_2,\vec \mu}(\vec m).
  \label{mudecompose3}
\end{equation}
Using \eqref{omegainter1}, \eqref{omegainter2}, and \eqref{omegainter3}, we have
$$\big| R\cap (\Omega_{\mmu,\kappa_1-1}^{f})^c\cap (\Omega_{\mmu,\kappa_2-1}^{g})^c\cap (\Omega_{\mmu,\kappa_3-1}^{h})^c \big|
 \ge \frac{97}{100}|R|$$
if $R\cap (\wt{\Omega}_{\mmu})^c\not= \emptyset$ and
$R\in\RR_{\kappa_1,\kappa_2,\kappa_3}$ for $\kappa_1,\kappa_2\ge 1$ and $\kappa_3\ge -\NN+1$.
Therefore, \eqref{mudecompose3} is bounded by
\begin{align*}
&\frac{100}{97}\sum_{\substack{\kappa_1,\kappa_2\ge 1\\ \kappa_3\ge-\NN+1}} \sum_{\substack{R_{j,k,\vec{M}}^{\vec m}\in \mathcal{R}_{\kappa_1,\kappa_2,\kappa_3}}}\tilde f_{jk}^{\vec{M}}(\vec m) \, \tilde {\tilde g}_{jk}^{\vec{M}}(\vec m)   \tilde{\tilde h}_{jk}^{\vec{M},\vec \mu}(\vec m)\\
&\qq\qq\qq\qq\qq\times \big| R_{j,k,\vec{M}}^{\vec m}\cap (\Omega_{\mmu,\kappa_1-1}^{f})^c\cap (\Omega_{\mmu,\kappa_2-1}^{g})^c\cap (\Omega_{\mmu,\kappa_3-1}^{h})^c \big|\\
&\approx \sum_{\substack{\kappa_1,\kappa_2\ge 1\\ \kappa_3\ge-\NN+1}} \int_{ (\Omega_{\mmu,\kappa_1-1}^f)^c \cap (\Omega_{\mmu,\kappa_2-1}^g)^c\cap (\Omega_{\mmu,\kappa_3-1}^h)^c}\sum_{\substack{R_{j,k,\vec{M}}^{\vec m}\in \mathcal{R}_{\kappa_1,\kappa_2,\kappa_3}}}  \tilde f_{jk}^{\vec{M}}(\vec{m})\\
&\qq\qq\qq\qq\qq\qq\qq\qq\qq\qq\times \tilde{\tilde g}_{jk}^{\vec{M}}(\vec{m}) \tilde{\tilde h}_{jk}^{\vec{M},\vec \mu}(\vec{m})  \chi_{R_{j,k,\vec{M}}^{\vec{m}}}(x)~dx  \\
   &\lesssim\sum_{\substack{\kappa_1,\kappa_2\ge 1\\ \kappa_3\ge-\NN+1}} \int_{ (\Omega_{\mmu,\kappa_1-1}^f)^c \cap (\Omega_{\mmu\kappa_2-1}^g)^c\cap (\Omega_{\mmu,\kappa_3-1}^h)^c\cap \Omega_{\RR_{\kappa_1,\kappa_2,\kappa_3}}}\MS_{\vec{M}}^{(1)}(f)(x)\\
   &\qq\qq\qq\qq\qq\qq\qq\qq\qq\qq\times  \SM_{\vec{M}}^{(u)}(g)(x) \DS_{\vec{M}}^{(u),\mmu}(h) (x) ~dx \\
   &\lesssim 2^{(N+20n)(\mu_1+\mu_2)}2^{(M_1+M_2)n({1\over q_0}-{1\over u})}  2^{(M_1+M_2)n(\frac{1}{u-\epsilon}-{1\over u})} \sum_{\substack{\kappa_1,\kappa_2\ge 1\\ \kappa_3\ge -\NN+1}} 2^{-\kappa_1}2^{-\kappa_2}2^{-\kappa_3}\big|\Omega_{\mathcal{R}_{\kappa_1,\kappa_2,\kappa_3}}\big|
\end{align*}
where $ \Omega_{\mathcal{R}_{\kappa_1,\kappa_2,\kappa_3}}:= \bigcup_{R\in \mathcal{R}_{\kappa_1,\kappa_2,\kappa_3}} R$.
Now we claim that
\begin{equation}\label{lastclaim}
\sum_{\substack{\kappa_1,\kappa_2\ge 1\\ \kappa_3\ge -\NN+1}} 2^{-\kappa_1}2^{-\kappa_2}2^{-\kappa_3}\big|\Omega_{\mathcal{R}_{\kappa_1,\kappa_2,\kappa_3}}\big|\lesssim 1,
\end{equation}
which clearly concludes \eqref{e021reduction}.

To verify \eqref{lastclaim}, we first observe that if $R\in\RR_{\mmu,\kappa_1}^{f}$ and $x\in R$, then
$$\mathcal{M}^{(n),\otimes}\big(\chi_{\Omega_{\mmu,\kappa_1}^{f}}\big)(x)\ge \frac{\big|R\cap \Omega_{\mmu,\kappa_1}^{f} \big|}{|R|}>\frac{1}{100},$$
and this implies that
$$\bigcup_{R\in \mathcal{R}_{\mmu,\kappa_1}^f} R\subset \Big\{ x\in\bbr^{2n}: \mathcal{M}^{(n),\otimes}\big( \chi_{\Omega_{\mmu,\kappa_1}^{f}}\big)(x)>\frac{1}{100}\Big\}.$$
Therefore, using Chebyshev's inequality, \eqref{fsinq0}, and Lemma \ref{lmmsl1}, we have
\begin{align*}
\big| \Omega_{\RR_{\kappa_1,\kappa_2,\kappa_3}}\big|&\le \bigg|\bigcup_{R\in\RR_{\mmu,\kappa_1}^f}R\,\bigg| \le \Big|\Big\{ x\in\bbr^{2n}: \mathcal{M}^{(n),\otimes}\big( \chi_{\Omega_{\mmu,\kappa_1}^{f}}\big)(x)>\frac{1}{100}\Big\}\Big|\\
  &\lesssim \big\Vert \mathcal{M}^{(n),\otimes}\big( \chi_{\Omega_{\mmu,\kappa_1}^{f}}\big) \big\Vert_{L^p(\bbr^{2n})}^p\lesssim \big|\Omega_{\mmu,\kappa_1}^f\big|\\
  &\le \big(2^{-\kappa_1}C2^{10n (\mu_1+\mu_2)}\big)^{-p}\big\Vert \MS_{\vec{M}}^{(1)}(f)\big\Vert_{L^p(\bbr^{2n})}^p\\
  &\lesssim 2^{\kappa_1 p}2^{-10n(\mu_1+\mu_2)p}.
\end{align*}
Similarly, applying Chebyshev's inequality, \eqref{fsinq0}, and Lemma \ref{lemmassmm},
\begin{align*}
\big| \Omega_{\RR_{\kappa_1,\kappa_2,\kappa_3}}\big|&\le \bigg|\bigcup_{R\in\RR_{\mmu,\kappa_2}^{g}}R\;\bigg| \lesssim   \big| \Omega_{\mmu,\kappa2}^{g}\big| \\
 &  \le\big(2^{-\kappa_2}C2^{10n (\mu_1+\mu_2)}2^{(M_1+M_2)n(\frac{1}{q_0}-\frac{1}{u})}\big)^{-q}\big\Vert \SM_{\vec{M}}^{(u)}(g)\big\Vert_{L^q(\bbr^{2n})}^q\\
  &\lesssim 2^{\kappa_2 q}2^{-10n(\mu_1+\mu_2)q}.
\end{align*}
Lastly, for any arbitrary $t>2$, we have
\begin{align*}
\big| \Omega_{\RR_{\kappa_1,\kappa_2,\kappa_3}}\big|&\le \bigg|\bigcup_{R\in\RR_{\mmu,\kappa_3}^{h}}R\;\bigg|\lesssim \big| \Omega_{\mmu,\kappa_3}^{h}\big|     \\
 &\le\big(2^{-\kappa_3}C2^{N (\mu_1+\mu_2)}2^{(M_1+M_2)n(\frac{1}{u-\epsilon}-\frac{1}{u})}\big)^{-t}\big\Vert \DS_{\vec{M}}^{(u),\mmu}(h)\big\Vert_{L^t(\bbr^{2n})}^t\\
  &\lesssim_t 2^{\kappa_3 t}
\end{align*} due to \eqref{fsinq0} and Lemma \ref{lmnoncompactbip}.

Averaging the above estimates, we finally obtain that
\begin{equation}\label{omegakeyest}
\big| \Omega_{\RR_{\kappa_1,\kappa_2,\kappa_3}}\big|\lesssim_t \big( 2^{\kappa_1 p}2^{-10n(\mu_1+\mu_2)p}\big)^{\theta_1}\big( 2^{\kappa_2 q}2^{-10n(\mu_1+\mu_2)q}\big)^{\theta_2}2^{t\kappa_3 \theta_3}
\end{equation}
for any $2<t<\infty$ and any $0\le \theta_1,\theta_2,\theta_3\le 1$ with $\theta_1+\theta_2+\theta_3=1$.

Now we separate the left-hand side of \eqref{lastclaim} into
$$\sum_{\kappa_1,\kappa_2,\kappa_3\ge 1} \cdots ~+\sum_{\substack{\kappa_1,\kappa_2 \ge 1\\ -\NN+1\le \kappa_3\le 0}}\cdots~=:\mathbb{I}_1+\mathbb{I}_2. $$

For the estimate of $\mathbb{I}_1$,
we recall that $0<\frac{1}{p}<1$, $0<\frac{1}{q}<1$, and $\frac{1}{p}+\frac{1}{q}=\frac 1r>1$. Thus, we may choose $0<\theta_1<\frac{1}{p}<1$ and $0<\theta_2<\frac{1}{q}<1$ satisfying $\theta_1+\theta_2=1$. Setting $\theta_3=0$ in \eqref{omegakeyest}, the term $\mathbb{I}_1$ can be estimated by
\begin{align*}
2^{-20n(\mu_1+\mu_2)}\Big(\sum_{\kappa_1\ge 1}2^{-\kappa_1(1-p\theta_1)} \Big)\Big(\sum_{\kappa_2\ge 1}2^{-\kappa_2(1-q\theta_2)} \Big)\Big(\sum_{ \kappa_3\ge  1}2^{-\kappa_3}\Big)\lesssim 2^{-20n(\mu_1+\mu_2)}\le 1,
\end{align*} as desired, because $0<p\theta_1<1$ and $0<q\theta_2<1$.

For the remaining term $\mathbb{I}_2$, we choose $0<\theta_1<\frac{1}{p}<1$ and $0<\theta_2<\frac{1}{q}<1$ satisfying $0<\theta_1+\theta_2<1$.
Setting $\theta_3=1-(\theta_1+\theta_2)$ and choosing $t>2$ sufficiently large that $t\theta_3>1$, we have
\begin{align*}
\mathbb{I}_2&\lesssim_t 2^{-20n(\mu_1+\mu_2)(\theta_1+\theta_2)}\Big(\sum_{\kappa_1\ge 1}2^{-\kappa_1(1-p\theta_1)} \Big)\Big(\sum_{\kappa_2\ge 1}2^{-\kappa_2(1-q\theta_2)} \Big)\Big(\sum_{-\NN+1\le \kappa_3\le 0}2^{\kappa_3(t\theta_3-1)}\Big)\\
&\lesssim 2^{-20n(\mu_1+\mu_2)(\theta_1+\theta_2)}\le 1,
\end{align*}
which completes the proof of \eqref{lastclaim}.
\end{proof}


\section{Multi-linear and multi-parameter extensions}
\label{mm}
Before we discuss the general multi-parameter case, we should remark that the argument in previous sections could be modified to the single-parameter setting easily, when the hybrid operators, introduced in Section \ref{prl}, are simpler and the corresponding key lemmas keep valid.

The extension to $d$-parameter setting ($d\geq 3$) follows in a natural way. Recall that when handling the essential difficulty in bi-parameter case, our main ingredients  include for instance, the decomposition in Lemma \ref{compact}, and the estimates for bi-parameter hybrid functions (Lemma \ref{lemmassmm}). Since  Lemma \ref{compact} is applied for each parameter separately, it also gives a straightforward way for using such decompositions in $d$-parameter setting. Then the estimates for the $d$-parameter hybrid functions in Lemma \ref{lmnoncompactbip} can be similarly obtained by using the natural extension of Lemma \ref{t301} and the Fefferman-Stein inequality.

The extension  of the theorem to $l$-linear case ($l>2$) is mainly based on a multi-linear analogue of Proposition \ref{pp42}, which helps us to obtain the boundedness with sharp indices by the multi-linear interpolation. Such an analogue originally appeared in \cite{LHH} in the single-parameter setting when $u=2$. As an extension of Proposition \ref{pp42} to the $l$-linear and $d$-parameter case, Theorem \ref{ml} can actually be reduced to the following via a multi-linear interpolation.
\begin{proposition}
\label{llinearpp}
	Suppose $\frac{1}{l}<p<\infty$, $1<p_1,\dots, p_l<\infty$, and $\frac{1}{p}=\sum_{1\leq i \leq l}\frac{1}{p_i}$. Let $1\leq k \leq l$ be any fixed integer, and for each $i\in J_l\setminus \{k\}$  let $1\leq q_i<\min{(p_i,u)}$. Then for any $\epsilon>0$ there exists a constant $C_\epsilon>0$ such that
	
	\begin{equation}
	\big\|T_{\vec A}^{\vec{M}}(f_1,\dots,f_l)\big\|_{L^p (\mathbb{R}^{dn})} \le  C_\epsilon   \mathcal{A}_{s_1,\dots,s_d}^{u}[m] \bigg( \prod_{j=1}^d 2^{-M_j(s_j-{n\over u-\epsilon}-\sum_{i\neq k}{n\over q_i})}\bigg)\bigg( \prod_{i=1}^{l}\|f_i\|_{L^{p_1}(\mathbb{R}^{dn})}. \label{ld}\bigg)
	\end{equation}	
\end{proposition}
Here $T_{\vec A}^{\vec{M}}$ is the $l$-linear and $d$-parameter version of $T_{A_{12}B_{21}}^{\vec{M}}$ in Proposition \ref{pp42}. To prove this proposition, one just needs to consider three different cases: $(a)$ $p>1$, and $ p_i>2$ for some $1\leq i\leq l$; $(b)$ $p>1$ and $p_i\leq 2$ for some $1\leq i\leq l$; $(c)$ $p<1$. These cases can be handled in a very similar way as the proof in Section \ref{pop}.
We omit a detailed proof of the proposition as it would be simply a repetition in the $l$-linear and $d$-parameter setting.

To see how Proposition \ref{llinearpp} implies Theorem \ref{ml} by interpolation, one can follow the ideas in \cite{LHH}, where the case  $u=2$ is discussed. In fact, this can be accomplished as long as we can reestablish Proposition 4.3 and Lemma 4.4 in \cite{LHH} for general $1< u\leq 2$. Now we give the outline of this procedure.

 For any $B\subset J_l:=\{1,2,\dots,l\}$, 
  $\frac{1}{u}+\frac{1}{u'}=1$,
  and $s=\min\{s_1,\dots,s_d\}>\frac{ln}{u}$,
 we define
\begin{eqnarray*}
  Q(s)&:=&\bigcap_{I\subset J_l} \bigg\{(x_1,\dots,x_l)\in (0,1)^l:\, \Big(\sum_{i\in J_l}x_i\Big)-{1\over u'}<{s\over n}+\sum_{i\in I}\Big(x_i-{1\over u}\Big)\bigg\}, \\
   R_B&:=&\Big\{(x_1,\dots,x_l)\in (0,1)^l:\, 0<x_i\leq{1\over u} \text{ if } i\in B, \ {1\over u}<x_i<1\text{ if } i\in B^c\Big\},\\
  S_B(s)&:=& \bigg\{(x_1,\dots,x_l)\in R_B:\, \Big(\sum_{i\in J_l}x_i\Big)-{1\over u'}<{s\over n}+\sum_{i\in B}\Big(x_i-{1\over u}\Big)\bigg\}.
\end{eqnarray*}

Note that
$$R_B\cap Q(s)=S_B(s)$$
as $$Q(s)= \bigg\{(x_1,\dots,x_1)\in (0,1)^l:\, \Big(\sum_{i\in J_l}x_i\Big)-{1\over u'}<{s\over n}+\sum_{i\in J_l}\min \Big\{0,x_i-{1\over u}\Big\}\bigg\}$$
and
$$\sum_{i\in J_l}\min \Big\{0,x_i-{1\over u}\Big\}=\sum_{i\in B}\Big(x_i-\frac{1}{u}\Big),$$
and thus
\begin{equation}
  \label{qsb}
  Q(s)=\bigcup_{B\subset J_l}\big(R_B\cap Q(s)\big)=\bigcup_{B\subset J_l}S_B(s)
\end{equation} since $\bigcup_{B\subset J_l}R_B=(0,1)^l$.
The followings are the replacements for Proposition 4.3 and Lemma 4.4 in \cite{LHH} for general $1<u\leq 2$.
\begin{proposition}
  \label{pp43u}
  Let $B$ be any subset of $J_l$ and $s>\frac{ln}{u}$. Then for any $({1\over p_1},\dots,{1\over p_l})\in S_B(s)$, we have
$$\big\|{T_m}(f_1,\dots,f_l)\big\|_{L^p(\mathbb{R}^{dn})}\lesssim  \mathcal{A}_{s_1,\dots,s_d}^{u}[m] \prod_{i=1}^l\|f_i\|_{L^{p_i}(\mathbb{R}^{dn})}.$$
\end{proposition}

\begin{proposition}
  \label{pp44u}
  Let $B$ be a proper subset of $J_l$ and $\#B$ denote the number of elements in $B$. Suppose that $s>\frac{ln}{u}$.  For each $l_0 \in B^c$, we define
  $${S_B^{l_0}(s):=\Big\{(x_1,\dots,x_l)\in R_B:\, \frac{s}{n}>{{\#B+1}\over u}+\sum_{i\in B^c\setminus \{l_0\}} x_i\Big\}.}$$
  Then  the convex hull of $\bigcup_{l_0\in B^c}S_B^{l_0}(s)$ equals $S_B(s)$, i.e.,
  $$ \text{convh} \Big(\bigcup_{l_0\in B^c}S_B^{l_0}(s)\Big)=S_B(s).$$
\end{proposition}

Note that assuming Proposition \ref{pp44u} is true, one can follow the ideas in the proof of Proposition 4.3 in \cite{LHH} to see that our Proposition \ref{llinearpp} and Proposition~\ref{pp44u} together would imply Proposition \ref{pp43u} by a multi-linear interpolation. Then according to \eqref{qsb}, Proposition \ref{pp43u} gives the boundedness of {$T_m$} for $(1/p_1,\dots,1/p_l)\in Q(s)$, which is exactly Theorem \ref{ml}.

Therefore, it remains to check Proposition~\ref{pp44u}, for which one can also follow the proof in \cite{LHH}. Here we will only list the key setups.

\begin{proof}[Outline of proof of Proposition \ref{pp44u}]
Without loss of generality one can assume $B^c=J_\kappa=\{1,2,\dots,\kappa\}$ ($1\leq \kappa\leq l$) and then can rewrite the sets $R_B,{S_B(s),S_B^{l_0}(s)}$ as
  \begin{eqnarray*}
   R_B&=& \Big({1\over u},1\Big)^\kappa \times \Big(0,{1\over u}\Big]^{l-\kappa} \\
  S_B(s)&=& \bigg\{(x_1,\dots,x_l)\in R_B:\, \sum_{i\in J_\kappa}x_i<\Big({s\over n}-{l\over u}\Big)+\Big(1+{\kappa-1\over u}\Big)\bigg\} \\
   S_B^{l_0}(s)&=& \bigg\{(x_1,\dots,x_l)\in R_B:\, \sum_{i\in J_\kappa\setminus \{l_0\}}x_i<\Big({s\over n}-{l\over u}\Big)+{\kappa-1\over u}\bigg\}, \q l_0\in J_{\kappa}.
\end{eqnarray*}
Making a change of variable $(x_1,\dots,x_l)\to(y_1,\dots,y_l)$ with
$$y_i=\begin{cases}
x_i, &  i\in B=\{\kappa+1,\dots,l\}\\
x_i-1/u, &  i\in B^c=\{1,\dots,\kappa\}
\end{cases},$$
the above three sets would be translated to
\begin{eqnarray*}
   \mathcal{R}_B&=& \Big(0,{1\over u'}\Big)^\kappa \times \Big(0,{1\over u}\Big]^{l-\kappa}=\mathcal{R}\times \Big(0,{1\over u}\Big]^{l-\kappa}, \\
  \mathcal{S}_B(s)&=& \Big\{(y_1,\dots,y_l)\in \mathcal{R}_B:\, \sum_{i\in J_\kappa}y_i< \beta+{1\over u'}\Big\}=\mathcal{S}\times \Big(0,{1\over u}\Big]^{l-\kappa}, \\
   \mathcal{S}_B^{l_0}(s)&=& \Big\{(y_1,\dots,y_l)\in \mathcal{R}_B:\, \sum_{i\in J_\kappa\setminus \{l_0\}}y_i< \beta \Big\}=\mathcal{S}^{l_0}\times \Big(0,{1\over u}\Big]^{l-\kappa},
\end{eqnarray*}
where we set $\beta={s\over n}-{l\over u}$ and
\begin{eqnarray*}
   \mathcal{R}&:=& (0,{1\over u'})^\kappa, \\
  \mathcal{S}&:=& \Big\{(y_1,\dots,y_\kappa)\in \mathcal{R}:\, \sum_{i\in J_\kappa}y_i< \beta+{1\over u'}\Big\}, \\
   \mathcal{S}^{l_0}&:=& \Big\{(y_1,\dots,y_\kappa)\in \mathcal{R}:\, \sum_{i\in J_\kappa\setminus \{l_0\}}y_i< \beta \Big\}.
\end{eqnarray*}
Then the conclusion in Proposition \ref{pp44u} will be equivalent to
$$\text{convh}\Big(\bigcup_{l_0\in J_{\kappa}}\mathcal{S}^{l_0}\Big)=\mathcal{S}.$$
{This clearly follows from  the ideas in the proof of Lemma 4.4 in \cite{LHH}, replacing ${1\over 2}$ with ${1\over u'}$.}
\end{proof}

\section{Counterexamples}
\label{ctexample}

In this section, we construct counterexamples indicating the sharpness of Theorem~\ref{ml}.
We start from the single-parameter case, as we will see below that the multi-parameter case follows immediately.

\subsection{Single-parameter examples}

Let us discuss the sharpness of Theorem~\ref{ml} in single-parameter setting first.
{For a function $m$ on $\bbr^{ln}$, $1\le u<\infty$, and $s\ge 0$, let
$$\mathcal{A}_{s}^{u}[m]:=\sup_{k\in\bbz}\big\Vert m(2^{k}\cdot)\wh{\Theta}(\cdot)\big\Vert_{L^u_{s}(\bbr^{ln})}.$$}
Our first result is as follows.

\begin{theorem}\label{t326}

Let $1<p_1,\dots, p_l<\infty$, $\tf1p=\tf1{p_1}+\cdots+\tf1{p_l}$, $1<u<\infty$, and $0<s<\infty$.
Suppose that {
\begin{equation}\label{e320}
\big\|{T_m}(f_1,\dots,f_l)\big\|_{L^p(\bbr^{n})}\lesssim  \mathcal{A}_{s}^{u}[m] \prod_{i=1}^l\|f_i\|_{L^{p_i}(\bbrn)}
\end{equation}
for all  $m$ satisfying $\mathcal{A}_{s}^{u}[m]<\infty$.
Then it is necessary to have
\begin{equation}\label{e321}
\tf sn\ge \max\Big\{\tf{l}u,\max_I\Big( \tf1p-\tf1{2}-\sum_{i\in I}\big(\tf1{p_i}-\tf12\big)\Big),\max_I\Big(\tf1p-\tf1{u'}-\sum_{i\in I}\big(\tf1{p_i}-\tf1u \big)\Big)\Big\},
\end{equation}
where the maximums in the bracket run over all subsets $I\sset\{1,2,\dots, l\}$.}

\end{theorem}

\begin{remark}\label{t328}
(i) The necessary condition
\begin{equation}\label{e334}
\tf sn\ge \max\Big\{\tf{l}u,{\max_I} \Big( \tf1p-\tf1{2}-\sum_{i\in I}\big(\tf1{p_i}-\tf12\big) \Big)\Big\}
\end{equation}
is obtained in
{\cite[Theorems 2 and 3]{Grafakos2016a}.}

{(ii) In order to prove the necessity of the condition
\begin{equation}\label{songe}
\frac{s}{n}\ge \max_I\bigg(\f1p-\f1{u'}-\sum_{i\in I}\Big(\f1{p_i}-\f1u \Big)\bigg),
\end{equation}
it suffices to consider only the case $1<u\le 2$ as the other cases can be done by \eqref{e334}.

(iii) When  $\#I\ge l-1$, it is straightforward that
$$\frac{l}{u}\ge \frac{1}{p}-\f1{u'}-\sum_{i\in I}\Big(\f1{p_i}-\f1u\Big)$$ for all $p_1,\dots,p_l>1$.}

\end{remark}

\begin{proof}[Proof of Theorem \ref{t326}]

By the above remark, it suffices to consider $1<u\le 2$  and prove \eqref{songe} when $0\le \#I\le l-2$. Here we take the convention $I=\emptyset$ when $\#I=0$. We discuss only the case $n=1$ as the higher dimensional extension could be obtained easily using tensor products; see \cite[Proposition 7]{Grafakos2016a} for details.

When the set $I$ is not empty, without loss of generality, we may assume $I=\{1,2,\dots,k\}$ for $1\le k\le l-2$ as other cases follow via symmetry.
 For a large number $N$, let
$$
E_k^N:=\begin{cases}\big\{\vec{j}:=(j_1,\dots,j_l)\in\mathbb N^l: j_{1}+\cdots+j_l=N \big\},& k=0\\
\big\{\vec{j}:=(j_1,\dots,j_l)\in\mathbb N^l: j_1=\cdots=j_k=\tf Nl, j_{k+1}+\cdots+j_l=\tf{l-k}lN \big\},& 1\le k\le l-2
\end{cases}
$$
and we define a multiplier
$$
m(\vec{\xi}\,):=\sum_{\vec{j}\in E_k^N}\prod_{i=1}^l\phi(N\xi_i-j_i), \qq \vec{\xi}:=(\xi_1,\dots,\xi_l)\in\bbr^l
$$
where $\phi\in\mathscr S(\mathbb R)$ is supported in $[-\tf1{10},\tf1{10}]$ and equal to  $1$ on $[-\tf1{20},\tf1{20}]$.
Applying complex interpolation between $s\in\mathbb N$, one can verify that
\begin{equation}\label{e322}
\mathcal{A}_{s}^{u}[m]\lesssim\|m\|_{L^u_s(\bbr^l)}\lesssim N^sN^{-\tf{k+1}u}
\end{equation}
as $m$ is supported in the annulus of a constant size and $\#E_k^N\approx N^{l-k-1}$  for $k\le l-2$.

Next we fix a function $\vp\in\mathscr S(\mathbb R)$ supported in $[-\tf1{100},\tf1{100}]$, and define
$${\wh{f_i}(\xi_i):=\begin{cases}
\sum_{j_i=1}^N \vp(N\xi_i-j_i),& i=k+1,\dots,l\\
\vp(N\xi_i-\tf Nl), & otherwise
\end{cases}.}$$
Then it follows from the calculation in \cite[Proposition 7]{Grafakos2016a} that
{\begin{equation}\label{e323}
\|f_i\|_{L^{p_i}(\bbr)}\approx
\begin{cases}1, & i=k+1,\dots,l\\
N^{\frac{1}{p_i}-1}, & otherwise
\end{cases}.
\end{equation}}

One can also verify that
{
\begin{equation}\label{e325}
\big\|T_m(f_1,\dots, f_l)\big\|_{L^p(\bbr)}\approx\bigg(\int_{\bbr}\Big|\sum_{j\in E_k^N}\tf1{N^l}\vp^\vee(\tf tN)^l \Big|^p \;dt\bigg)^{\frac{1}{p}}\approx N^{\tf1p-k-1}.
\end{equation}}

Applying the estimates \eqref{e322}-\eqref{e325} to \eqref{e320}, we obtain
$$
N^{\tf1p-k-1}\lesssim N^{s-\tf{k+1}u}\prod_{i=1}^kN^{\tf1{p_i}-1}{=N^{s-\tf {k+1}{u}+\tf 1{p_1}+\cdots+\tf{1}{p_k}-k}.}
$$
Letting $N\to\infty$, we obtain the necessary condition
{$$s\ge \frac{1}{p}-\Big(\frac{1}{p_1}+\cdots+\frac{1}{p_k} \Big)+\frac{k+1}{u}-1,$$
which is equivalent to}
\eqref{songe} {when $I=\{1,\dots,k\}$.}
\end{proof}


Even  in the bilinear setting, the above examples do not give the full range of unboundedness up to endpoints.
{Let us discuss this in detail below.
The condition
$$
\tf sn\ge  \max_I\tf1p-\tf1{2}-\sum_{i\in I}(\tf1{p_i}-\tf12)
$$
can be written as
\begin{equation}\label{e401}
\tf sn\ge \tf1p-\tf1{2}+\sum_{i=1}^2\max\{\tf12-\tf1{p_i},0\}.
\end{equation}
In particular,
\begin{equation}\label{e400}
\tf sn\ge \tf1{p_1}.
\end{equation}
As observed in \cite{GraNg}, when $p>2$, we have $p_1,p_2>2$, therefore by duality any bilinear H\"ormander multiplier bounded from $L^{p_1}(\mathbb R^n)\times L^{p_2}(\mathbb R^n)$ to $L^{p}(\mathbb R^n)$ is bounded from
$L^{p'}(\mathbb R^n)\times L^{p_2}(\mathbb R^n)$ to $L^{p_1'}(\mathbb R^n)$.
Applying \eqref{e400}, we obtain
$$
\tf sn\ge\tf1{p'} =\sum_{i=1}^2\max\{\tf12-\tf1{p_i},0\}>\tf1p-\tf1{2}+\sum_{i=1}^2\max\{\tf12-\tf1{p_i},0\}.
$$
Therefore we come to a  necessary condition better than \eqref{e401}, i.e.
\begin{equation}\label{e335}
\tf sn\ge \max\{\tf12-\tf1{p'},0\}+\sum_{i=1}^2\max\{\tf12-\tf1{p_i},0\}. 
\end{equation}
}

In  last theorem, we cannot rule out the endpoints, which means that we do not know what happens when
$$
\tf sn= \max\Big\{\tf{l}u,{\max_I}\Big(\tf1p-\tf1{2}-\sum_{i\in I}(\tf1{p_i}-\tf12)\Big), {\max_I}\Big(\tf1p-\tf1{u'}-\sum_{i\in I}(\tf1{p_i}-\tf1u) \Big)\Big\}.
$$
Some partial results however can be obtained using   variants of Bessel kernels introduced by the fourth author in \cite{Park}. The exact statement is as follows.
{
\begin{theorem}\label{t327}
Let  $0<p_1,\dots, p_l<\infty$, $\tf1p=\tf1{p_1}+\cdots+\tf1{p_l}$, $1<u<\infty$, and  $0<s<\infty$.
Suppose that 
for all $m$ satisfying $\mathcal{A}_s^{u}[m]<\infty$, $T_m$ satisfies \eqref{e320}.
Then it is necessary to have
$$\frac{s}{n}>\max\Big\{ \frac{l}{u},\frac{1}{p}-\frac{1}{u'}\Big\}.$$
\end{theorem}
}
This theorem follows from the following two propositions.

\begin{proposition}\label{t320}

Let $1<u<\infty$ and $0<p_1,\dots, p_l<\infty$ with $\tf1p=\tf1{p_1}+\cdots+\tf1{p_l}$. Suppose that
\begin{equation}
s\le\tf{ln}u.
\end{equation}
Then 
{there exists a multiplier $m$ on $\bbr^{ln}$ such that
$\mathcal A_s^u[m]<\infty$, but
\begin{equation}\label{e327}
\|T_m\|_{H^{p_1}\times \cdots \times H^{p_l}\to L^p}=\infty  .
\end{equation}}
\end{proposition}

\begin{proposition}\label{t321}

Let $1<u<\infty$ and $0<p_1,\dots, p_l<\infty$ with $\tf1p=\tf1{p_1}+\cdots+\tf1{p_l}$.
Suppose that
\begin{equation}\label{e338}
\tf{ln}u<s\le \tf np-\tf n{u'}.
\end{equation}
Then  there exists a  multiplier $m$  on $\bbr^{ln}$ such that
{$\mathcal A_s^u[m]<\infty$, but
\begin{equation}
\|T\|_{H^{p_1}\times \cdots \times H^{p_l}\to L^p}=\infty.
\end{equation}}
\end{proposition}

To prove these two results, we need to prepare several lemmas.
We first introduce a variant of the Bessel kernel, which is the key ingredient in constructing counterexamples for the sharpness.
For any $N\in\bbn$, and $t,\gamma>0$, let
$$
\mathcal{H}_{(t,\gamma)}^{(N)}(x):=\frac{1}{(1+4\pi^2 |x|^2)^{t/2}} \frac{1}{(1+\ln(1+4\pi^2|x|^2))^{\gamma/2}}, \qq x\in\mathbb{R}^N.
$$
Concerning the $L^u$-norm of $\mathcal{H}_{(t,\gamma)}^{(N)}$, we have the following lemma from \cite[Section 4]{Gr_Park}.
\begin{lemma}
Suppose that $N\in\bbn$, and $t,\gamma>0$. Then
\begin{equation}\label{hntcondition1'}
\Vert \mathcal{H}_{(t,\gamma)}^{(N)}\Vert_{L^u(\mathbb{R}^N)}<\infty\qq \Leftrightarrow \qq t>\frac{N}{u}\q \text{or}\q t=\frac{N}{u}, \gamma>\frac{2}{u},
\end{equation}
\begin{equation}\label{hntcondition2'}
\Vert \wh{\mathcal{H}_{(t,\gamma)}^{(N)}}\Vert_{L^u(\mathbb{R}^N)}<\infty\qq \Leftrightarrow \qq t>\frac{N}{u'}\q \text{or}\q t=\frac{N}{u'}, \gamma>\frac{2}{u}.
\end{equation}

\end{lemma}

We need also the following estimate on Sobolev norms, which can be found in \cite[Lemma~7.5.8]{G250}.
\begin{lemma}\label{t322}
Let $s\ge 0$, $1<u<\infty$, and $\vp\in\mathscr S(\mathbb R^N)$. Then there exists a constant $C_{\vp,s}$ such that
$$
\big\|f\cdot \vp\big\|_{L^u_s(\bbr^N)}\le C_{\vp,s}\|f\|_{L^u_s(\bbr^N)}
$$
for any $f\in L^u_s(\mathbb R^N)$.

\end{lemma}

Next we introduce some auxiliary functions.
Let $\theta$ and $\wt{\theta}$ be Schwartz functions on $\mathbb{R}^n$ satisfying
$$\theta\ge 0, \q \theta(0)\not= 0, \q  \supp(\wh{\theta})\subset \Big\{\xi\in\mathbb{R}^n: |\xi|\le \frac{1}{200l}\Big\},$$
$$\wh{\wt{\theta}}(\xi)=1~\text{for}~|\xi|\le \frac{1}{100l}, \q \text{and}\q \supp(\wh{\wt{\theta}})\subset \Big\{\xi\in\mathbb{R}^n: |\xi|\le \frac{1}{10l}\Big\}.$$
Then we observe that for any $t,\gamma>0$
\begin{equation}\label{htgthetaest}
\mathcal{H}_{(t,\gamma)}^{(n)}\ast \theta(x)\gtrsim_{\theta}\mathcal{H}_{(t,\gamma)}^{(n)}(x), \q x\in\bbrn.
\end{equation}
{Choose a nonnegative Schwartz function $\wh\Gamma$ on $\mathbb R^{ln}$ such that $\wh\Gamma(\vec y)=1$ on the unit ball $\big\{\vec{y}=(y_1,\dots,y_l)\in \mathbb R^{ln}: |\vec{y}|\le 1 \big\}$ and $\wh\Gamma$ is supported in a larger ball.}
For sufficiently large $M$, let $\Gamma_M(\vec{\xi}\,):=M^{ln}\Gamma(M\vec{\xi}\,)$.
Let $1<t, \gamma<\infty$ to be chosen later and set
{$$\mathcal{H}_{(t,\gamma)}^{(ln,M)}(\vec{y}\,):=\mathcal{H}_{(t,\gamma)}^{(ln)}(\vec{y}\,)\wh{\Gamma_M}(\vec{y}\,),\qq \vec{y}\in \bbr^{ln}.$$}

Now we are ready to prove Proposition~\ref{t320}.

\begin{proof}[Proof of Proposition~\ref{t320}]
Write $\nu:=(l^{-{1/2}},0,\dots,0)\in\mathbb{R}^n$
and  define
 $$
 m_{(t,\gamma)}^{(M)}(\vec{\xi}\,):=\wh{\mathcal{H}_{(t,\gamma)}^{(ln,M)}}(\xi_1-\nu,\dots,\xi_l-\nu)\prod_{i=1}^l\wh{\wt{\theta}}(\xi_i-\nu),
 $$
 where we recall the notation $\vec{\xi}=(\xi_1,\dots,\xi_l)\in(\mathbb R^{n})^l$.
Then it follows from the compact support condition of $\wh{\wt{\theta}}$ that $m_{(t,\gamma)}^{(M)}$ is supported in the annulus
$$ \Big\{\vec{\xi}\in\bbr^{ln}:   ~ \frac{4}{5}\le \big| \vec{\xi}\big|\le \frac{6}{5}  \Big\},$$
which implies that $$m_{(t,\gamma)}^{(M)}(2^{k}\vec{\xi}\,)\wh{\Theta}(\vec{\xi}\,)$$
 vanishes unless $-1\le k\le 1$.
Therefore, we have
\begin{equation}\label{aammtg}
\mathcal{A}_{s}^u\big[m_{(t,\gamma)}^{(M)}\big]\lesssim \max_{-1\le k \le 1}\big\Vert m_{(t,\gamma)}^{(M)}(2^{k }\vec{\cdot}\, )\wh{\Theta}  \big\Vert_{L^u_{s}(\bbr^{ln})}
\end{equation}
and, using Lemma \ref{t322}, this can be dominated by
\begin{equation}\label{mmmtgest}
{ \max_{-1\le k \le 1}\big\Vert m_{(t,\gamma)}^{(M)}(2^{k}\vec{\cdot}\,)\big\Vert_{L^u_{s}(\bbr^{ln})}} \lesssim \Vert m_{(t,\gamma)}^{(M)} \Vert_{L^u_{s}(\bbr^{ln})}.
\end{equation}
Moreover, 
\begin{equation}\label{mmtgest}
\big\Vert m_{(t,\gamma)}^{(M)} \big\Vert_{L^u_{s}(\bbr^{ln})}= \Big\Vert \wh{\mathcal{H}^{(ln,M)}_{(t,\gamma)}}\; \cdot\; \otimes_{i=1}^l\wh{\wt{\theta}}\,\Big\Vert_{L_{s}^u(\bbr^{ln})}
\lesssim \Big\Vert \wh{\mathcal{H}_{(t,\gamma)}^{(ln,M)}}\Big\Vert_{L^u_{s}(\bbr^{ln})}
\end{equation}
where Lemma \ref{t322} is applied in the last inequality as $\otimes_{i=1}^l\wh{\wt{\theta}}$ is a Schwartz function in $\bbr^{ln}$.
We observe that
$$(I-\Delta_{\vec{y}})^{s/2}\wh{\mathcal{H}_{(t,\gamma)}^{(ln,M)}} =\wh{\mathcal{H}_{(t-s,\gamma)}^{(ln,M)}}= \wh{\mathcal{H}_{(t-s,\gamma)}^{(ln)}}\ast \Gamma_M $$
and thus Young's convolution inequality yields that
$$\Big\Vert \wh{\mathcal{H}_{(t,\gamma)}^{(ln,M)}}\Big\Vert_{L^u_{s}(\bbr^{ln})}  \lesssim \Big\Vert \wh{\mathcal{H}_{(t-s,\gamma)}^{(ln)}}\Big\Vert_{L^u(\bbr^{ln})}\big\Vert \Gamma_M\big\Vert_{L^1(\bbr^{ln})} \approx \Big\Vert \wh{\mathcal{H}_{(t-s,\gamma)}^{(ln)}}\Big\Vert_{L^u(\bbr^{ln})} $$
uniformly on $M$.
This shows
\begin{equation}\label{amestfinal'}
\mathcal{A}_{s}^{u}\big[m_{(t,\gamma)}^{(M)}\big]\lesssim \Big\Vert \wh{\mathcal{H}_{(t-s,\gamma)}^{(ln)}}\Big\Vert_{L^u(\bbr^{ln})} \qq \text{ uniformly on }~M.
\end{equation}

\hfill

On the other hand, in order to estimate the operator norm of $T_{m_{(t,\gamma)}^{(M)}}$, let $\epsilon>0$ be sufficiently small and for each $i\in \{1,\dots,l\}$,
$$
f_i^{(\epsilon)}(y_i):=\epsilon^{\tf n {p_i}}\theta(\epsilon y_i)e^{2\pi i\langle y_i,\nu\rangle}
$$
for $y_i\in\bbrn$.
Using the square function characterization of Hardy spaces, we obtain
\begin{equation}\label{filpiest}
\big\Vert f_i^{(\epsilon)}\big\Vert_{H^{p_i}(\bbrn)}\approx\big\Vert f_i^{(\epsilon)}\big\Vert_{L^{p_i}(\bbrn)}=\Vert \theta\Vert_{L^{p_i}(\bbrn)} \,\lesssim 1\qq \text{uniformly on }~\epsilon.
\end{equation}
Since
$$
\wh{f_i^{(\epsilon)}}(\xi_i)=\epsilon^{\tf n {p_i}}\frac{1}{\epsilon^{n}}\wh{\theta}\Big(\frac{\xi_i-\nu}{\epsilon}\Big)
$$
we have
\begin{align*}
&\Big| T_{m_{(t,\gamma)}^{(M)}}\big(f_1^{(\epsilon)},\dots,f_l^{(\epsilon)}\big)(x)\Big|\\
&=\epsilon^{\tf n p}\bigg|\int_{\bbr^{ln}}     \wh{\mathcal{H}_{(t,\gamma)}^{(ln,M)}}(\xi_1-\nu,\dots,\xi_l-\nu)\Big[ \prod_{i=1}^l\frac{1}{\epsilon^{n}}\wh{\theta}\Big(\frac{\xi_i-\nu}{\epsilon} \Big)\Big]
e^{2\pi i\langle x,\xi_1+\cdots+\xi_l\rangle}       ~d\vec{\xi}\,\bigg|\\
&=\epsilon^{n/p}\big| \mathcal{H}_{(t,\gamma)}^{(ln,M)}\ast \big(\theta(\epsilon \cdot)\otimes\cdots\otimes\theta(\epsilon\cdot) \big)(x,\cdots,x) \big|. 
\end{align*}
This yields that
\begin{align*}
&\Big\Vert T_{m_{(t,\gamma)}^{(M)}}\big(f_1^{(\epsilon)},\cdots, f_l^{(\epsilon)}\big)\Big\Vert_{L^p{(\bbr^{n})}}=\Big\Vert  \mathcal{H}_{(t,\gamma)}^{(ln,M)}\ast \big(\theta(\epsilon \cdot)\otimes\cdots\otimes\theta(\epsilon\cdot) \big)(x/\epsilon,\cdots,x/\epsilon) \Big\Vert_{L^p(x)}
\\
&\approx \bigg( \int_{\bbrn}  \Big| \int_{\bbr^{ln}}  \mathcal{H}_{(t,\gamma)}^{(ln,M)}(\vec{y}\,)\theta(x-\epsilon y_1)\cdots\theta(x-\epsilon y_l)    ~d\vec{y}\,\Big|^p   ~dx\bigg)^{1/p}.
\end{align*}
Applying Fatou's lemma, we have
\begin{align}
&\big\Vert T_{m_{(t,\gamma)}^{(M)}} \big\Vert_{H^{p_1}(\bbr^{n})\times\cdots\times H^{p_l}(\bbr^{n})\to L^p(\bbr^{n})}\gtrsim \liminf_{\epsilon\to 0}\Big\Vert T_{m_{(t,\gamma)}^{(M)}}\big(f_1^{(\epsilon)},\dots,f_l^{(\epsilon)}\big)\Big\Vert_{L^p(\bbr^{n})}\nonumber\\
&\q \gtrsim {\bigg( \int_{\bbrn} \liminf_{\epsilon\to 0} \Big| \int_{\bbr^{ln}}  \mathcal{H}_{(t,\gamma)}^{(ln,M)}(\vec{y}\,)\theta(x-\epsilon y_1)\cdots\theta(x-\epsilon y_l)    ~d\vec{y}\, \Big|^p   ~dx\bigg)^{1/p}. }\label{liminfep}
\end{align}
Since
$$ \Big| \mathcal{H}_{(t,\gamma)}^{(ln,M)}(\vec{y}\,)\theta(x-\epsilon y_1)\cdots\theta(x-\epsilon y_l) \Big|\lesssim \mathcal{H}_{(t,\gamma)}^{(ln,M)}{(\vec{y})}\q \text{ uniformly in }~\epsilon>0,~x\in\bbrn$$
and
$$\big\Vert \mathcal{H}_{(t,\gamma)}^{(ln,M)}\big\Vert_{L^1{(\bbr^{ln})}}\le \big\Vert  \wh{\Gamma_M}\big\Vert_{L^1({\bbr^{ln})}}\approx  M^{ln}<\infty,$$
the Lebesgue dominated convergence theorem proves that the right-hand side of \eqref{liminfep} is equal to
\begin{equation*}
\big\Vert |\theta|^l\big\Vert_{L^p(\bbrn)}\big\Vert \mathcal{H}_{(t,\gamma)}^{(ln,M)}\big\Vert_{L^1({\bbr^{ln})}} \approx \big\Vert \mathcal{H}_{(t,\gamma)}^{(ln,M)}\big\Vert_{L^1{(\bbr^{ln})}}.
\end{equation*}
Finally, taking $\liminf_{M\to \infty}$ {and applying Fatou's lemma}, we have
\begin{align}
&\liminf_{M\to\infty}\big\Vert T_{m_{(t,\gamma)}^{(M)}} \big\Vert_{H^{p_1}(\bbr^{n})\times\cdots\times H^{p_l}(\bbr^{n})\to L^p(\bbr^{n})}\nonumber\\
&\gtrsim \liminf_{M\to \infty}\big\Vert \mathcal{H}_{t,\gamma}^{(ln,M)}\big\Vert_{L^1(\bbr^{ln})}\ge\Big\Vert \mathcal{H}_{(t,\gamma)}^{(ln)} \cdot \liminf_{M\to \infty}\wh{\Gamma_M}(\vec{\cdot}\,)\Big\Vert_{L^1(\bbr^{ln})} =\big\Vert   \mathcal{H}_{(t,\gamma)}^{(ln)}  \big\Vert_{L^1({\bbr^{ln}})}.\label{limlowerbd'}
\end{align}

Now choosing $$t=ln ~\text{ and }~ \gamma=2 $$
and applying \eqref{hntcondition2'} and \eqref{hntcondition1'} to \eqref{amestfinal'} and  \eqref{limlowerbd'}, respectively, we conclude that
$$\mathcal{A}_{s}^u\big[m_{(ln,2)}^{(ln,M)}\big]\lesssim 1 \q \text{ uniformly on }~M,$$
but
$$\liminf_{M\to\infty}\big\Vert T_{m_{(ln,2)}^{(ln,M)}} \big\Vert_{H^{p_1}(\bbr^{n})\times\cdots\times H^{p_l}(\bbr^{n})\to L^p(\bbr^{n})}=\infty,$$
which verify \eqref{e327} when $s\le ln/u$.
\end{proof}

Let us prove Proposition~\ref{t321} below.

\begin{proof}[Proof of Proposition~\ref{t321}]
In this case, we need a more complicated multiplier, whose construction is  inspired by \cite{GP}.

Let
$$
{m_{(t,\gamma)}(\vec{\xi}\,):=\wh{\mathcal{H}_{(t,\gamma)}^{(n)}}}(\tf{\xi_1+\cdots+\xi_l}l-\nu)\wh\theta(\tf{\xi_1+\cdots+\xi_l}l-\nu)
\wh\theta(\tf{\xi_1+\cdots+\xi_l-l\xi_2}l)\cdots \wh\theta(\tf{\xi_1+\cdots+\xi_l-l\xi_l}l).
$$
We define the non-degenerate matrix  $\mathcal{R}$ by
\begin{equation}\label{e328}
\mathcal{R}(\vec{\xi}\,):=(\tf{\xi_1+\cdots+\xi_l}l-\nu,\tf{\xi_1+\cdots+\xi_l-l\xi_2}l,\cdots,\tf{\xi_1+\cdots+\xi_l-l\xi_l}l).
\end{equation}
so that
$${m_{(t,\gamma)}(\vec{\xi}\,)}=m_1\big(\mathcal{R}(\vec{\xi}\,)\big)$$
where
$$
m_1(\vec{\xi}\,):={\wh{\mathcal{H}_{(t,\gamma)}^{(n)}}}(\xi_1-\nu)\wh\theta(\xi_1-\nu)
\wh\theta(\xi_2)\cdots \wh\theta(\xi_l).
$$
 This implies that
${\wh{m_{(t,\gamma)}}(\vec{y})}=\wh {m_1}\big((\mathcal{R}^{-1})^t(\vec{y}\,)\big)$
where $(\mathcal{R}^{-1})^t$ is the permutation matrix of $\mathcal{R}^{-1}$.
As a result, we have
\begin{align*}
{\big\|m_{(t,\gamma)}\big\|_{L^u_s(\bbr^{ln})}}&=\Big\|\Big(\big(1+4\pi^2|\vec{\cdot}\,|^2\big)^{s/2}\wh {m_1}\big((\mathcal{R}^{-1})^t(\vec{\cdot}\,)\big)\Big)^\vee\Big\|_{L^u(\bbr^{ln})}\\
&\approx\Big\|\Big(\big(1+4\pi^2\big|\mathcal{Q}(\vec{\cdot}\,)\big|^2\big)^{s/2}\wh {m_1}\Big)^\vee\Big\|_{L^u(\bbr^{ln})}
\end{align*}
where $\mathcal{Q}:=\big((\mathcal{R}^{-1})^t\big)^{-1}$.

As {$m_{(t,\gamma)}$} is supported in the annulus $\big\{\vec{\xi}\in\bbr^{ln}:   ~ \frac{4}{5}\le \big| \vec{\xi}\big|\le \frac{6}{5}  \big\}$,
we obtain
\begin{equation*}
\mathcal{A}_{s}^u\big[m_{t,\gamma} \big]\lesssim\Big\|\Big(\big(1+4\pi^2\big|\mathcal{Q}(\vec{\cdot}\,)\big|^2\big)^{s/2}\wh {m_1}\Big)^\vee\Big\|_{L^u(\bbr^{ln})}.
\end{equation*}

To complete the proof, the following lemma is in order.
\begin{lemma}\label{t323}
Let $1<u<\infty$ and $s>0$. Let $\mathcal{Q}$  be a non-degenerate matrix. Then
$$
a(\vec{y}\,):=\f{\big(1+4\pi^2|\mathcal{Q}(\vec{y}\,)|^2\big)^{s/2}}{(1+4\pi^2|y_1|^2)^{s/2}\cdots(1+4\pi^2|y_l|^2)^{s/2}},\q \vec{y}=(y_1,\dots,y_l)
$$
is an $L^u$-multiplier.

\end{lemma}

We leave the proof of this lemma to the end.

Applying Lemma~\ref{t323}, we control $\mathcal{A}_{s}^u[m_{(t,\gamma)} ]$ by
\begin{equation}
\Big\|\Big(\big(1+4\pi^2|\cdot_1|^2\big)^{s/2}\cdots(1+4\pi^2|\cdot_l|^2)^{s/2}\wh {m_1}\Big)^\vee\Big\|_{L^u(\bbr^{ln})}\lesssim {\big\| \wh{\mathcal{H}_{(t,\gamma)}^{(n)}}\big\|_{L^u_s(\bbrn)}\lesssim \big\|\wh{\mathcal{H}_{(t-s,\gamma)}^{(n)}}\big\|_{L^u(\bbrn)}}.\label{e329}
\end{equation}

Next we define
$$ f_i(y_i):=\wt{\theta}(y_i)e^{2\pi i\langle y_i,\nu\rangle}\q \text{ for $1\le i\le l$.}$$
{ Similar to \eqref{filpiest}, it is clear that
$$
\big\Vert f_i\big\Vert_{H^{p_i}(\bbrn)}\approx\big\Vert f_i\big\Vert_{L^{p_i}(\bbrn)} =\big\Vert \wt{\theta}\big\Vert_{L^{p_i}(\bbrn)}\lesssim 1.
$$}
As $\otimes_{i=1}^l\wh f_i=1$ on the support of {$m_{(t,\gamma)}$},
one can easily verify that
{\begin{align*}
\big|T_{m_{(t,\gamma)}}\big(f_1,\dots, f_l\big)(x)\big|&=\bigg|\int_{\bbr^{ln}} m_1\big(\mathcal{R}(\vec{\xi}\,)\big) e^{2\pi i\langle x,\xi_1+\cdots+\xi_l\rangle}\, d\vec{\xi}\, \bigg|\\
&\approx \bigg|\int_{\bbr^{ln}} m_1(\vec{\xi}\,)e^{2\pi i \langle lx, \xi_1\rangle}\,d\vec{\xi}\, \bigg|\\
&\approx \big(\mathcal{H}_{(t,\gamma)}^{(n)}\ast\theta\big)(lx)\big|\theta(0)\big|^{l-1}
\end{align*}}
where  we use \eqref{e328} in the second line. Consequently, by \eqref{htgthetaest}, we have
\begin{equation}\label{e330}
\big\|T_{m_{(t,\gamma)}}(f_1,\dots, f_l)\big\|_{L^p(\bbrn)}\gtrsim \big\|{\mathcal{H}_{(t,\gamma)}^{(n)}}\big\|_{L^p(\bbrn)}.
\end{equation}

By assumption \eqref{e338}, $\tf1p>\tf{l-1}u\ge\tf1u$. Taking $t=\tf np$, $\gamma=\tf2p>\tf2u$, we have
$$
\mathcal{A}_{s}^u[m_{(t,\gamma)} ]<\infty
$$
by \eqref{e329}, and
$$
\big\Vert T_{m_{(t,\gamma)}} \big\Vert_{H^{p_1}(\bbr^{n})\times\cdots\times H^{p_l}(\bbr^{n})\to L^p(\bbr^{n})}=\infty
$$
by \eqref{e330}.

We now turn to the proof of Lemma~\ref{t323}.

\begin{proof}[Proof of Lemma~\ref{t323}]

As $\mathcal{Q}$ is a fixed matrix, we have
$$
\big|\mathcal{Q}(\vec{y})\big|\lesssim_{\mathcal{Q}} |y_1|+|y_2|+\cdots+|y_l|,
$$
hence $a\in L^\infty$.

To show that $a$ is an $L^u$-multiplier, it suffices to verify that it satisfies
\begin{equation}\label{e331}
|\partial^\alpha a(\vec{y})|\lesssim |y_1|^{-|\al_1|}\cdots |y_l|^{-|\al_l|}\qq \al=(\al_1,\dots, \al_l)\in\bbn^{ln},
\end{equation}
since such a multiplier is a Marcinkiewicz multiplier which is bounded on $L^u$.

Obviously
$$
\partial^\al a(y)=\sum_{\al'+\al''=\al}\partial^{\al'}_y[(1+4\pi^2|\mathcal{Q}(\vec{y})|^2)^{s/2}]\prod_{i=1}^l\partial_{y_i}^{\al''_i}[(1+4\pi^2|y_i|^2)^{-s/2}].
$$
We claim that
\begin{equation}\label{e332}
|\partial^{\al'}_y[(1+4\pi^2|\mathcal{Q}(\vec{y})|^2)^{s/2}]|\lesssim_Q(1+|\vec{y}|^2)^{\tf{s-|\al'|}2}
\end{equation}
and
\begin{equation}\label{e333}
|\partial_{y_i}^{\al''_i}[(1+4\pi^2|y_i|^2)^{-s/2}]|\lesssim_Q(1+|\vec{y}|^2)^{-\tf{s+|\al''_i|}2}.
\end{equation}
It follows from these two inequalities that
\begin{align*}
\big|\partial^\al a(\vec{y})\big|\lesssim& \left(\f{1+|\vec{y}|^2}{\prod_{i=1}^l(1+|y_i|^2))}\right)^{s/2}(1+|\vec{y}|^2)^{-|\al'/2|}\prod_{i=1}^l(1+|y_i|^2)^{-|\al_i''|/2}\\
\lesssim & |y_1|^{-|\al_1|}\cdots |y_l|^{-|\al_l|},
\end{align*}
which is \eqref{e331}.

It remains to verify the claims, which follow from the argument in \cite[p450]{G249}. We present the proof of \eqref{e332} in detail, while the proof of \eqref{e333}  is similar.

To verify \eqref{e332}, we notice that
$$
m(t,\vec{y}):=(t^2+4\pi^2|\mathcal{Q}(\vec{y})|^2)^{s/2}
$$
is homogeneous of degree $s$. Hence
$$
\lambda^{|\al|}\partial^\al m(\lambda t,\lambda \vec{y})=\lambda^s\partial^\al m(t,\vec{y})
$$
for $\lambda>0$, which, by taking $\lambda=|(t,y)|^{-1}$, implies that
$$
|\partial^\al m(t,\vec{y})|\lesssim \Big[\sup_{|(t,\vec{y})=1|}|\partial^\al m(t,\vec{y})|\Big] (t^2+|\vec{y}|^2)^{\tf{s-|\al|}2}.
$$
Then \eqref{e332} follows if we can show
\begin{equation}
\sup_{|(t,\vec{y})|=1}|\partial^\al m(t,\vec{y})|<\infty.
\end{equation}

It is easy to see that {$\partial^{\alpha} m$} is a finite sum of terms of the form
$$
(t^2+|\mathcal{Q}(\vec{y})|^2)^{\tf s2-k}P(t,\vec{y}),
$$
where $P(t,\vec{y})$ is a finite product of partial derivatives of $t^2+|\mathcal{Q}(\vec{y})|^2$ which is finite on the unit sphere.
To show that $(t^2+|\mathcal{Q}(\vec{y})|^2)^{\tf s2-k}<\infty$, we need to show that there exists a constant $C\in(0,\infty)$ such that
$$
1/C\le (t^2+|\mathcal{Q}(\vec{y})|^2)^{\tf s2-k}\le C.
$$
That $(t^2+|\mathcal{Q}(\vec{y})|^2)^{\tf s2-k}\le C$ is trivial. To prove the converse, we notice that on the unit sphere, at least one of $t^2$ and $|\vec{y}|^2$ is positive. Therefore, by that $\mathcal{Q}$ is non-degenerate $t^2+|\mathcal{Q}(\vec{y})|^2>0$ for any $|(t,\vec{y})|=1$. Using that the unit sphere is compact, we see that
$$
(t^2+|\mathcal{Q}(\vec{y})|^2)^{\tf s2-k}\ge 1/C
$$
on the unit sphere. This prove \eqref{e332}, and therefore Lemma~\ref{t323}.

\end{proof}

We complete the proof of Proposition~\ref{t321}.

\end{proof}

Let us summarize the necessary conditions on bilinear H\"ormander multipliers so far.

\begin{theorem}

Let $1<p_1, p_2<\infty$, $\tf1p=\tf1{p_1}+\tf1{p_2}$, $1<u<\infty$, and $0<s<\infty$.
Suppose that
\begin{equation}
\|T_m(f_1,f_2)\|_{L^p}\lesssim  \mathcal{A}_{s}^{u}[m] \prod_{i=1}^2\|f_i\|_{L^{p_i}},
\end{equation}
for all $m$ satisfying $\mathcal{A}_{s}^{u}[m]<\infty$.
Then it is necessary to have
\begin{equation}\label{e337}
\f sn> \max\Big\{\f{2}u,\f1p-\f1{u'}\Big\}
\end{equation} and

\begin{equation}\label{e336}
\tf sn\ge \max\{\tf12-\tf1{p'},0\}+\sum_{i=1}^2\max\{\tf12-\tf1{p_i},0\},
\end{equation}

\end{theorem}

\begin{proof}

\eqref{e336} is eaxctly \eqref{e335}. \eqref{e337} follows from Theorem~\ref{t327}.
The condition $\tf sn\ge \tf1p-\tf1{u'}-\sum_{i\in I}(\tf1{p_i}-\tf1u)$ obtained in Theorem~\ref{t326} is superfluous in the bilinear setting. Actually the case when $\#I\ge 1$ is discussed in Remark~\ref{t328} (iii), while the case when $\#I=0$ is improved in Proposition~\ref{t321}.

\end{proof}
{
Actually, we conjecture that the condition \eqref{e336} would be improved to
$$\tf sn> \max\{\tf12-\tf1{p'},0\}+\sum_{i=1}^2\max\{\tf12-\tf1{p_i},0\}$$
but we could not verify it.
}

\subsection{Proof of Theorem \ref{t329}}

{
We now complete the proof of Theorem \ref{t329}, by constructing multi-parameter examples based on the single-parameter version.}

We will discuss only
$$
\min\Big\{\f {s_1}n,\dots,\f{s_d}n\Big\}>\f{l}u,
$$
while other conditions can be obtained similarly from corresponding single-parameter results.

Without loss of generality, we may assume that $s_1\le \tf{ln}u$.
Now we define, for an appropriate large number $M$,
$$
m_1(\xi_{11},\xi_{21},\dots,\xi_{l1}):=\wh{\mathcal{H}_{(t,\gamma)}^{(ln,M)}}(\xi_{11}-\nu,\dots,\xi_{l1}-\nu)\prod_{i=1}^l\wh{\wt{\theta}}(\xi_{i1}-\nu)
$$
and
$$
m_k(\xi_{1k},\xi_{2k},\dots,\xi_{lk}):=\prod_{i=1}^l\wh{\wt{\theta}}(\xi_{ik}-\nu)
$$
for $2\le k\le d$.
The   multiplier we need is
$$
m(\xi_1,\dots,\xi_l):=\prod_{k=1}^d m_k(\xi_{1k},\xi_{2k},\dots,\xi_{lk}), \q { (\xi_1,\dots,\xi_l)\in (\bbr^{nd})^l.}
$$
By the argument in Proposition~\ref{t320},
\begin{equation}\label{e340}
\mathcal{A}_{s_1,\dots,s_d}^{u}[m]\lesssim \prod_{k=1}^d\|m_k\|_{L^u_{s_k}(\mathbb R^{ln})}
\lesssim \big\|\wh{\mathcal{H}_{(t-s_1,\gamma)}^{(ln,M)}}\big\|_{L^u(\bbr^{ln})}
\end{equation}
as $\|m_k\|_{L^u_{s_k}(\mathbb R^{ln})}\lesssim 1$ for $2\le k\le d$.

Let
$$
\wh{f_{j1}}(\xi_{j1}):=\epsilon^{n/p_j'}\wh{\theta}\Big(\frac{\xi_{j1}-\nu}{\epsilon}\Big)
$$
for $1\le j\le l$.

{Then by \eqref{limlowerbd'} in the proof of Proposition~\ref{t320},
\begin{equation}\label{ke1ca}
\liminf_{M\to \infty}\big\Vert T_{m_1}(f_{11},\dots,f_{l1})\big\Vert_{L^p(\bbrn)}\gtrsim \big\Vert \mathcal{H}_{(t,\gamma)}^{(ln)}\big\Vert_{L^1(\bbr^{ln})}
\end{equation}
}

For $2\le k\le d$, we define
$$
\wh{f_{jk}} (\xi_{jk}):=\epsilon^{n/p_j'}\wh{\theta}\big(\xi_{jk}-\nu\big)\qq \text{for }1\le j\le l.
$$
Then
\begin{equation}\label{kge2ca}
\big\Vert T_{m_k}(f_{1k},\dots, f_{lk}) \big\Vert_{ L^p(\bbr^{n})}=\||\theta|^l\|_{ L^p(\bbr^{n})}\approx 1.
\end{equation}
{Combining \eqref{ke1ca} and \eqref{kge2ca} together,}
if we take
$$\wh {f_j}(\xi_j):=\prod_{k=1}^d\wh{f_{jk}}(\xi_{jk}),$$
 then
{\begin{align}\label{htglnrlnest}
\liminf_{M\to \infty}\big\Vert T_{m }(f_{1},\dots, f_{l}) \big\Vert_{ L^p(\bbr^{dn})}&=\liminf_{M\to \infty}\Big\Vert \prod_{k=1}^dT_{m_k }(f_{1k},\dots, f_{lk}) \Big\Vert_{ L^p(\bbr^{dn})}\nonumber\\
&=\liminf_{M\to \infty}\prod_{k=1}^d\big\Vert T_{m_k }(f_{1k},\dots, f_{lk}) \big\Vert_{ L^p(\bbr^{n})}\nonumber\\
&\gtrsim \big\Vert \mathcal{H}_{(t,\gamma)}^{(ln)}\big\Vert_{L^1(\bbr^{ln})}.
\end{align}

So far, we have proved in \eqref{e340} and  \eqref{htglnrlnest} that
\begin{equation}
\mathcal{A}_{s_1,\dots,s_d}^{u}[m]\lesssim \big\|\wh{\mathcal{H}_{(t-s_1,\gamma)}^{(ln,M)}}\big\|_{L^u(\bbr^{ln})}, \text{uniformly in }~M
\end{equation}
and $$\liminf_{M\to \infty}\Vert T_m\Vert_{L^{p_1}\times\cdots\times L^{p_l}\to L^p}\gtrsim \big\Vert \mathcal{H}_{(t,\gamma)}^{(ln)}\big\Vert_{L^1(\bbr^{ln})}$$
and thus, by taking $t=ln$ and $\gamma=2$, we finally obtain the desired result.
}

\end{document}